\documentclass[11pt]{article}

\usepackage{latexsym}
\usepackage{amsmath}
\usepackage{amssymb}
\usepackage{amsbsy}
\usepackage{amscd}
\usepackage[all]{xy}
\usepackage[final]{graphicx}

\newtheorem{thm}{Theorem}[section]
\newtheorem{prop}[thm]{Proposition}
\newtheorem{defi}[thm]{Definition}
\newtheorem{lem}[thm]{Lemma}
\newtheorem{cor}[thm]{Corollary}

         \newcommand{\cN}{\mathcal N}
         \newcommand{\cO}{\mathcal O}
\newcommand{\cC}{\mathcal C}         
         \newcommand{\cQ}{\mathcal Q}
\newcommand{\cE}{\mathcal E}         
\newcommand{\cF}{\mathcal F}         
         \newcommand{\cT}{\mathcal T}
\newcommand{\cH}{\mathcal H}         \newcommand{\cU}{\mathcal U}
\newcommand{\cI}{\mathcal I}

\newcommand{\cM}{\mathcal M}

\newcommand{\bE}{\mathbb E}

\newcommand{\bP}{\mathbb P}

\newcommand{\bZ}{\mathbb Z}

\newcommand{\bo}{\hspace*{\fill} $\Box $ }

\newcommand{\barbeta}{\overline { \cM } _{0,0} \bigl( X, \beta \bigr)}

\newcommand{\barbieta}{\overline { \cM } _{bir} \bigl( X, \beta \bigr) }

\begin{document}

\pagestyle{headings}

\title{The irreducibility of the spaces of \\
rational curves on del Pezzo surfaces}
\author{Damiano Testa}

\maketitle

\begin{abstract}
Let $X$ be a del Pezzo surface of degree $d$, and assume that $X$ 
is general if $d=1$.  We prove that the spaces 
$\cM _{0,0} \bigl( X, \beta \bigr)$ are either empty or irreducible, if 
$(d, \beta ) \neq (1, -K_X)$.  When $(d, \beta ) = (1, -K_X)$ it is 
well known that $\cM _{0,0} \bigl( X, \beta \bigr)$ consists of twelve 
reduced points.
\end{abstract}

\subsection*{Introduction}
\addcontentsline{toc}{section}{Introduction}

Let $X$ be a del Pezzo surface.  If the degree of $X$ is one, assume that 
$X$ is general.  Let $\beta \in {\rm H}_2 \bigl( X, \bZ \bigr)$ be the class of 
a curve on $X$.  Denote by $R(\beta )$ the subscheme of 
the linear system $|\beta |$ consisting of the integral nodal curves 
of geometric genus zero.  The Kontsevich mapping space 
$\overline \cM _{0,0} \bigl( X, \beta \bigr)$ is a natural 
compactification of the space $R(\beta )$.  Some care is 
required, since the mapping spaces in general have more irreducible 
components than the corresponding spaces $R(\beta )$, arising from 
degenerate configurations of curves on the surface.  In fact it may 
happen that $R(\beta ) = \emptyset $, while 
$\overline \cM _{0,0} \bigl( X, \beta \bigr) \neq \emptyset $.

The mapping spaces 
parametrize the set of all (stable) maps to the surface $X$ 
from possibly reducible curves.  The domain curves of the maps 
in $\overline \cM _{0,0} \bigl( X, \beta \bigr)$ are connected 
and nodal, have all components isomorphic to $\bP^1$ and the 
components are attached in such a way that the resulting 
topological space is simply connected.  We refer to such a 
domain curve as a ``rational tree.''  Taking the image of a map 
yields a morphism $FC$ from (the semi-normalization of) 
$\overline \cM _{0,0} \bigl( X, \beta \bigr)$ to the closure in 
$|\beta |$ of $R(\beta )$ (see \cite{Ko} Section I.6).

Let $\barbieta $ be the closure of the subspace of $\barbeta $ 
consisting of morphisms $f : C \rightarrow X$, with $C \simeq \bP^1$ 
and $f$ birational onto its image.  In many interesting cases it is true 
and easy to check that the map $FC$ defined above is in fact 
birational, when restricted to $\barbieta $.

The main result of this paper is that the space $\barbieta $ is 
irreducible or empty, except in the case where $X$ has degree 
one and $\beta = - K_X$.

The idea of the proof is straightforward.  First, prove that in the 
boundary of all the irreducible components of $\barbieta $ there are 
special morphisms of a given type (called in what follows ``morphisms 
in standard form'').  Second, show that the locus of morphisms in standard 
form is connected and contained in the smooth locus of $\barbieta $.

From these two facts we conclude immediately that the smooth locus 
of $\barbieta $ is connected.  Since the smooth locus is dense, we deduce 
that $\barbieta $ is irreducible.

The methods used in the proof are of two different kinds.  First, 
there are general techniques, mainly Mori's Bend and Break Theorem, 
to break curves into components with low anticanonical degree.  
In the case where $X$ is the projective plane, this shows that we 
may specialize a morphism in $\barbieta $ so that its image is a 
union of lines.  Second, we need explicit geometric arguments to 
deal with the low degree cases.  Again, in the case of the 
projective plane, this step is used to bring the domain to a 
standard form (a chain of rational curves, rather than a general 
rational tree), while preserving the property that the image of 
the morphism consists of a union of lines.

To analyze the curves of low anticanonical degree on a del Pezzo 
surface, we need a detailed description of their divisor classes in 
${\rm Pic} (X)$.  In particular, we use the group of symmetries of 
the Picard lattice to reduce the number of cases to treat.  Section 
\ref{sporco} is devoted to this analysis.

Two technical deformation-theoretic tools prove useful.  The first 
is a description of the obstruction space of a stable map to a 
smooth surface in terms of combinatorial invariants of the map.  This 
is proved in Section \ref{conormale}.  The second is a lifting result 
that allows us, given a deformation 
of a component of a curve, to get a deformation of the whole curve.  
The statement is proved in Lemma \ref{ovvio} and the construction 
following it is the way in which we are typically going to use it.  This 
is specific to the surface case.  The lifting result allows us 
to deform a map with reducible domain by deforming only a few 
components at a time.  This is done systematically in Section 
\ref{rompisezione}.  We are therefore able to reduce the general 
problem to relatively few special cases, cf. Theorem \ref{passo}.  
The explicit computation of the obstruction spaces allows us to 
prove that in the deformations performed we never move to a 
different irreducible component of the moduli space.

The connectedness of the locus of morphisms in standard form is a 
consequence of some explicit computations, some of which are 
reformulations of classical geometric statements, such as the fact 
that the ramification locus of the projection from a general point 
on a smooth cubic surface in $\bP^3$ is a smooth plane quartic curve.  
This is the content of Section \ref{pochipochi}, but see also 
Section \ref{conichette}.

\tableofcontents

\section{Cohomology Groups and Obstruction Spaces}

\subsection{Rational Trees}

The purpose of this section is to prove some general results which 
are useful to compute the cohomology groups of coherent sheaves on 
rational trees.

\begin{defi}
A {\rm rational tree} $C$ is a connected, projective, nodal curve of 
arithmetic genus zero.  If $C$ is a rational tree, we call a component 
$E$ of $C$ an {\rm end} if $E$ contains at most one node of $C$.
\end{defi}

\begin{defi}
Given a connected projective nodal curve $C$, define the dual graph of $C$ 
to be the graph $\Gamma _C$ whose vertices are indexed by the components 
$C_i$ of $C$ and whose edges between the distinct vertices $[C_i]$ and 
$[C_j]$ are indexed by $\{ p\in C_i \cap C_j \}$.
\end{defi}

\noindent
{\it Remark}.  A connected projective nodal curve $C$ is a rational tree if 
and only if all its components are smooth rational curves and its dual graph 
$\Gamma _C$ is a tree (for a proof see \cite{De}).

\begin{lem} \label{basspl}
Let $C$ be a rational tree, and let $\nu : \tilde C \rightarrow C$ be the normalization of 
$C$ at the points $\{p_1, \ldots , p_r \} \subset Sing(C)$; denote by $\iota : \{p _1, 
\ldots , p_r \} \hookrightarrow C$ the inclusion morphism.  For any locally free sheaf 
$\cF$ of finite rank on $C$ we have the following short exact sequence of sheaves on $C$:
$$ 0 \longrightarrow \cF \longrightarrow \nu _* \nu ^* \cF \longrightarrow \iota _* \cF 
|_{\{ p_1, \ldots , p_r \}} \longrightarrow 0 $$
\end{lem}

\noindent
{\it Remark}.  From now on, we may sometimes denote the sheaf 
$\iota _* \cF |_{\{ p_1, \ldots , p_r \}}$ simply by $\oplus \cF _{p_i}$, 
and similarly for the pushforwards of sheaves on irreducible components 
of a curve.

\noindent
{\it Proof.}  Consider the sequence defining the sheaf $\cQ$:
\begin{equation} \label{defqu}
0 \longrightarrow \cO _C \longrightarrow \nu _* \cO _{\tilde C} \longrightarrow \cQ 
\longrightarrow 0
\end{equation}

Since $\nu $ is an isomorphism away from the inverse image of the $p_i$'s, it follows that 
$\cQ$ is supported at the union of the $p_i$'s.  We now want to prove that $\cQ _{p_i}$ is 
a skyscraper sheaf, i.e. that it has length one.  Since this is a local property, it is 
enough to check it when $C$ has a unique node.  In this case, $C$ is the nodal union of two 
smooth $\bP ^1$'s, and $\tilde C$ is their disjoint union.  Since the normalization map is 
finite, it is affine, and therefore ${\rm H}^j (\nu _* \cO _{\tilde C}, C) \simeq {\rm H}^j 
(\cO _{\tilde C}, \tilde C)$.  Therefore the long exact sequence defining $\cQ$ is given by
$$ 0 \longrightarrow k \longrightarrow k + k \longrightarrow \cQ \longrightarrow 0 $$
and we deduce that the length of $\cQ $ is 1.  Thus it follows in general that $\cQ = 
\oplus \cO _{p_i}$, the direct sum of the skyscraper sheaves of the nodes $p_1, \ldots , 
p_r$.

Let us now go back to the sequence (\ref{defqu}).  Since 
the sheaf $\cF$ is locally free, we may tensor the sequence 
by $\cF$, preserving exactness.  To identify the tensor 
product in the middle we use the projection formula:
$$ \nu _* \cO _{\tilde C} \otimes \cF \cong \nu _* (\cO _{\tilde C} \otimes \nu ^* \cF ) 
\cong \nu _* ( \nu ^* \cF ) $$
and we may therefore write the tensored sequence as
$$ 0 \longrightarrow \cF \longrightarrow \nu _* \nu ^* \cF \longrightarrow 
\mathop {\oplus } \limits _i \cF _{p_i} \longrightarrow 0 $$
thus proving the lemma.  \bo

Given a rational tree $C$ and a node $p \in C$, construct a new curve $C'$ as follows: 
consider the normalization of $\nu : \tilde C \rightarrow C$ of $C$ at the point $p$, 
and let $\{p_1, p_2 \} = \nu ^{-1} (p)$.  Attach to $\tilde C$ a smooth rational curve 
$E$ so that $\tilde C \cap E = \{p_1, p_2 \}$ and $C' := \tilde C \cup E$ is a nodal 
curve.  Clearly we have a morphism $\pi : C' \rightarrow C$, which is an isomorphism 
away from $E$ and contracts $E$ to the node $p$.  We call the morphism $\pi$ the 
contraction of $E$.  The curve $C'$ so obtained is called the ``total transform'' of 
$C$ at the node $p$ and $E$ the ``exceptional component.''

\begin{lem} \label{blodo}
Let $C$ be a rational tree, and let $\cF$ be a locally free coherent sheaf on $C$.  Let 
$\pi : C' \rightarrow C$ denote the total transform of $C$ at a node $p \in C$; then 
${\rm H}^1 (C, \cF) \cong {\rm H}^1 (C', \pi ^* \cF )$.
\end{lem}
{\it Proof.}  The result follows immediately from the Leray's spectral sequence associated 
to the map $\pi $ and the fact that $R ^i \pi _* (\cO_{C'}) = 0$, for $i > 0$.  \bo

\begin{lem} \label{pieces}
Let $C$ be a rational tree and let $\cF $ be a locally free sheaf on $C$.  Suppose 
that $C = C_1 \cup C_2$, where $C_1$ and $C_2$ are unions of components 
having no components in common.  Let $\{ p_1, \ldots , p_r \} = C_1 \cap C_2$ be 
the nodes of $C$ lying on $C_1$ and $C_2$.  If 
$h^1 (C_1, \cF | _{C_1} (-p_1 - \ldots -p_r)) = 0$, then 
${\rm H}^1 (C, \cF ) \cong {\rm H}^1 (C_2, \cF | _{C_2})$.
\end{lem}
{\it Proof.}  Simply consider the long exact sequence associated to the 
``component sequence''
$$ 0 \longrightarrow \cF | _{C_1} (-p_1 - \ldots -p_r) \longrightarrow \cF 
\longrightarrow \cF | _{C_2} \longrightarrow 0 $$
where the first map is extension by zero and the second map is restriction.  \bo

\begin{cor} \label{easyredu}
Let $C$ be a rational tree and let $R \subset C$ be a connected union of irreducible 
components of $C$.  Let $\cF$ be a locally free sheaf on $C$ such that the restriction 
of $\cF$ to each irreducible component of $C$ which is not in $R$ is generated by 
global sections.  Then $h^1 (C, \cF ) = h^1 (R, \cF |_R )$.
\end{cor}
{\it Proof.}  Proceed by induction on the number $\ell $ of 
irreducible components of $C$ not in $R$.  If 
$\ell = 0$, there is nothing to prove.  Suppose 
$\ell \geq 1$.  Let $C_1$ be an end of $C$ which is 
not an end of $R$, and let $p\in C_1$ be the node.  
The existence of such a component is easy to prove: 
since $R$ is a proper subcurve of $C$, there must be 
a node of $C$ where $R$ meets a component not in $R$.  
Removing this node disconnects $C$ into a connected 
component containing $R$ and a connected component $K$ 
disjoint from $R$.  Clearly an end $C_1$ of the component 
$K$ (different from the one meeting $R$ if there is more 
than one end) is then also an end of $C$ not contained in 
$R$.  Since $C_1$ is a smooth rational curve, 
$\cF| _{C_1} \simeq \bigoplus \cO (a_j)$, with $a_j \geq 0$, 
thanks to the fact that $\cF| _{C_1}$ is globally generated.  
In particular $h^1 (C_1, \cF (-p)) = 1$, and it is clear 
that we can now apply Lemma \ref{pieces} to remove the 
component $C_1$ without changing $h^1$ and conclude using 
induction.  \bo

The last lemma of this section is an explicit computation of the 
cohomology of a locally free sheaf on a curve which will be extremely 
useful in the later sections.

\begin{lem} \label{tretre}
Let $C$ be a rational tree and $f: C \rightarrow S$ a morphism to a smooth surface.  Let $p 
\in C$ be a node, denote by $C_a$ and $C_b$ the two irreducible components of $C$ meeting 
at $p$.  Let $\nu : \tilde C \rightarrow C$ be the normalization of $C$ at $p$ and let 
$\tilde f = f \circ \nu $.  Suppose that:
\begin{enumerate}
\item the valences of the vertices $C_a$ and $C_b$ in the dual graph of $C$ are at most 3, and
\item the map $f_* : \cT _{C_a, p} + \cT _{C_b, p} \longrightarrow \cT _{S, f(p)} $ is 
surjective.
\end{enumerate}

Then 
$$ {\rm H}^1 (C, f^* \cT _S ) \cong {\rm H}^1 (\tilde C, \tilde f ^* \cT _S )$$
\end{lem}
{\it Proof.}  Consider the sequence on $C$
$$ 0 \longrightarrow f^* \cT _S \longrightarrow \tilde f ^* \cT _S \stackrel {\varepsilon} 
{\longrightarrow } \cT _{S, f(p)} \longrightarrow 0 $$

Because $\cT _{S, f(p)}$ is supported in dimension 0, ${\rm H}^1 (C, \cT _{S, f(p)}) = 0$, and 
it is enough to prove that the sequence is exact on global sections.  Let $\{ p, q_a, r_a \}$ 
contain all the nodes of $C$ on $C_a$ and let $\{p, q_b, r_b \}$ contain the nodes on $C_b$.  
Consider now the following diagram:
$$\xymatrix {0 \ar[d] & 0 \ar[d] \\
\cT _{C_a} (-p -q_a -r_a) \oplus \cT _{C_b} (-p -q_b -r_b) \ar[d] & f^* \cT _S \ar[d] \\
\cT _{C_a} (-q_a -r_a) \oplus \cT _{C_b} (-q_b -r_b) \ar[d] ^{\alpha } \ar[r] & 
\tilde f ^* \cT _S \ar[d] ^{\varepsilon } \\
\cT _{C_a, p} \oplus \cT _{C_b, p} \ar[d] \ar[r] ^{f_*} & \cT _{S, f(p)} \ar[d] \\
0 & 0}$$
where the unlabeled horizontal map is extension by zero.  Since $C_a$ and $C_b$ are rational 
curves, their tangent bundles have degree 2 and $\alpha $ is surjective on global sections; 
$f_*$ is surjective by assumption.  It follows that $\varepsilon $ is also surjective on 
global sections.  \bo

\noindent
{\it Remark}.  The second condition in the lemma is certainly satisfied if $f|_{C_a}$ and 
$f|_{C_b}$ are birational and the intersection of $f(C_a)$ and $f(C_b)$ is transverse at 
$f(p)$.

\subsection{The Conormal Sheaf} \label{conormale}

Let $f:C \rightarrow X$ be a morphism from a connected, projective, 
at worst nodal curve $C$ to a smooth projective variety $X$.

\begin{defi}
The morphism $f : C \rightarrow X$ is called a stable map if $C$ is 
a connected, projective, at worst nodal curve and every contracted 
component of geometric genus zero contains at least three singular points 
of $C$ and every contracted component of geometric genus one contains 
at least one.
\end{defi}

We are interested in computing the obstruction space to deforming the 
stable map $f : C \rightarrow X$.  Let $f^* \Omega ^1 _X \rightarrow
\Omega ^1 _C$ be the natural complex of sheaves associated with the 
differential of $f$ and where the sheaf $f^* \Omega ^1 _X$ is in 
degree -1 and the sheaf $\Omega ^1 _C$ is in degree 0.  We know that 
the stability condition is equivalent to the vanishing of the group 
${\rm Hom} \bigl( f^* \Omega ^1 _X \rightarrow
\Omega ^1 _C , \cO _C \bigr)$.  The tangent space to 
$\overline \cM _{0,0} \bigl( X, \beta \bigr)$ at $f$ is given by 
the hypercohomology group 
${\rm \bE xt } ^1 \bigl( f^* \Omega ^1 _X \rightarrow
\Omega ^1 _C , \cO _C \bigr)$.  The obstruction space 
is a quotient of the hypercohomology group 
${\rm \bE xt } ^2 \bigl( f^* \Omega ^1 _X \rightarrow \Omega ^1 _C , 
\cO _C \bigr)$.  Denote by $L _f ^\bullet $ the complex 
$f^* \Omega ^1 _X \rightarrow \Omega ^1 _C$, where the first sheaf is 
in degree -1 and the second one is in degree 0.  Our strategy to 
compute these groups is to use the short exact sequence of complexes 
of sheaves:
$$\xymatrix@R=17pt { \left( 0 \rightarrow 0 \right) \ar[d] \\
\left( 0 \rightarrow \Omega ^1 _C \right) \ar[d] \\
\left( f^* \Omega ^1 _X \rightarrow \Omega ^1 _C \right) \ar[d] \\
\left( f^* \Omega ^1 _X \rightarrow 0 \ar[d] \right) \\
\left( 0 \rightarrow 0 \right) }$$

Applying the functor ${\rm Hom } ( - , \cO _C)$ and using the long 
exact hypercohomology sequence we obtain:
$$\xymatrix @R=3pt {& 0 \ar[r] & 
{\rm Hom } \bigl( L _f ^\bullet , \cO _C \bigr) \ar[r] & 
{\rm Hom } (0 \rightarrow \Omega ^1 _C , \cO _C) \\ \ar[r] &
{\rm \bE xt } ^1 (f^* \Omega ^1 _X \rightarrow 0 , \cO _C) \ar[r] & 
{\rm \bE xt } ^1 \bigl( L _f ^\bullet , \cO _C \bigr) \ar[r] & 
{\rm \bE xt } ^1 (0 \rightarrow \Omega ^1 _C , \cO _C) \\ \ar[r] &
{\rm \bE xt } ^2 (f^* \Omega ^1 _X \rightarrow 0 , \cO _C) \ar[r] & 
{\rm \bE xt } ^2 \bigl( L _f ^\bullet , \cO _C \bigr) \ar[r] & 
{\rm \bE xt } ^2 (0 \rightarrow \Omega ^1 _C , \cO _C)}$$

We can now rewrite many of these terms.  First of all, the stability 
condition is equivalent to 
${\rm Hom } \bigl( L _f ^\bullet , \cO _C \bigr) = 0$.  
Also, remembering the fact that all the complexes are concentrated in 
degrees $-1$ and 0, and using the fact that $f^* \Omega _X ^1$ is 
locally free, that its dual is $f^* \cT _X$ and the 
isomorphisms ${\rm Ext} ^i (f^* \Omega _X ^1 , \cO _C) \simeq 
{\rm H} ^i (C, \cT _X)$ we obtain the sequence
\begin{equation} \label{preserra}
\xymatrix@C=20pt @R=3pt {0 \ar[r] & {\rm Hom } (\Omega ^1 _C , \cO _C) \ar[r] &
{\rm H}^0 (C, f^* \cT _X) \ar[r] & 
{\rm \bE xt } ^1 \bigl( L _f ^\bullet , \cO _C \bigr) \ar[r] &\\
\ar[r] & {\rm Ext } ^1 (\Omega ^1 _C , \cO _C) \ar[r] &
{\rm H}^1 (C, f^* \cT _X) \ar[r] & 
{\rm \bE xt } ^2 \bigl( L _f ^\bullet , \cO _C \bigr) \ar[r] & 0}
\end{equation}

In particular we see that if ${\rm H}^1 (C, f^* \cT _X)=0$, then the 
obstruction group ${\rm \bE xt } ^2 \bigl( L _f ^\bullet , \cO _C \bigr)$ 
vanishes as well, i.e. the map is unobstructed, the space of stable maps has 
the expected dimension at $f$ and the point represented by 
$f$ is smooth (for the stack).

If we consider the dual sequence of (\ref{preserra}) and use Serre duality 
we obtain the sequence 
$$ \xymatrix @R=3pt {0 \ar[r] & \bigl( {\rm \bE xt } ^2 \bigr) ^{\vee } \ar[r] & 
{\rm H}^0 (C, f^* \Omega ^1 _X \otimes \omega _C ) \ar[r] ^\alpha & 
{\rm H} ^0 (C , \Omega ^1 _C \otimes \omega _C ) \ar[r] & \\
\ar[r] & \bigl( {\rm \bE xt } ^1 \bigr) ^{\vee } \ar[r] & 
{\rm H} ^1 (C, f^* \Omega ^1 _X \otimes \omega _C ) \ar[r] & 
{\rm H} ^1 (C, \Omega ^1 _C \otimes \omega _C) \ar[r] & 0} $$

It is easy to convince oneself that the morphism $\alpha $ is the 
morphism induced by the differential map 
$df : f^* \Omega ^1 _X \longrightarrow \Omega _C$, by tensoring with 
the dualizing sheaf and taking global sections.  Associated to $f$ 
we may define the sheaves $\cC _f$ and $\cQ _f$ on $C$, by requiring 
the following sequence to be exact:
\begin{equation} \label{conor}
\xymatrix{ 0 \ar[r] & \cC_f \ar[r] & f^* \Omega ^1 _X \ar[r]^{df} & 
\Omega ^1 _C \ar[r] & \cQ _f \ar[r] & 0 }
\end{equation}
Since the dualizing sheaf $\omega _C$ is locally free, tensoring by 
$\omega _C$ is exact and taking global sections is left exact.  From 
these remarks we deduce that 
$$ {\rm H}^0 \bigl( C, \cC_f \otimes \omega _C \bigr) \simeq 
{\rm \bE xt } ^2 \bigl( L _f ^\bullet , \cO _C \bigr) ^\vee $$
and we conclude that in order to compute the obstruction space of $f$, 
it is enough to compute the global sections of the sheaf 
$\cC_f \otimes \omega _C$.

\begin{defi}
The sheaf $\cC _f$ defined in (\ref{conor}) is the {\rm conormal 
sheaf of $f$}.
\end{defi}

We drop the subscript $f$, when the morphism is clear from 
the context.

\begin{defi}
A sheaf $\cF$ on a scheme of pure dimension one is {\rm pure} 
if the support of every non-zero section has pure dimension one.
\end{defi}

It is clear that a locally free sheaf is pure.  In fact, any subsheaf 
of a locally free sheaf is pure, and more generally any subsheaf of a 
pure sheaf is pure.  In particular, the sheaves $\cC _f$ defined in 
(\ref{conor}) are pure.

\begin{defi} 
A point $p \in C$ is called a {\rm break for the morphism $f$} (or 
simply a {\rm break}), if the sheaf $\cC _f$ is not locally free at 
$p$.  We say that the morphism $f$ has no breaks if the sheaf $\cC_f$ 
is locally free.
\end{defi}

It is clear from the definition that a smooth point of $C$ is never a break.

\begin{defi}
Suppose the morphism $f$ is finite.  A point $p \in C$ is called a {\rm ramification 
point}, if it belongs to the support of the sheaf $\cQ _f$.  We call {\rm 
ramification divisor of $f$} the (Weil) divisor whose multiplicity at $p\in C$ is the 
length of $\cQ$ at $p$.
\end{defi}

Let $f_1 : C_1 \rightarrow X$ and $f_2 : C_2 \rightarrow X$ be non-constant morphisms from 
two smooth curves to a smooth surface.  Suppose $p_1 \in C_1$ and $p_2 \in C_2$ are points 
such that $f_1 (p_1) = f_2 (p_2) = q$, let $u$ and $v$ be local coordinates on $X$ near 
$q$ and let $x_1$ and $x_2$ be local parameters for $C_1$ and $C_2$ near $p_1$ and $p_2$ 
respectively.  Since $f_1$ and $f_2$ are not constant, there exist integers $k_1$ and 
$k_2$ such that
$$ f_1 ^* : \left\{ \begin{tabular}{c@{ $\longmapsto $ }c}
$u$&$x_1 ^{k_1} U_1 (x_1)$ \vphantom{$\frac{\frac{1}{1}}{\frac{1}{1}}$} \\
$v$&$x_1 ^{k_1} V_1 (x_1)$ \vphantom{$\frac{\frac{1}{1}}{\frac{1}{1}}$} 
\end{tabular} \right. \hspace{20pt}
f_2 ^* : \left\{ \begin{tabular}{c@{ $\longmapsto $ }c}
$u$&$x_2 ^{k_2} U_2 (x_2)$ \vphantom{$\frac{\frac{1}{1}}{\frac{1}{1}}$} \\
$v$&$x_2 ^{k_2} V_2 (x_2)$ \vphantom{$\frac{\frac{1}{1}}{\frac{1}{1}}$} 
\end{tabular} \right. $$
and $\bigl( U_1 (0), V_1 (0) \bigr) , \bigl( U_2 (0), V_2 (0) \bigr) \neq (0,0)$.  
We call a {\it tangent vector to $C_i$ at $p_i$} any non-zero vector in $\cT _q X$ 
proportional to $\bigl( U_i (0), V_i (0) \bigr)$, and {\it tangent direction to 
$C_i$ at $p_i$} the point in $\bP \left( \cT _q X \right)$ determined by a tangent 
vector to $C_i$ at $p_i$.  Geometrically, we may easily associate to each smooth 
point of $f_i (C_i)$ a tangent vector in the same way we did above, and then the 
tangent direction at any point is simply the limiting position of the tangent 
directions at the smooth points.

We say that {\it $C_1$ and $C_2$ are transverse at the point $q = f_i (p_i) \in X$} 
if their respective tangent directions at $p_1$ and $p_2$ are distinct and we will 
say that {\it $C_1$ and $C_2$ are not transverse at the point $q \in X$} if the 
tangent directions coincide.

Finally, we say that the morphism $f_i$ is ramified at $p_i$ on $C_i$ if $k_i > 1$ 
and we say it is unramified at $C_i$ if $k_i = 1$.

\begin{lem}
Let $f_i : C_i \rightarrow X$, $i \in \{ 1, 2 \}$ be two non-constant morphisms from two 
smooth curves to a smooth surface $X$ and let $p_1 \in C_1$ and $p_2 \in C_2$ be points 
such that $f_1 (p_1) = f_2 (p_2) = q$.  Denote by $\tilde f _1$ and $\tilde f _2$ the 
morphisms induced by $f_1$ and $f_2$ from each curve to the blow-up of $X$ at $q$, and 
assume $\tilde f _1 (p_1) = \tilde f _2 (p_2) = \tilde q$.  Then the following conditions 
are equivalent:
\begin{enumerate}
\item $\tilde f _1$ and $\tilde f _2$ are unramified at $\tilde q$ and $C_1$ and $C_2$ 
are transverse at $\tilde q$;
\item after possibly 
renumbering the curves $C_1$ and $C_2$, there are coordinates $u,v$ on $X$ near $q$ and 
$x_i$ on $C_i$ near $p_i$ such that
\begin{equation} \label{bln2}
\begin{array} {ll}
f_1 ^* : 
\left\{ \begin{tabular}{c@{ $\longmapsto$ }c}
$u$ & $x_1 U_1 (x_1)$ \\[5pt]
$v$ & $x_1 ^3 V_1 (x_1)$ \\[5pt]
\end{tabular} \right. & U_1 (0) \neq 0 \vspace{10pt} \\ 
f_2 ^* : \left\{ 
\begin{tabular}{c@{ $\longmapsto $ }c}
$u$&$x_2 ^{m} U_2 (x_2)$ \vphantom{$\frac{\frac{1}{1}}{\frac{1}{1}}$} \\
$v$&$x_2 ^{m+1} V_2 (x_2)$ \vphantom{$\frac{\frac{1}{1}}{\frac{1}{1}}$} 
\end{tabular} \right. & 
\begin{tabular} {l} $U_2 (0)$, $V_2 (0) \neq 0 $ \\ $m \geq 1$ \end{tabular}
\end{array}
\end{equation}
\end{enumerate}
\end{lem}
{\it Proof.}  Suppose we are given coordinates so that the $f_i$'s are given by (\ref{bln2}).  
Let $b : \tilde X \rightarrow X$ be the blow-up morphism.  Let $\tilde u := b^* u$, and 
note that near the point $\tilde q$ the function $\tilde u$ is a local equation for 
$E := b^{-1} (q)$, since the tangent vector to the curve locally defined by the vanishing of 
$u$ is $(0,1)$, while a tangent vector to the curve $C_1$ at $q$ is $(1,0)$.  It follows that 
we may write $b^* v = \tilde u \cdot \tilde v$, and 
$\tilde u , \tilde v$ is a local system of parameters on $\tilde X$ at $\tilde q$ 
such that $b$ and its rational inverse $b ^{-1}$ are given by:
\begin{equation} \label{bloco}
b ^* : \left\{ \begin{tabular}{c@{ $\longmapsto $ }c}
$u$& $\tilde u$ \\$v$& $\tilde u \tilde v$ \end{tabular} \right. \hspace{20pt} 
(b ^*) ^{-1} : \left\{ \begin{tabular}{c@{ $\longmapsto $ }c}
$\tilde u$& $u$ \\$\tilde v$& $v / u$ \end{tabular} \right.
\end{equation}

Thus the morphisms $\tilde f _i : C_i \rightarrow \tilde X$ are given by
\begin{eqnarray*}
\tilde f_1 ^* : & \left\{ \begin{tabular}{c@{ $\longmapsto $ }c}
$\tilde u$&$x_1 U_1 (x_1)$ \vphantom{$\frac{\frac{1}{1}}{\frac{1}{1}}$} \\
$\tilde v$&$x_1 ^2 \frac {V_1 (x_1)} {U_1 (x_1)}$
\end{tabular} \right. & U_1 (0) \neq 0 \\[5pt]
\tilde f_2 ^* : & \left\{ \begin{tabular}{c@{ $\longmapsto $ }c}
$\tilde u$&$x_2 ^{m} U_2 (x_2)$ \vphantom{$\frac{\frac{1}{1}}{\frac{1}{1}}$} \\
$\tilde v$&$x_2 \frac {V_2 (x_2)} {U_2 (x_2)}$ 
\end{tabular} \right. & U_2 (0), V_2 (0) \neq 0
\end{eqnarray*}

Clearly these maps are unramified at $x_i = 0$ and since $(1,0)$ 
and $(\star ,1)$ are tangent vectors at $\tilde q$  to $C_1$ and $C_2$ respectively, 
the maps are also transverse at $\tilde q$.  This simple computation proves the first 
half of the lemma.

Suppose conversely that in the blow-up $\tilde X$ of $X$ at $q$, the curves $C_1$ and 
$C_2$ meet transversely at the point $\tilde q = \tilde f_i (p_i) \in \tilde X$.  Fix 
coordinates $x_1$ on $C_1$ at $p_1$ and $x_2$ on $C_2$ at $p_2$, and choose coordinates 
$u,v$ near $q$ and $\tilde u , \tilde v$ near $\tilde q$ such that (\ref{bloco}) are 
the equations of the blow-up morphism.  We have
$$ \tilde f_1 ^* : \left\{ \begin{tabular}{c@{ $\longmapsto $ }c}
$\tilde u $&$x_1 ^{k_1} U _1 (x_1)$ \vphantom{$\frac{\frac{1}{1}}{\frac{1}{1}}$} \\
$\tilde v $&$x_1 ^{k_1} V _1 (x_1)$
\end{tabular} \right. \hspace{20pt} 
\tilde f_2 ^* : \left\{ \begin{tabular}{c@{ $\longmapsto $ }c}
$\tilde u $&$x_2 ^{k_2} U _2 (x_2)$ \vphantom{$\frac{\frac{1}{1}}{\frac{1}{1}}$} \\
$\tilde v $&$x_2 ^{k_2} V _2 (x_2)$
\end{tabular} \right. $$
with $\bigl( U _1 (0), V _1 (0) \bigr) , \bigl( U _2 (0), V _2 (0) \bigr)$ 
linearly independent.  By changing $v$ to $v - \frac{V_1 (0)} {U_1 (0)} u$ 
and $\tilde v$ to $\tilde v - \frac{V_1 (0)} {U_1 (0)} \tilde u$, we may 
assume that $V_1 (0) = 0 $, while preserving the equations of $b$.  With 
these assumptions, $(1,0)$ and $(\star , 1)$ are tangent vectors at 
$\tilde q$ to $C_1$ and $C_2$ respectively.  Moreover, since $\tilde f _i$ 
is not ramified at $p_i$, necessarily $k_i = 1$.  We have therefore
\begin{eqnarray*}
\tilde f_1 ^* : \left\{ \begin{tabular}{c@{ $\longmapsto $ }c}
$\tilde u $&$x_1 U _1 (x_1)$ \vphantom{$\frac{\frac{1}{1}}{\frac{1}{1}}$} \\
$\tilde v $&$x_1 ^2 \overline V _1 (x_1)$
\end{tabular} \right. & &
f_1 ^* = \tilde f_1 ^* \circ b ^*: \left\{ \begin{tabular}{c@{ $\longmapsto $ }c}
$u $&$x_1 U _1 (x_1)$ \vphantom{$\frac{\frac{1}{1}}{\frac{1}{1}}$} \\
$v $&$x_1 ^3 U _1 (x_1) \overline V _1 (x_1)$
\end{tabular} \right. \\  &  \Longrightarrow & \\
\tilde f_2 ^* : \left\{ \begin{tabular}{c@{ $\longmapsto $ }c}
$\tilde u $&$x_2 U _2 (x_2)$ \vphantom{$\frac{\frac{1}{1}}{\frac{1}{1}}$} \\
$\tilde v $&$x_2 V _2 (x_2)$
\end{tabular} \right. & &
f_2 ^* = \tilde f_2 ^* \circ b ^* : \left\{ \begin{tabular}{c@{ $\longmapsto $ }c}
$u $&$x_2 U _2 (x_2)$ \vphantom{$\frac{\frac{1}{1}}{\frac{1}{1}}$} \\
$v $&$x_2 ^2 U _2 (x_2) V _2 (x_2)$
\end{tabular} \right.
\end{eqnarray*}
where $U_1 (0), V_2 (0) \neq 0$.

In order to conclude we still need to show that $U_2 (x_2)$ is not identically zero, 
but this is clear, since otherwise the morphism $f _2$ would be constant (i.e. the 
morphism $\tilde 
f_2$ would map $C_2$ to the exceptional divisor $E$).  \bo

\begin{defi} \label{semtang}
In the situation described by the previous lemma, the two curves $C_1$ and $C_2$ 
are {\rm simply tangent at $q$}.
\end{defi}

We will see later (Lemma \ref{conota}) that being simply tangent is closely related to the 
local structure of the conormal sheaf.

\begin{lem} \label{conotra}
Suppose that $X$ is a smooth surface and let $f: C \rightarrow X$ be a morphism from a 
curve $C$ consisting of two irreducible components $C_1$ and $C_2$, meeting in a node $p$.  
Denote by $f_i$ the restriction of $f$ to $C_i$ and by $p_i \in C_i$ the point $p \in C$, 
and suppose that $f$ does not contract any component of $C$ and that $C_1$ and $C_2$ meet 
transversely at $f(p)$.  Then there are the following cases:
\begin{enumerate}
\item \label{tuu} Both maps $f_1$ and $f_2$ are unramified at $p$.

Then $\cC _f$ is locally free and the following sequence is exact
$$\xymatrix{ 0 \ar[r] & \cC _f \ar[r] & \cC _{f_1} (-p) \oplus \cC _{f_2} (-p) \ar[r] & 
\cC _{f,p} \ar[r] & 0 } $$
\item \label{tur} $f_i$ is unramified at $p$ on $C_i$ and $f_{3-i}$ 
is ramified at $p$ on $C_{3-i}$ ($i \in \{ 1, 2 \}$)

Then $\cC _f$ is not locally free (i.e. $p$ is a break point) and
$$ \cC _f \cong \cC _{f_i} (-p) \oplus \cC _{f_{3-i}} (-2p) $$
\item \label{trr} Both maps $f_1$ and $f_2$ are ramified at $p$.

Then $\cC _f$ is not locally free and
$$ \cC _f \cong \cC _{f_1} (-p) \oplus \cC _{f_2} (-p) $$
\end{enumerate}
\end{lem}
{\it Proof.}  We can write
$$ f ^* : \left\{ \begin{tabular}{c@{ $\longmapsto $ }c}
$u$&$x^{k_1} U_1 (x) + y^{k_2} U_2 (y)$ \vphantom{$\frac{\frac{1}{1}}{\frac{1}{1}}$} \\
$v$&$x^{l_1} V_1 (x) + y^{l_2} V_2 (y)$ \vphantom{$\frac{\frac{1}{1}}{\frac{1}{1}}$} 
\end{tabular} \right. $$
where $l_1 > k_1$, $k_2 > l_2$ and $U_1(0), V_2(0) \neq 0$.  We thus have
$$ \xymatrix @R=3pt @C=-30.5pt { 
\cO _{C,p} \cdot du + \cO _{C,p} \cdot dv \ar[r]^{df \hspace{30pt}} & 
\raisebox {3pt} {$\Bigl( \cO _{C,p} \cdot dx + \cO _{C,p} \cdot dy \Bigr) $} / \raisebox 
{-3pt} {$\bigl( ydx + xdy \bigr)$} \\
du \ar[r] & x^{k_1-1} \Bigl( k_1 U_1(x) + x U_1 '(x) \Bigr) dx + y^{k_2-1} 
\Bigl( k_2 U_2(y) + y U_2 '(y) \Bigr) dy \\
dv \ar[r] & x^{l_1-1} \Bigl( l_1 V_1(x) + x V_1 '(x) \Bigr) dx + y^{l_2-1} 
\Bigl( l_2 V_2(y) + y V_2 '(y) \Bigr) dy } $$

In order to simplify this expression, let us define $\alpha _1 $ to be the invertible 
function $k_1 U_1(x) + x U_1 '(x)$ and $\alpha _2$ to be the invertible function 
$l_2 V_2(y) + y V_2 '(y)$.  Choosing $\frac {du } {\alpha _1}$ and $\frac {dv} {\alpha _2}$ 
as a basis for the $\cO _{C,p} -$module $f^* \Omega ^1 _{X,p}$ we may write
$$ \xymatrix{ 
\frac{du} {\alpha _1} \ar[r] & x^{k_1-1} dx + y^{k_2-1} \varphi (y) dy \\
\frac{dv} {\alpha _2} \ar[r] & x^{l_1-1} \psi (x) dx + y^{l_2-1} dy } $$

Note that
\begin{eqnarray*}
y^{k_2-1} \varphi (y) & = & \frac {y^{k_2-1}} {k_1 U_1(0)} 
\Bigl( k_2 U_2(y) + y U_2 ' (y) \Bigr) \\
x^{l_1-1} \psi (x) & = & \frac {x^{l_1-1}} {l_2 V_2(0)} 
\Bigl( l_1 V_1(x) + x V_1 ' (x) \Bigr)
\end{eqnarray*}

The elements of the kernel of $df$ are determined by the condition
$$ f_1 (x,y) \frac{du} {\alpha _1} + f_2 (x,y) \frac{dv} {\alpha _2} \longmapsto 
r(x,y) \bigl( ydx + xdy \bigr) $$
which translates to
\begin{eqnarray} \nonumber
x^{k_1-1} \Bigl(f_1 (x,y) + x^{l_1-k_1} f_2 (x,y) \psi (x) \Bigr) & = & y r(x,y) = 
y r(0,y) \\ \label{equara} \\ \nonumber
y^{l_2-1} \Bigl(y^{k_2-l_2} f_1 (x,y) \varphi (y) + f_2 (x,y) \Bigr) & = & x r(x,y) = 
x r(x,0)
\end{eqnarray}

We are now going to split the three cases.

\noindent
{\bf Case \ref{tuu}}.  In this case $k_1 = l_2 = 1$, and equation (\ref{equara}) 
becomes
\begin{eqnarray*}
f_1 (x,y) + x^{l_1-1} f_2 (x,y) \psi (x) & = & y r(x,y) = y r(0,y) 
\\ \\
y^{k_2-1} f_1 (x,y) \varphi (y) + f_2 (x,y) & = & x r(x,y) = x r(x,0)
\end{eqnarray*}

This clearly implies that neither $f_1$ nor $f_2$ have constant term and hence we may 
write $f_1(x,y) = x g_1(x) + y h_1 (y)$ and $f_2(x,y) = x g_2 (x) + y h_2(y)$ and we have
\begin{eqnarray*}
y h_1(y) + x \Bigl( g_1 (x) + x^{l_1-1} g_2 (x) \psi (x) \Bigr) & = & y r(x,y) \\[7pt]
x g_2(x) + y \Bigl( y^{k_2-1} h_1 (y) \varphi (y) + h_2 (y) \Bigr) & = & x r(x,y)
\end{eqnarray*}

Therefore
\begin{eqnarray*}
x g_1 (x) & = & - x^{l_1} g_2 (x) \psi (x) \\
y h_2(y) & = & - y^{k_2} h_1 (y) \varphi (y) \\
y h_1(y) & = & y r(x,y) \\
x g_2(x) & = & x r(x,y)
\end{eqnarray*}
and near $p$ all elements of the kernel of $df$ are multiples of
$$ \kappa := \Bigl( - x^{l_1} \psi (x) + y \Bigr) \frac {du} {\alpha _1} + 
\Bigl( x - y^{k_2} \varphi (y) \Bigr) \frac {dv} {\alpha _2} $$

It is very easy to check that $\kappa $ is also in the kernel of $df$.  This in 
particular implies that $\cC _f$ is locally free near $p$.  The restriction of $\kappa $ 
to $C_1$ (which is defined near $p$ by $y=0$) is
$$ \frac { - x^{l_1} \psi (x)} {\alpha _1} du + \frac {x} {l_2 V_2 (0)} dv = \frac {x} 
{l_2 V_2 (0)} \left( \frac { - x^{l_1-1} \Bigl( l_1 V_1 (x) + x V_1 ' (x) \Bigr)} 
{U_1 (x) + x U_1 ' (x)} du +  dv \right) $$

On the other hand, the restriction of $f$ to $C_1$ is given by
$$ f_1 ^* : \left\{ \begin{tabular}{c@{ $\longmapsto $ }c}
$u$&$x U_1 (x)$ \vphantom{$\frac{\frac{1}{1}}{\frac{1}{1}}$} \\
$v$&$x^{l_1} V_1 (x)$ \vphantom{$\frac{\frac{1}{1}}{\frac{1}{1}}$} 
\end{tabular} \right. $$
and the kernel of $df_1$ is clearly generated by
$$ \frac { - x^{l_1-1} \Bigl( l_1 V_1 (x) + x V_1 ' (x) \Bigr)} 
{U_1 (x) + x U_1 ' (x)} du + dv $$

Thus $\cC _{f_1} (-p)$ is generated near $p$ by the same generator of $\cC _f |_{C_1}$.  
Similarly, $\cC _{f_2} (-p)$ and $\cC _f |_{C_2}$ have the same generator near $p$.  Hence 
\ref{tuu} follows.

\noindent
{\bf Case \ref{tur}}.  Let us go back to equation (\ref{equara}) and substitute $k_1=1$ 
and $l_2 \geq 2$:
\begin{eqnarray*}
\vphantom{\frac{1}{1}} f_1 (x,y) + x^{l_1-1} f_2 (x,y) \psi (x) & = & y r(x,y) = y r(0,y) \\[7pt]
y^{l_2-1} \Bigl( y^{k_2-l_2} f_1 (x,y) \varphi (y) + f_2 (x,y) \Bigr) & = & x r(x,y) = x r(x,0)
\end{eqnarray*}

This implies that $r(x,0) = 0$, i.e. $r(x,y) = yr(y)$.  Thus we have 
$f_1 (x,y) = x g_1(x) + y^2 r (y)$, and finally 
$f_2 (x,y) = x g_2(x) - y^{k_2- l_2 +2} \varphi (y) r (y)$.  
Substituting back, we find
$$x g_1 (x) + x^{l_1} g_2 (x) \psi (x) = 0 $$

Therefore $x g_1 (x) = - x^{l_1} g_2 (x) \psi (x) $ and near $p$ the kernel of $df$ is 
generated by
\begin{eqnarray*}
x \left( - x^{l_1-1} \psi (x) \frac{du} {\alpha _1} + \frac{dv} {\alpha _2} \right) 
& \hspace{20pt} {\rm and} \hspace{20pt}
& y^2 \left( \frac{du} {\alpha _1} - y^{k_2- l_2} \varphi (y) \frac{dv} {\alpha _2} \right)
\end{eqnarray*}
(as before, it is very easy to check that these elements lie indeed in the kernel of $df$).

We thus see that $\cC _f$ is not locally free near $p$; since clearly the terms in brackets in 
the previous expression are local generators for $\cC _{f_1}$ and $\cC _{f_2}$ respectively near 
$p$, we deduce that $\cC _f \cong \cC _{f_1} (-p) \oplus \cC _{f_2} (-2p)$.  Thus \ref{tur} 
follows.

\noindent
{\bf Case \ref{trr}}.  Once more we refer to (\ref{equara}), now with $k_1, l_2 \geq 2$.  
In this case, the equations imply that $r(x,y) = 0$, and thus $f_1 (x,y) = - x^{l_1-k_1} 
f_2 (x,y) \psi (x) + y h_1 (y)$.  Substituting back in (\ref{equara}), we find
$$y^{l_2-1} \Bigl( y^{k_2-l_2 +1} h_1 (y) \varphi (y) + f_2 (x,y) \Bigr) = 0 $$
i.e. $f_2 (x,y) = x g_2 (x) - y^{k_2-l_2 +1} h_1 (y) \varphi (y)$ and therefore 
$$ f_1 (x,y) = - x^{l_1-k_1 +1} g_2 (x) \psi (x) + y h_1 (y)$$

By inspection we see that choosing $\bigl( g_2(x), h_1 (y) \bigr) = (1, 0)$ or 
$(0, 1)$ yields elements of the kernel of $df$.  Thus near $p$ the kernel of 
$df$ is generated by
\begin{eqnarray*}
x \left( - x^{l_1-k_1} \psi (x) \frac{du} {\alpha _1} + \frac{dv} {\alpha _2} \right) 
& \hspace{20pt} {\rm and} \hspace{20pt}
& y \Bigl( \frac{du} {\alpha _1} +  y^{k_2-l_2} \varphi (y)\frac{dv} {\alpha _2} \Bigr)
\end{eqnarray*}

We thus see again that $\cC _f$ is not locally free near $p$.  Since clearly the terms in 
brackets in the previous expression are local generators for $\cC _{f_1}$ and $\cC _{f_2}$ 
respectively near $p$, it follows that $\cC _f \cong \cC _{f_1} (-p) \oplus \cC _{f_1} (-p)$.  
Thus \ref{trr} is established, and with it the lemma.  \bo

Now that we have treated the transverse case, we need an analogous lemma for the 
non-transverse case.

\begin{lem} \label{conota}
Suppose that $X$ is a smooth surface and let $f: C \rightarrow X$ be a morphism from a 
curve $C$ consisting of two irreducible components $C_1$ and $C_2$, meeting in a node $p$.  
Denote by $f_i$ the restriction of $f$ to $C_i$ and let $p_i \in C_i$ be the point $p \in C$.  
Suppose that $f$ does not contract any component of $C$ and that $C_1$ and $C_2$ do not meet 
transversely at $f(p)$.  Then there are the following cases:
\begin{enumerate}
\item \label{n2} $C_1$ and $C_2$ are simply tangent at $f(p)$.

Then $\cC _f$ is not locally free and
$$ \cC _f \cong \cC _{f_1} (-p) \oplus \cC _{f_2} (-p) $$
\item \label{nl} $C_1$ and $C_2$ are not simply tangent at $f(p)$.

Then $\cC _f$ is locally free and the following sequence is exact
\begin{equation} \label{sesnf}
\xymatrix{ 0 \ar[r] & \cC _f \ar[r] & \cC _{f_1} \oplus \cC _{f_2} \ar[r] & \cC _{f,p} \ar[r] & 0 }
\end{equation}
\end{enumerate}
\end{lem}
{\it Proof.}  We proceed as before, and we can write
$$ f ^* : \left\{ \begin{tabular}{c@{ $\longmapsto $ }c}
$u$&$x^{k_1} U_1 (x) + y^{k_2} U_2 (y)$ \vphantom{$\frac{\frac{1}{1}}{\frac{1}{1}}$} \\
$v$&$x^{l_1} V_1 (x) + y^{l_2} V_2 (y)$ \vphantom{$\frac{\frac{1}{1}}{\frac{1}{1}}$} 
\end{tabular} \right. $$
where $l_1 > k_1$, $l_2 > k_2$ and $U_1 (0), U_2 (0) \neq 0$.  By exchanging if necessary the roles 
of $C_1$ and $C_2$, we may further assume that $k_1 \leq k_2$.  Then we have
\begin{equation*}
\xymatrix @C=-30.5pt @R=3pt { 
\cO _{C,p} \cdot du + \cO _{C,p} \cdot dv \ar[r]^{df \hspace{30pt}} & 
\raisebox {3pt} {$\Bigl( \cO _{C,p} \cdot dx + \cO _{C,p} \cdot dy \Bigr) $} / \raisebox 
{-3pt} {$(ydx + xdy)$} \\
du \ar[r] & x^{k_1-1} \Bigl( k_1 U_1(x) + x U_1 '(x) \Bigr) dx + 
            y^{k_2-1} \Bigl( k_2 U_2(y) + y U_2 '(y) \Bigr) dy \\
dv \ar[r] & x^{l_1-1} \Bigl( l_1 V_1(x) + x V_1 '(x) \Bigr) dx + 
            y^{l_2-1} \Bigl( l_2 V_2(y) + y V_2 '(y) \Bigr) dy }
\end{equation*}

Let $\alpha _1 := k_1 U_1(x) + x U_1 '(x)$ and 
$\alpha _2 := l_2 V_2(y) + y V_2 '(y)$.  We may write
$$ \xymatrix{ 
du \ar[r] & \alpha _1 x^{k_1-1} dx + \alpha _2 y^{k_2-1} dy \\
dv \ar[r] & x^{l_1-1} \psi (x) dx + y^{l_2-1} \varphi (y) dy } $$
Note that $\alpha _1 (0)$, $\alpha _2 (0) \neq 0$.  
The kernel of this morphism is determined by the condition
$$ f_1 (x,y) du + f_2 (x,y) dv \longmapsto r(x,y) \bigl( ydx + xdy \bigr) $$
which translates to
\begin{equation} \label{equata}
f ^* : \left\{ \begin{tabular}{c@{ $=$ }c}
$x^{k_1-1} \left( \alpha _1 f_1 (x,y) + x^{l_1-k_1} \psi (x) f_2 (x,y) 
 \vphantom{\frac {\frac{1}{1}} {\frac{1}{1}}} \right) $&$y r(x,y) = y r(0,y) $ \vspace{10pt} \\
$y^{k_2-1} \left( \alpha _2 f_1 (x,y) + y^{l_2-k_2} \varphi (y) f_2 (x,y) 
 \vphantom{\frac{\frac{1}{1}}{\frac{1}{1}}} \right) $&$x r(x,y) = x r(x,0) $
\end{tabular} \right.
\end{equation}

Let us now consider separately some cases.

\noindent
{\bf Case \ref{n2}}.  $k_1 = 1$ and $l_2 = k_2 + 1$.  (i.e. $f$ is not ramified on $C_1$ 
and $f(C_1)$ and $f(C_2)$ are simply tangent).  We know we may also assume $l_1 \geq 3$.  
Equations (\ref{equata}) imply (multiplying the second one by $y$ if $k_2=1$)
\begin{eqnarray*}
\alpha _1 f_1 (x,y) + x^{l_1-1} \psi (x) f_2 (x,y) & = & y r(x,y) \\[7pt]
y \Bigl( \alpha _2 f_1 (x,y) + y \varphi (y) f_2 (x,y) \vphantom{y^{l_2-k_2}} \Bigr) & = & 0
\end{eqnarray*}

From the second equation we deduce that $f_1 (x,y) = x g(x) - \frac {y \varphi (y)} 
{\alpha _2} f_2 (x,y)$, and substituting in the first equation we find
$$ x g(x) = - \frac {x^{l_1-1} \psi (x)} {\alpha _1} f_2 (x,y) + \frac {y \varphi (y)} 
{\alpha _2} f_2 (x,y) + \frac {y r(x,y)} {\alpha _1} $$

This gives us the equations
\begin{eqnarray*}
x g(x) & = & - \frac {x^{l_1-1} \psi (x)} {\alpha _1} f_2 (x,y) \\[7pt]
\frac {y \varphi (y)} {\alpha _2} f_2 (x,y) & = & - \frac {y r(x,y)} {\alpha _1}
\end{eqnarray*}

It follows that $f_2 (x,y) = x h(x) - \frac{\alpha _2} {\alpha _1 \varphi (y)} r(x,y)$.  
Observe now that choosing $h(x)=0$ and $r (x,y) = 1$ gives the element $\left( \frac 
{x^{l_1-1} \alpha _2 \psi (x)} {\alpha _1 ^2 \varphi (y)} + \frac {y} {\alpha _1} \right) 
du - \frac{\alpha _2} {\alpha _1 \varphi (y)} dv $ whose image under $df$ is $y dx$ 
(remember we are assuming $l_1 \geq 3$), which is not zero.  Therefore, $r(x,y)$ (and 
hence $f_2$) cannot have a constant term, which implies that all elements of the kernel 
are combinations of
\begin{eqnarray*}
x \Bigl( - \frac {x^{l_1-1} \psi (x)} {\alpha _1} du + dv \Bigr)
& \hspace{20pt} {\rm and} \hspace{20pt} &
y \Bigl( - \frac {y \varphi (y)} {\alpha _2} du + dv \Bigr)
\end{eqnarray*}

Clearly these elements are also in the kernel of $df$ and the terms in the brackets are local 
generators for $\cC _{f_1}$ and $\cC _{f_2}$.  Thus $\cC _f \cong \cC _{f_1} (-p) \oplus 
\cC _{f_2} (-p)$.

\noindent
{\bf Case \ref{nl}a}. $k_1 = 1$, $l_2 \geq k_2 + 2$.  Equations (\ref{equata}) imply
\begin{eqnarray*}
\alpha _1 f_1 (x,y) + x^{l_1-1} \psi (x) f_2 (x,y) & = & y r(x,y) \\[7pt]
y \Bigl( \alpha _2 f_1 (x,y) + y^{l_2-k_2} \varphi (y) f_2 (x,y) \Bigr) & = & 0
\end{eqnarray*}

From the second equation we deduce that
$$f_1 (x,y) = x g(x) - \frac {y^{l_2-k_2} \varphi (y)} {\alpha _2} f_2 (x,y)$$
and substituting in the first equation we find
$$ x g(x) = - \frac {x^{l_1-1} \psi (x)} {\alpha _1} f_2 (x,y) + \frac {y^{l_2-k_2} 
\varphi (y)} {\alpha _2} f_2 (x,y) + \frac {y r(x,y)} {\alpha _1} $$

This gives us the equations
\begin{eqnarray*}
x g(x) & = & - \frac {x^{l_1-1} \psi (x)} {\alpha _1} f_2 (x,y) \\[7pt]
\frac {y^{l_2-k_2} \varphi (y)} {\alpha _2} f_2 (x,y) & = & - \frac {y r(x,y)} {\alpha _1}
\end{eqnarray*}

Therefore all elements of the kernel are multiples of
$$ - \left( \frac {x^{l_1-1} \psi (x)} {\alpha _1} + \frac {y^{l_2-k_2} \varphi (y)} 
{\alpha _2} \right) du + dv $$

By inspection these elements are also in the kernel of $df$ and the restrictions $-\frac 
{x^{l_1-1} \psi (x)} {\alpha _1} du + dv $ and $-\frac {y^{l_2-k_2} \varphi (y)} {\alpha _2} 
du + dv $ are generators for $\cC _{f_1}$ and $\cC _{f_2}$.  In particular there is a short 
exact sequence as in (\ref{sesnf}).

\noindent
{\bf Case \ref{nl}b}.  $k_1, k_2 \geq 2$.  Then (\ref{equata}) implies $r(x,y) = 0$ and 
from the first equation we may write $f_1 (x,y) = -\frac {x^{l_1-k_1} \psi (x)} 
{\alpha _1} f_2 (x,y) + yh (y)$ and substituting in the second equation we obtain
$$ y^{k_2-1} \Bigl( y \alpha _2 h (y) + y^{l_2-k_2} \varphi (y) f_2 (x,y) \Bigr) = 0 
\hspace{4pt} \Longrightarrow \hspace{4pt} y h (y) = - \frac {y^{l_2-k_2} \varphi (y)} 
{\alpha _2 } f_2 (x,y) $$

Therefore near $p$ any element of the kernel of $df$ is a multiple of
$$ - \left( \frac {x^{l_1-k_1} \psi (x)} {\alpha _1} + \frac {y^{l_2-k_2} \varphi (y)} 
{\alpha _2 } \right) du + dv$$
and it is easy to check that this element lies indeed in the kernel of $df$.  Thus 
$\cC _f$ is locally free and since $- \frac {x^{l_1-k_1} \psi (x)} {\alpha _1} du + dv $ 
and $- \frac {y ^{l_2-k_2} \varphi (y)} {\alpha _2 } du + dv $ are the local generators 
for $\cC _{f_1}$ and $\cC _{f_2}$, we deduce that we have a short exact sequence as in 
(\ref{sesnf}).  This completes the proof of the lemma.  \bo

Let $f: C \rightarrow X$ be a morphism from a connected, projective, 
nodal curve of arithmetic genus zero to a smooth surface $X$.  In view 
of the two previous lemmas, we partition the set of nodes of $C$ in five 
disjoint sets:
\begin{description}
\item{$\tau _{uu} $} is the set of nodes $p$ such that the two components 
of $C$ meeting at $p$ are transverse at $f(p)$ and both are unramified;
\item{$\tau _{ur} $} is the set of nodes $p$ such that the two components 
of $C$ meeting at $p$ are transverse at $f(p)$ and one is unramified and 
the other one is ramified;
\item{$\tau _{rr} $} is the set of nodes $p$ such that the two components 
of $C$ meeting at $p$ are transverse at $f(p)$ and both are ramified;
\item{$\nu _2 $} is the set of nodes $p$ such that the two components of 
$C$ meeting at $p$ are simply tangent at $f(p)$;
\item{$\nu _l $} is the set of nodes $p$ such that the two components of 
$C$ meeting at $p$ are not transverse and not simply tangent at $f(p)$.
\end{description}

Thus it follows from the lemmas that the sheaf $\cC _f$ is locally free at the nodes $\tau _{uu}$ 
and $\nu _l$, while it is not free at the others.  Let $C_1, \ldots , C_\ell $ be the components of $C$.  
Then we let $\tau _{uu} ^i $ denote the divisor on $C_i$ of nodes lying in $\tau _{uu}$, and similarly 
for the other types of nodes.  Note that only one of the definitions above is not symmetric, namely 
$\tau _{ur}$ (and $\tau _{ur} ^i$).  To take care of this, let us introduce one more divisor on each 
component of $C$: 
let $\tau _{ru} ^i$ be the divisor on $C_i$ consisting of all nodes $p$ of $C$ on $C_i$, such that the 
two components of $C$ through $p$ are transverse at $f(p)$, and the restriction of $f$ to these two 
components is ramified only on $C_i$.

Often we will denote by the same symbol a divisor on a curve and its degree.  
For instance we write equations like
\begin{eqnarray*}
\sum _i \left( \tau _{uu} ^i + \tau _{ur} ^i + \tau _{rr} ^i + \nu _2 ^i + \nu _l ^i \right) 
& = & 2 \# \{ \text{nodes of } C \} \\
\sum _i \left( \tau _{uu} ^i + \tau _{ur} ^i + \tau _{ru} ^i + \tau _{rr} ^i + \nu _2 ^i + \nu _l ^i \right) 
& = & 2 \# \{ \text{nodes of } C \} + \tau _{ur}
\end{eqnarray*}

Given a coherent sheaf $\cF $ on a curve $C$, let $\tau (\cF)$ denote 
the subsheaf generated by the sections whose support has dimension at 
most 0 and let $\cF ^{free}$ be the sheaf $\cF / \tau (\cF )$.  By 
definition the sheaf $\cF ^{free}$ is pure.

\begin{prop} \label{grafico}
Let $f: C \rightarrow X$ be a stable map of genus zero with no contracted 
components to a smooth surface $X$, with canonical divisor $K_X$.  Let 
$C_1, \ldots , C_\ell $ be the irreducible components of $C$.  
Then we have
\begin{eqnarray} \label{graco}
\deg \left( \left. \left( \cC _f \otimes \omega _C \right) \vphantom{^{2}} \right| ^{free} _{C_i} 
\right) & = & 
f_* [C_i] \cdot K_X - \deg \tau _{ru} ^i + \deg \nu _l ^i + \deg \cQ _i \\[5pt]
\nonumber
\chi \left( \cC _f \otimes \omega _C \right) & = & 
f_* [C] \cdot K_X + \deg \tau _{rr} + \deg \nu _2 + \\\nonumber
& & + 2 \deg \nu _l + \sum \deg \cQ _i + 1
\end{eqnarray}

Moreover, let $\nu : \tilde C \rightarrow C$ be the normalization of $C$ at the nodes in 
$\tau _{ur} \cup \tau _{rr} \cup \nu _2$.  Then, the sheaf $\cC _f$ is the pushforward of a 
locally free sheaf on $\tilde C$.
\end{prop}
{\it Proof.}  This is simply a matter of collecting the information 
we already proved in the previous lemmas.  Thanks to Lemma 
\ref{conotra} and Lemma \ref{conota} we have the following 
short exact sequence of sheaves on $C$
$$\xymatrix@C=15pt { 
0 \ar[r]& \cC _f \ar[r]& \oplus _i \cC _{f_i} \left( -\tau _{uu} ^i - \tau _{ur} ^i 
- 2 \tau _{ru} ^i - \tau _{rr} ^i - \nu _2 ^i \right) \ar[r]&
\cC _f|_{\tau _{uu}} \oplus \cC _f|_{\nu _l} \ar[r]& 0 }$$

Note that the sheaf in the middle on the component $C_i$ is 
twisted down by all nodes of $C$ on $C_i$, with the exception 
of the nodes in $\nu _l ^i$, which do not appear, and the 
nodes in $\tau _{ru} ^i$, which ``appear twice.''  Hence we 
can write the divisor by which we are twisting $\cC _{f_i}$ 
as $- val [C_i] - \tau _{ru} ^i + \nu _l ^i$ (we denote by $val [C_i]$ 
the valence of the vertex $[C_i]$ in the dual graph of $C$).

To compute the degree of the sheaf $\cC _{f_i}$, remember 
that there is an exact sequence
$$\xymatrix{ 0 \ar[r]& \cC _{f_i} \ar[r]& f_i ^* \Omega ^1 _X \ar[r]& \Omega ^1 _{C_i} \ar[r]& 
\cQ _i \ar[r]& 0 } $$

Therefore we have $\deg \cC _{f_i} = f_*[C_i] \cdot K_X + 2 + \deg \cQ _i$.  
Thus, we may rewrite the previous sequence as follows
\begin{eqnarray*}
& \hspace{-2pt} \xymatrix@C=15pt { 0 \ar[r]& \cC _f \ar[r]& \oplus _i 
\cO _{C_i} \Bigl( f_*[C_i] \cdot K_X +2-val [C_i] 
- \tau _{ru} ^i + \nu _l ^i + \deg \cQ _i \Bigr) 
\ar[r]& } \\[7pt]
& \hspace{-2pt} \xymatrix@C=15pt {\ar[r] & 
\cC _f|_{\tau _{uu}} \oplus \cC _f|_{\nu _l} \ar[r]& 0 }
\end{eqnarray*}

The dualizing sheaf $\omega _C$ is invertible and on the component 
$C_i$ has degree equal to $-2+val [C_i]$.  Thus twisting the previous 
sequence by $\omega _C$ we obtain (using the isomorphisms 
$\cC _f|_{\tau _{uu}} \otimes \omega _C \simeq \cC _f|_{\tau _{uu}}$ 
and $\cC _f|_{\nu _l} \otimes \omega _C \simeq \cC _f|_{\nu _l}$)
\begin{eqnarray}
& \xymatrix@C=15pt { 0 \ar[r] & \cC _f \otimes \omega _C \ar[r] & \oplus _i 
\cO _{C_i} \Bigl( f_*[C_i] \cdot K_X - \tau _{ru} ^i + \nu _l ^i + \deg \cQ _i \Bigr) 
\ar[r] & } \nonumber \\[7pt] \label{gallo}
& \xymatrix@C=15pt { \ar[r] & \cC _f|_{\tau _{uu}} \oplus \cC _f|_{\nu _l} \ar[r] & 0 }
\end{eqnarray}

The first identity in (\ref{graco}) follows.  For the second 
one, note that $\sum \tau _{ru}^i = \tau _{ur}$ and 
$\sum \nu _l ^i = 2 \nu _l$ and compute Euler characteristics 
of (\ref{gallo}):
\begin{eqnarray*}
\chi \left( \cC _f \otimes \omega _C \right) & \hspace{-8pt} = & \hspace{-8pt} 
\sum _i \left( f_*[C_i] \cdot K_X - \tau _{ru} ^i + \nu _l ^i + \deg \cQ _i  +1 \right) \!-\! 
\deg \tau _{uu} \!-\! \deg \nu _l \!= \\
& \hspace{-8pt} = & \hspace{-8pt} f_*[C] \cdot K_X - \tau _{ur} + 2 \nu _l + \sum _i \deg \cQ _i + \\
& & \hspace{-8pt} + \# \bigl\{ \text{components of $C$} \bigr\} - \tau _{uu} - \nu _l
\end{eqnarray*}

\noindent
Remember that the dual graph of $f$ is a tree and hence 
$\# \bigl\{ \text{components}\bigr \} = \# \bigl\{ \text{nodes of $C$} \bigr\} + 1 = 
\tau _{uu} + \tau _{ur} + \tau _{rr} + \nu _2 + \nu _l + 1$.  
Using this, we conclude
$$ \chi \left( \cC _f \otimes \omega _C \right) = 
f_*[C] \cdot K_X + \tau _{rr} + \nu _2 + 2 \nu _l + 1 + \sum _i \deg \cQ _i $$
and the proposition is proved.  \bo

The next proposition has a similar proof, but deals with morphisms with 
contracted components.  As for the previous case, it is useful to 
introduce two more subsets of the nodes on contracted components, 
depending on the behaviour of $f : \bar C \rightarrow X$ near the node.  
We let
\begin{description}
\item{$\rho _{u} $} be the set of nodes $p$ such that $f$ is constant on 
one of the two components, and it is unramified on the other;
\item{$\rho _{r} $} be the set of nodes $p$ such that $f$ is constant on 
one of the two components, and it is ramified on the other.
\end{description}

\begin{prop} \label{graficone}
Let $f: \bar C \rightarrow X$ be a stable map of genus zero to a smooth surface $X$, with 
canonical divisor $K_X$.  Let $\bar C = C \cup R$, where $C=C_1 \cup \ldots \cup C_\ell$ 
is the union of all components of $\bar C$ which are not contracted by $f$, and $R$ is 
the union of all components of $\bar C$ contracted by $f$.  Let $r$ be the number of 
connected components of the curve $R$ (equivalently, $r = \chi (\cO _R)$).  Then we have
\begin{eqnarray} \label{gracone}
\deg \left( \left. \left( \cC _f \otimes \omega _{\bar C} \right) \vphantom{^{2}} 
\right| _{C_i} ^{free} \right) & = & 
f_* [C_i] \cdot K_X + \cQ _i - \tau _{ru} ^i + \nu _l ^i + \rho _u ^i + \rho _r ^i \\[5pt]
\nonumber
\chi \left( \cC _f \otimes \omega _{\bar C} \right) & = & 
f_* [C] \cdot K_X + \sum \cQ _i + 1 + \tau _{rr} + \nu _2 + 2 \nu _l + \\ \nonumber
& & + \rho _u + 2 \rho _r - 3r
\end{eqnarray}
\end{prop}

\noindent
{\it Proof.}  For the first formula in (\ref{gracone}), we only need to 
check the local behaviour of $\cC_f$ near a node between $C_i$ 
and a contracted component $R_j$.  As before, let $x$ be a 
local coordinate on $C_i$ near the node $p$ between $C_i$ and 
$R_j$ and let $y$ be a local coordinate on $R_j$ near $p$.  
Let $u,v$ be local coordinates on $X$ near $f(p)$ and suppose 
that the tangent direction to the vanishing set of $u$ near 
$f(p)$ is the tangent direction to $C_i$ at $f(p)$.  We have 
$$ f ^* : \left\{ \begin{array} {rcl}
u & \longmapsto & x^{k} U(x) \\
v & \longmapsto & x^{k+1} V(x) 
\end{array} \right. ~~~ U(0) \neq 0 $$
for some $k \geq 1$.  The sheaf $\cC_f$ near $p$ is the kernel 
of the map 
$$ \xymatrix @C=30pt { df : 
\cO _{\bar C, p} \cdot du + \cO _{\bar C, p} \cdot dv \ar[r] &
\raisebox {3pt} {$\Bigl( \cO _{\bar C, p} \cdot dx + 
\cO _{\bar C, p} \cdot dy \Bigr) $} / \raisebox 
{-3pt} {$\bigl( ydx + xdy \bigr)$} \\
du \ar[r] & x^{k-1} \Bigl( k U(x) + x U '(x) \Bigr) dx \\
dv \ar[r] & x^{k} \Bigl( (k+1) V(x) + x V '(x) \Bigr) dx } $$
It is readily seen that 
$$ x \Bigl( (k+1) V(x) + x V '(x) \Bigr) du - 
\Bigl( k U(x) + x U '(x) \Bigr) dv $$
is a local generator for the kernel of $df$.  Note that this 
means that we may pretend that the component $R_j$ is not there 
for the purpose of computing the contribution of the node $p$, 
regardless of whether $f|_{C_i}$ ramifies or not at $p$.  
This is enough to prove the first formula in (\ref{gracone}).

To prove the second one, we carry the previous analysis slightly 
further, and note that the image of $df$ contains the torsion 
section $y dx$ if and only if $f$ does not ramify at $p$.  Remember 
that we have the diagram 
$$ \xymatrix { & & & 0 \ar[d] \\
& & 0 \ar[d] & \tau \ar[d] \\
0 \ar[r] & \cC_f \ar[r] & f^* \Omega ^1_X \ar[r] \ar[d] & 
\Omega ^1 _{\bar C} \ar[r] \ar[d] & \cQ _{\bar C} \ar[r] & 0 \\
& & \mathop {\bigoplus } \limits _{C' \subset \bar C} 
\left( f|_{C'} \right) ^* \Omega ^1 _X \ar[d] & 
\mathop {\bigoplus } \limits _{C' \subset \bar C} \Omega ^1 _{C'} \ar[d] \\
& & \mathop {\bigoplus } \limits _{\nu \in Sing(\bar C)} 
\hspace{-10pt} \Omega ^1 _{X , \nu } \ar[d] & 0 \\
& & 0 } $$
where $C'$ ranges over the irreducible components of $\bar C$ and 
$\tau $ denotes the torsion subsheaf of $\Omega ^1 _{\bar C}$.  
We deduce that 
$$ \chi \left( \cC _f \otimes \omega _{\bar C} \right) 
= \chi \left( f^* \Omega ^1_X \otimes \omega _{\bar C} \right) 
- \chi \left( \Omega ^1 _{\bar C} \otimes \omega _{\bar C} \right) 
+ \chi \left( \cQ _{\bar C} \otimes \omega _{\bar C} \right) $$
and we know that 
\begin{eqnarray*}
\chi \left( f^* \Omega ^1_X \otimes \omega _{\bar C} \right) & \hspace{-5pt}= & \hspace{-5pt}
f_*[\bar C] \cdot K_X - 4 \# \{ {\text {components of }}\bar C \} 
+ 2 \sum _{C' \subset \bar C} val (C') + \\
& & \hspace{-5pt} + 2 \# \{ {\text {components of }}\bar C \} 
- 2 \# \{ {\text {nodes of }}\bar C \} = \\
& \hspace{-5pt} = & \hspace{-5pt} f_*[\bar C] \cdot K_X - 2
\end{eqnarray*}
$$ \chi \left( \Omega ^1 _{\bar C} \otimes \omega _{\bar C} \right) = 
\# \{ {\text {nodes of }}\bar C \} - 3 \# \{ {\text {components of }}\bar C \} + 
\sum _{C' \subset \bar C} val (C') = -3 $$
$$ \chi \left( \cQ _{\bar C} \otimes \omega _{\bar C} \right) = 
\chi \left( \cQ_C \otimes \omega _{\bar C} \right) + 
\chi \left( \Omega ^1_R \otimes \omega _{\bar C} \right) + \rho _r $$
where $\cQ_C$ is the cokernel of the differential of the restriction of 
$f$ to the union $C$ of the non-contracted components.  By what we saw 
above, the sheaf $\cQ_C$ behaves like when there are no contracted 
components.  The Euler characteristic of 
$\Omega ^1_R \otimes \omega _{\bar C}$ is given by 
\begin{eqnarray*}
\chi \left( \Omega ^1_R \otimes \omega _{\bar C} \right) & = & 
\# \{ {\text {nodes of }}\bar R \} 
- 3 \# \{ {\text {irreducible components of }}\bar R \} + \\[5pt]
& & + \sum _{R' \subset R} val (R') = 
- 3 \# \{ {\text {connected components of }}\bar R \} = \\[7pt]
& = & -3r
\end{eqnarray*}

We collect all these numbers as we did before and conclude.  \bo

\subsection{Deformations of Stable Maps}

\subsection{Dimension Estimates}

In what follows we refer to the integer $- C \cdot K_X$ as 
the anticanonical degree (or simply as the degree) of a curve 
$C$ in $X$, where $K_X$ is the canonical divisor of $X$.

We consistently use the following notational convention: if 
$f: \bar C \rightarrow X$ is a morphism and $\bar C_1$ denotes a component 
of $\bar C$, we denote the image of $\bar C_1$ by $C_1$, and in general, 
a symbol with a bar over it denotes an object on the source curve $\bar C$, 
while the same symbol without the bar over it denotes the image of the object 
in $X$.

\begin{defi} (\cite{Ko} II.3.6).  
Let $f,g \in {\rm Hom } (\bar C , X)$; we say $g$ is a deformation of $f$, 
if there is an irreducible subscheme of ${\rm Hom } (\bar C , X)$ containing 
$f$ and $g$.  We say that a general deformation of $f$ has some property if 
there is an open subset $U \subset {\rm Hom } (\bar C , X)$ containing $f$ 
and a dense open subset $V \subset U$ such that all $f' \in V$ have that 
property.
\end{defi}

When we choose a general deformation $g$ of a morphism $f$, we assume 
that $g$ is a deformation of $f$, i.e. that $f$ and $g$ lie in the same 
irreducible component of ${\rm Hom } (\bar C , X)$.

\begin{lem} \label{immersione}
Let $f: \bP^1 \rightarrow X$ be a free morphism; then 
$-f(\bP^1) \cdot K_X \geq 2$.  If moreover $f$ is 
birational onto its image, then a general deformation 
of $f$ is free and it is an immersion.
\end{lem}
{\it Proof.}  Since $f$ is free, $f^*\cT_X$ is globally generated, and hence the 
normal sheaf $\cN _f$ is also.  Thus we have 
$$ 0 \leq \deg \cN_f = \deg f^*\cT_X - 2 = - f(\bP^1) \cdot K_X - 2 $$

For the second assertion, by \cite{Ko} Complement II.3.14.4, a general 
deformation of $f$ is of the form $f_t : \bP^1 \stackrel{g_t} {\rightarrow } 
\bP^1 \stackrel {h_t} {\rightarrow } X$, where $h_t$ is an immersion.  
Since it is also true that a general deformation of a birational map 
is still birational, we see that for a general deformation $f_t$ 
of $f$, $g_t$ is an isomorphism, and $f_t$ is an immersion.  Clearly being 
free is also an open property.  \bo

Fix a free rational curve $\beta \subset X$ and let $d = - \beta \cdot K_X$.

\begin{defi}
Denote by $\barbieta $ the closure in $\barbeta $ of the set of free 
morphisms $f: \bP^1 \rightarrow X$ such that $f$ is birational onto 
its image.
\end{defi}

We want to prove that given $r \leq d-1$ general points 
$p_1, \ldots , p_r \in X$, in all irreducible components of $\barbieta$ 
there is an $f$ whose image contains all the $p_i$'s.

\begin{prop} \label{rpunti}
Let $f: \bP^1 \rightarrow X$ be an immersion, and let $d$ be the degree of the 
image of $f$.  Let $c_1, \ldots , c_r $ be distinct points where $f$ is an 
embedding.  The natural morphism 
$$ \xymatrix @R=10pt { F^{(r)} : (\bP^1)^r \times {\rm Hom} (\bP^1,X) \ar[r] & X^r \\
\bigl( d_1, \ldots , d_r ; [g] \bigr) \ar@{|->} [r] & 
\bigl( g(d_1) , \ldots , g(d_r) \bigr) } $$
is smooth at the point $\bigl( c_1, \ldots , c_r ; [f] \bigr)$ if and only 
if $r\leq d-1$.
\end{prop}
{\it Proof.}  Recall the sequence defining $\cN _f$:
\begin{equation} \label{dinuovo}
\xymatrix{ 0 \ar[r] & \cT \bP^1 \ar[r] ^{df^t} & f^* \cT X \ar[r]& 
\cN _f \ar[r] & 0 }
\end{equation}

Let us prove first of all that ${\rm Hom} (\bP^1,X)$ is smooth at $[f]$.  From 
(\ref{dinuovo}), it follows that $\deg \cN _f = -f_* (\bP^1) 
\cdot K_X -2 = d-2 \geq -1$, and since $f$ is an immersion, $\cN _f$ is locally 
free.  Thus from the long exact sequence associated to (\ref{dinuovo}) we deduce 
that ${\rm H}^1 (\bP^1 , f^* \cT X ) = 0$, and by \cite{Ko} Theorem II.1.2 it 
follows that ${\rm Hom} (\bP^1,X)$ is smooth at $[f]$.

Consider now the following commutative diagram with exact rows
$$ \xymatrix @C=15.3pt { 0 \ar[r] & \oplus \cT _{c_i} \bP^1 \ar[r] \ar [d]^{\sim} & 
\cT _{(\underline {c}; [f])} (\bP^1)^r \times {\rm Hom} (\bP^1,X) \ar[r] \ar[d]^{dF^{(r)}} & 
\cT _{[f]} {\rm Hom} (\bP^1,X) \ar[r] \ar[d]^{\delta} & 0 \\
0 \ar[r] & \oplus \cT _{f(c_i)} f(\bP^1) \ar[r] & \oplus \cT _{f(c_i)} X 
\ar[r] & \oplus \cN _{f,c_i} \ar[r] & 0} $$
The top row is clear, since we have the isomorphism 
$$\cT _{(\underline {c}; [f])} (\bP^1)^r \times 
{\rm Hom} (\bP^1,X) \simeq \oplus \cT _{c_i} \bP^1 \oplus \cT _{[f]} {\rm Hom} (\bP^1,X)$$
For the second row, restrict the sequence (\ref{dinuovo}) to 
$\{ c_1 , \ldots , c_r \}$ and note that $f$ induces an isomorphism $\cT _{c_i} \bP^1 
\simeq \cT _{f(c_i)} f(\bP^1)$, since $f$ is an embedding at the $c_i$.  The first 
vertical arrow is induced by $f$, while $\delta $ 
is the quotient map, followed by the evaluation map (\cite{Ko} Proposition II.3.5):
$$\xymatrix{ \cT _{[f]} {\rm Hom} (\bP^1,X) \simeq {\rm H}^0 (\bP ^1, f^* \cT _X) 
\ar[r] ^{\hspace{45pt}q} \ar[dr] _\delta & {\rm H}^0 (\bP ^1, \cN _f) \ar[d] ^{ev} \\
& \oplus \cN _{f,c_i} }$$

The morphism $q$ is induced by the long exact sequence associated to (\ref{dinuovo}), 
and the next term in the sequence is ${\rm H}^1 (\bP ^1, \cT _{\bP^1}) = 0$.  
Therefore $q$ is surjective.  Observe that $dF^{(r)}$ is surjective if and only if 
$\delta $ is surjective, and finally, $\delta $ is surjective if and only if the 
evaluation map $ev$ is surjective.  Consider the exact sequence of sheaves
\begin{equation} \label{sse}
\xymatrix{ 0 \ar[r] & \cN _f (-c_1 - \ldots - c_r) \ar[r] & \cN _f \ar[r] & 
\oplus \cN _{f,c_i} \ar[r] & 0 }
\end{equation}

Remember that $\deg \cN _f = d-2$, and since $f$ is an immersion, 
$\cN _f \simeq \cO _{\bP ^1} (d-2)$.  Thus ${\rm H}^1 ( \bP^1 , \cN _f) = 0$, 
and the sequence on global sections induced by (\ref{sse}) is exact if and only 
if ${\rm H}^1 ( \bP^1 , \cN _f (-c_1 - \ldots - c_r)) = 0$, i.e. if and only 
if $\deg \cN _f (-c_1 - \ldots - c_r) = d-2-r \geq -1$.  Therefore 
${\rm H}^0 (\bP ^1, \cN _f) \rightarrow \oplus \cN _{f,c_i}$ 
is surjective if and only if $r \leq d-1$, and hence $dF^{(r)}$ is surjective 
if and only if $r \leq d-1$.  \bo

Let $f: \bP^1 \rightarrow X$ be an immersion representing an element of $\barbieta$, 
and denote by $_f \barbieta $ the irreducible component of $\barbieta $ containing 
$f$.  Denote by $\cH ^f \subset {\rm Hom} (\bP^1 , X)$ the irreducible component of 
${\rm Hom} (\bP^1 , X)$ containing $[f]$ (remember that ${\rm Hom} (\bP^1 , X)$ 
is smooth at $[f]$).

There is an action
$$ \xymatrix @R=0pt {
{\rm Aut} \bigl( \bP^1 \bigr) \times \left( \bP^1 \right) ^r \times 
{\rm Hom} \bigl( \bP^1 , X \bigr) \ar[r] 
& \left( \bP^1 \right) ^r \times {\rm Hom} \bigl( \bP^1 , X \bigr) \\
\bigl( \varphi , (c_1, \ldots , c_r \,;\, [g]) \bigr) \ar@{|->} [r] & 
\bigl( \varphi (c_1), \ldots , \varphi (c_r) \,;\, [g \circ \varphi ^{-1}] \bigr) }$$
which clearly preserves the irreducible components of 
${\rm Hom} \bigl( \bP^1 , X \bigr)$.  Since $f$ is not constant, the action of 
${\rm Aut} \bigl( \bP^1 \bigr)$ has finite stabilizers.

Consider the diagram
$$\xymatrix{ & \left( \bP^1 \right) ^r \times \cH ^f 
\ar[dr] _{M \vphantom{^{F^{(r)}}}} \ar[dl] ^{F^{(r)}} \\
X^r & & _f \barbieta }$$
where $M$ is the projection onto the factor $\cH ^f $ followed by the natural map that 
quotients out the action of ${\rm Aut} (\bP^1)$.

Let us compute the dimensions of some of these spaces.  The 
morphism $M$ is obviously dominant, while Proposition 
\ref{rpunti} (together with Lemma \ref{immersione}) implies 
that $F^{(r)}$ is dominant if $r\leq d-1$.  Thus we may compute
$$ \dim \left( _f \barbieta \right) = \dim \Bigl( \left( \bP^1 \right) ^r \times \cH \Bigr) 
- r - 3 = -f(\bP^1) \cdot K_X - 1 = d - 1 $$

Let $c_1, \ldots , c_r \in \bP^1$ be $r\leq d-1$ distinct points where $f$ 
is an isomorphism onto its image and let $p_i = f(c_i)$.

Let $p:= (c_1, \ldots , c_r; [f]) \in \left( \bP^1 \right) ^r \times {\rm Hom} 
(\bP^1 , X)$; it follows from Proposition \ref{rpunti} that
$$\dim \bigl( F^{(r)} \bigr) ^{-1} (p_1, \ldots , p_r) = r + \dim \cH ^f - 2r = 
-f(\bP^1) \cdot K_X + 2 - r = d-r+2 $$

Denote by $\barbieta (p_1, \ldots , p_r)$ the (closure of the) image under 
$M$ of $\bigl( F^{(r)} \bigr) ^{-1} (p_1, \ldots , p_r)$, alternatively 
$$ \barbieta (p_1, \ldots , p_r) := \left\{ [f] \in \barbieta ~ \Bigl| 
~ Image(f) \supset \{p_1, \ldots , p_r\} \right\} $$

Since ${\rm Aut} (\bP^1)$ acts with finite stabilizers on 
$\bigl( F^{(r)} \bigr) ^{-1} (p_1, \ldots , p_r)$, we may compute
\begin{equation} \label{dimdibarbi}
\dim \barbieta (p_1, \ldots , p_r) = d-r-1
\end{equation}

\subsection{Independent Points}

The next lemma analyzes the case of curves through $d-1$ general points.

\begin{lem} \label{dimme}
For a general $(d-1)-$tuple $(p_1, \ldots , p_{d-1})$ of points of $X^{d-1}$, 
all the morphisms in $\barbieta (p_1, \ldots , p_{d-1})$ are immersions.
\end{lem}
{\it Proof.}  Let $\cI \subset \left( \bP^1 \right) ^{d-1} \times \cH ^f$ be the set 
of all $d-$tuples $(c_1, \ldots , c_{d-1}; [g])$ such that $g$ is not an 
immersion; Lemma \ref{immersione} implies that $\cI $ is a proper closed 
subset of $\left( \bP^1 \right) ^{d-1} \times \cH ^f$.  Note that $\cI$ is 
${\rm Aut }(\bP^1)-$invariant.  Consider the morphism $F^{(d-1)}$.  
By Proposition \ref{rpunti} and Lemma \ref{immersione} 
this morphism is dominant, hence the general fiber of this morphism has dimension 
$d-1-f(\bP^1) \cdot K_X + 2 - 2(d-1) = d+2-d+1 = 3$, thus the fibers of this morphism 
are ${\rm Aut }(\bP^1)-$orbits, since they are stable under the action of 
${\rm Aut }(\bP^1)$.  If the general fiber of $F^{(d-1)}$ met $\cI$, then 
we would have
$$ \dim \cI \geq 2(d-1) + 3 = 2d +1 = (d-1)+(d+2) = 
\dim \left( \left( \bP^1 \right) ^{d-1} \times \cH ^f \right) $$
and $\cI$ would equal $\left( \bP^1 \right) ^{d-1} \times \cH ^f $, which 
contradicts Lemma \ref{immersione}.  Thus there is an open dense subset 
$\cU $ in $X^{d-1}$ not meeting the image of $\cI $.  For any 
$(d-1)-$tuple $(p_1, \ldots , p_{d-1}) \in \cU $ we have that 
$$\barbieta (p_1, \ldots , p_{d-1}) := 
M \left( \bigl( F^{(d-1)} \bigr) ^{-1} (p_1, \ldots , p_{d-1}) \right) 
\subset \barbieta $$
consists only of (finitely many) immersions.  \bo

We now want to prove that for a general choice of $d-2$ points 
on $X$, all the resulting morphisms in $\barbieta$ through 
them have reduced image.  To achieve this, let us first 
introduce the following notion.

\begin{defi} \label{indip}
We say that $r$ points $p_1, \ldots , p_r$ in $X$ are {\rm independent} 
if the following conditions hold:
\begin{enumerate}
\item no $k$ of them are contained in a rational curve of degree $k$; \label{finoak}
\item the normalization of a rational curve of degree $k$ in $X$ through $k-1$ of them 
is an immersion. \label{nok}
\end{enumerate}
\end{defi}

Proposition \ref{rpunti}, Lemma \ref{dimme} and the dimension estimates 
(\ref{dimdibarbi}) easily imply that for any $r \geq 1$ there are $r-$tuples 
of independent points if there are free rational curves of anticanonical degree 
$d \geq r+1$, and that there are rational curves of anticanonical degree $d$ 
through $r$ independent points if $d \geq r+1$.

We are ready to prove the following result.

\begin{lem} \label{nonnonred}
Let $C \subset X$ be a divisor of anticanonical degree $d \geq 3$ such that 
each reduced irreducible component is rational.   Let 
$p_1, \ldots , p_{d-2} \in C$ be a $(d-2)-$tuple of independent points.  The 
divisor $C$ has at most two irreducible components and it is reduced.
\end{lem}
{\it Proof.}  Denote by $C_1, \ldots , C_\ell$ the reduced irreducible components of 
$C$.  For each curve $C_i$ let $d_i$ be the degree of $C_i$, $m_i$ be the 
multiplicity of $C_i$ in $C$ and $\delta _i$ be the number of points 
$p_1, \ldots , p_{d-2}$ lying on $C_i$.  Then we have $ \sum m_i d_i = d$ and 
$\delta _i \leq d_i-1$.  Therefore
$$ d-2 = \sum \delta _i \leq \sum d_i - \ell \leq \sum m_i d_i - \ell = d-\ell $$

Thus $\ell \leq 2$, and if $\ell = 2$, then all inequalities are equalities 
and hence $m_1 = m_2 = 1$.  If $\ell = 1$, then $C_1$ is a rational curve of 
degree $d_1$ on $X$ containing $d-2$ independent points.  It follows that 
$d_1 \geq d-1$ and $m_1 d_1 = d$ and hence $d \geq m_1 (d-1)$, or 
$(m_1-1) d \leq m_1$.   Since $d \geq 3$ this implies $d_1 = d$ and $m_1=1$.  \bo

\begin{lem} \label{niette}
Let $p_1, \ldots , p_r \in X$ be $r \geq 2$ independent points, and let 
$\alpha \subset X$ be an integral curve of degree $r+2$ of geometric genus zero 
containing $p_1 , \ldots , p_r$.  Let $B$ be a smooth connected projective 
curve and let $F: B \rightarrow \overline \cM _{bir} ^\alpha (p_1, \ldots , p_r)$ 
be a non-constant morphism.  The reducible curves in the family parametrized by 
$B$ cannot always contain a component mapped isomorphically to a curve of 
anticanonical degree strictly smaller than two.
\end{lem}
{\it Proof.}  Consider the following fiber product diagram
$$\xymatrix{ S \ar[rr] \ar[d] & & \overline \cM _{0,1} (X, \alpha ) \ar[d] \\
B \ar[r] & \overline \cM _{bir} ^\alpha (p_1, \ldots , p_r) \ar@{^(->} [r] & 
\overline \cM _{0,0} (X, \alpha) }$$
thus $S \rightarrow B$ is the pull-back of the universal family.

It follows that $S \rightarrow B$ is a surface whose general fiber over $B$ 
is a smooth rational curve and with a finite number of fibers consisting of 
exactly two smooth rational curves (Lemma \ref{nonnonred}) meeting 
transversely at a point, corresponding to the reducible curves in the 
family $B$.  By hypothesis $S \rightarrow B$ admits $r$ contractible 
sections.  Suppose that in all reducible fibers of $S$ one component is 
mapped to a curve of anticanonical degree strictly smaller than two.  
Denote the components in $S$ mapped to such curves by $L_1, \ldots , L_t$, 
and the other components in the respective fiber by $Q_1, \ldots , Q_t$ 
(thus $L_i + Q_i$ represents the numerical class of a fiber, for all $i$'s).  
By definition of independent points, the sections of $S \rightarrow B$ cannot 
meet the components $L_i$.  Since $L_i \subset S$ is a smooth rational curve 
of self-intersection $L_i^2 = L_i \cdot (Q_i + L_i) - L_i \cdot Q_i = -1$, we 
may contract all the $L_i$ to obtain a smooth surface $S' \rightarrow B$, 
which is a $\bP ^1-$bundle over the curve $B$.  Since the contracted curves 
did not meet the $r$ sections, there still are $r \geq 2$ negative sections 
of $S' \rightarrow B$, but there can be at most one negative section in a 
$\bP^1-$bundle.  Thus there must be reducible fibers in the family $B$ all of 
whose components are mapped to curves of anticanonical degree at least two.  \bo

\begin{lem} \label{sovrappo}
Let $f: \bP^1 \longrightarrow X$ be a non-constant morphism to a smooth 
surface $X$ and suppose that $f^* \cT _X$ is globally generated.  
Denote by $\cM _f$ the irreducible component of 
$\overline \cM _{0,0} (X, f_*[\bP^1])$ containing $[f]$ and by 
$C \subset X$ the integral curve $f(\bP^1)$.  Let $\cM _{f,C}$ be the 
locus of stable maps
$$ \cM _{f,C} := \biggl\{ [g] \in \cM _f ~ \bigl| ~ 
\text{image}(g) = C \biggr\} $$

Then we have 
$${\rm codim} \bigl( \cM _{f,C} , \cM _f \bigr) \leq 1$$
Equality holds if and only if $f_*[\bP^1] = \delta C$ for 
some positive integer $\delta $ and $K_X \cdot C = -2$.
\end{lem}
{\it Proof.}  Using \cite{Ko} Proposition II.3.7, we may deform $f$ so that the 
image of the resulting morphism avoids a point on $C$.  It follows that 
$\cM _{f,C} \subsetneq \cM _f$, and hence, $\cM _{f,C}$ being closed, that 
it has codimension at least one.

To prove the second assertion, note that any morphism $\phi : R \rightarrow X$ 
from a rational tree with image contained in $C$ is such that $\phi ^* \cT_X$ 
is globally generated.  This is obvious on each irreducible component of $R$: 
the morphism factors through the normalization of $C$ and a multiple cover, 
and under the normalization the pull-back of $\cT_X$ is globally generated.  
Thus $\phi ^* \cT_X$ is globally generated on each component of $R$, and hence 
it is globally generated on $R$.

Let $\Gamma $ be the dual graph of some morphism in $\cM _f$.  Let 
$\cM _f ^\Gamma $ be the subscheme of $\cM _f$ consisting of morphisms with 
dual graph $\Gamma $; then
\begin{equation} \label{disestu}
{\rm codim} \bigl( \cM _{f,C} \cap \cM _f ^\Gamma , \cM _f ^\Gamma \bigr) \geq 1
\end{equation}

\noindent
Indeed, let $n$ be the number of vertices of $\Gamma $ and consider the scheme 
$\widetilde \cM _f ^\Gamma $:
$$ \left\{ \bigl[ g : K \rightarrow X \,;\, p_1 , \ldots , p_n \bigr] \in 
\cM _{0,n} \bigl( X, f_*[\bP^1] \bigr) \, \left| \begin{tabular}{ll} 
$[g] \in \cM _f ^\Gamma  $ and the points \\
$p_1$, \ldots , $p_n$ lie in different \\
components of $K$
\end{tabular} \right. \hspace{-6pt} \right\} $$

Clearly there is a surjective morphism 
$\widetilde \cM _f ^\Gamma \longrightarrow \cM _f ^\Gamma $, and let 
$$\widetilde \cM _{f,C} ^\Gamma := \left( \cM _{f,C} \cap \cM _f ^\Gamma \right) 
\times _{\cM _f ^\Gamma }\widetilde \cM _f ^\Gamma $$
Let $g : K \rightarrow X$ represent a morphism in 
$\widetilde \cM _{f,C} ^\Gamma $; again by \cite{Ko} Proposition II.3.7 
we may deform $g$ to miss a point of $C$, while still lying in 
$\widetilde \cM _f ^\Gamma $ and thus (\ref{disestu}) follows.

Suppose that ${\rm codim} (\cM _{f,C} , \cM _f) = 1$.  It is clear that 
$f_* [\bP^1] = \delta C$ for some positive integer $\delta $.

Using (\ref{disestu}) it follows that the general morphism in every component 
of maximal dimension of $\cM _{f,C}$ has irreducible domain, and hence 
these components of $\cM _{f,C}$ are dominated by $\cM _{0,0} (\bP^1 , \delta )$, 
where the morphisms are induced by composition with the normalization map 
$\nu : \bP ^1 \longrightarrow C$.  We have 
$\dim \cM _{f,C} \leq \dim \cM _{0,0} (\bP^1 , \delta ) = 2 \delta - 2$, 
and also $\dim \cM _f = \left( - K_X \cdot C \right) \delta - 1$.  We 
already know (Lemma \ref{immersione}) that $- K_X \cdot C \geq 2$, and 
hence we must have $- K_X \cdot C = 2$ and 
$\dim \cM _{f,C} = 2 \delta - 2$.  \bo

\subsection{Sliding moves} \label{slittino}

The next lemma and its corollary allow us to construct irreducible 
subschemes in the boundary of the spaces $\barbeta $.  First, let us 
introduce some notation that will be used in the lemma.

Let $f : \bar C \rightarrow X$ be a stable map of genus zero to the 
smooth surface $X$.  Let $\bar C_0$ be a connected subcurve, let 
$\bar C_1 , \ldots , \bar C_\ell $ be the connected components of 
the closure of $\bar C \setminus \bar C_0$.  Let $\bar C_{0i}$ 
be the irreducible component of $\bar C_0$ meeting $\bar C_i$, and 
let $\bar C_{i,1}$ be the irreducible component of $\bar C_i$ 
meeting $\bar C_0$ and let the intersection point of $\bar C_{0i}$ 
and $\bar C_{i,1}$ be $\bar p_i$.  
Denote by $f_i$ the restriction of $f$ to $\bar C_i$, for 
$i \in \{ 0, \ldots , \ell \}$.

Let $V \subset \overline \cM _{0,\ell} \bigl( X , f_* [\bP^1]\bigr)  \times 
\left( \bar C_1 \times \cdots \times \bar C_\ell \right) $ be the 
subscheme consisting of all points 
$\bigl( [g \,;\, \bar c_1 , \ldots , \bar c_\ell ] \,;\, 
\bar c_1 ' , \ldots , \bar c_\ell '\bigr) $, 
such that $g(\bar c_i) = f (c_i ')$ and 
$[g \,;\, \bar c_1 , \ldots , \bar c_\ell ]$ is in the 
same irreducible component of 
$\overline \cM _{0,\ell} \bigl( X , f_* [\bP^1]\bigr) $ as 
$[f \,;\, \bar p_1 , \ldots , \bar p_\ell ]$.

\begin{lem} \label{ovvio}
With notation as above, assume also that a general deformation of 
$f_0$ is generated by global sections, $\bar C_{0i}$ is not contracted 
by $f$ and 
$f \bigl( \bar C_{0i}\bigr)  \not \supset f \bigl( \bar C_{i,1}\bigr) $, 
for all $i$'s.  It follows that every irreducible component of $V$ 
containing $\bigl( [f _0 \,;\, \bar p_1 , \ldots , \bar p_\ell ] \,;\, 
\bar p_1 , \ldots , \bar p_\ell \bigr) $ 
surjects onto the irreducible component of 
$\overline \cM _{0,0} \bigl( X , f_* [\bP^1]\bigr) $ containing $[f]$.
\end{lem}
{\it Proof.}  Let $\Phi $ be an irreducible component of 
$\overline \cM _{0,0} \bigl( X , f_* [\bP^1]\bigr) $ containing (the 
stable reduction of) $[f]$.  Define $\cC $ by the Cartesian square 
on the left and $\underline {ev}$ as the composite of the maps in 
the diagram
$$ \xymatrix { \cC \ar[rr] \ar[d] 
\ar@/^2pc/ ^{\underline {ev} := (ev _1 , \ldots , ev _\ell )} [rrrr] && 
{\overline \cM _{0,\ell} \bigl( X , f_* [\bP^1]\bigr) } \ar[d] \ar[rr] && X ^\ell \\ 
\Phi \ar@{^{(}->} [rr] && {\overline \cM _{0,0} \bigl( X , f_* [\bP^1]\bigr) }} $$

Clearly, $V$ is then defined by the diagram 
$$ \xymatrix { V \ar[rr] \ar[d] ^\iota \ar@/_2pc/ ^{\pi } [dd] && 
{\left( \bar C_1 \times \cdots \times \bar C_\ell \right) } 
\ar[d] ^{(f_1, \ldots , f_\ell )} \\
\cC \ar[rr] ^{\underline {ev} } \ar[d] && X ^\ell \\ 
\Phi } $$
and we have
$$\xymatrix { {}\save[]+<-78pt,0pt>*{ V \subset W := 
\cC \times \left( \bar C_1 \times \cdots \times \bar C_\ell \right) } \restore 
\ar[rr] ^{P} && \cC }$$

Obviously $P$ is flat and since $\cC \longrightarrow \Phi$ is flat, it 
follows that $W \longrightarrow \Phi $ is flat.  The fiber of $\pi $ 
at the point $[g]$ is given by 
$$\pi ^{-1} \bigl( [g] \bigr) = \biggl\{ 
\bigl( [\tilde g \,;\, \bar c_1 , \ldots , \bar c_\ell ] \,;\, 
\bar c_1 ' , \ldots , \bar c_\ell ' \bigr) 
~\Bigl|~ \tilde g (\bar c_i) = f_i (\bar c_i ') \biggr\}$$
where the stable reduction of $\tilde g$ is $g$.  If $g$ has irreducible 
domain, and if the image of $g$ does not contain any singular point of 
(the reduced scheme) $f \bigl(  \bar C_1 \cup \ldots \cup \bar C_\ell \bigr) $, 
nor does it contain any component of $f\bigl( \bar C_i\bigr) $, then the 
scheme $\pi ^{-1} \bigl( [g] \bigr) $ is finite.  Thanks to \cite{Ko} 
Theorem II.7.6 and Proposition II.3.7, a general deformation $g$ of $f _0 $ 
satisfies the previous conditions; thus the general fiber of $\pi $ in 
a neighbourhood of $[f]$ is finite and hence, letting 
$v_0 := \bigl( [f _0 \,;\, \bar p_1 , \ldots , \bar p_\ell ] \,;\, 
\bar p_1  , \ldots , \bar p_\ell \bigr) $, we conclude that 
$\dim _{v_0} V = \dim \Phi = \dim \cC - \ell $.

Let $\kappa _i \in \cO _{X, f_(\bar p_i)}$ be a local equation of 
$f_i \bigl( \bP^1\bigr) $; clearly the $\ell $ equations 
$P^* ev_1 ^* (\kappa _1 ) $, \ldots , $P^* ev_\ell ^* (\kappa _\ell )$ 
define $V$ near $v_0$.  Since $\dim V = \dim \cC - \ell $, it follows that 
$\cO _{V, v_0}$ is a Cohen-Macaulay $\cO _{W , v_0}-$module.  Using 
\cite{EGA4} Proposition 6.1.5, we deduce that $\cO _{V, v_0}$ is a flat 
$\cO _{\Phi , [f _0 ]}-$module, and the result follows.  \bo

\noindent
{\bf Construction.}  Suppose $f: \bar C \rightarrow X$ is a stable map, 
and suppose $\bar C = \bar C _0 \cup \ldots \cup \bar C_\ell $, where 
$\bar C_i$ is a connected union of components for all $i$'s, such that 
${\rm H}^1 \bigl( \bar C_0 , f^* \cT_X|_{\bar C_0} \bigr) = 0$ and all 
the irreducible components of $\bar C_0$ meeting $\bar C_i$ are not 
contracted by $f$ and the image of the component of $\bar C_0$ meeting 
$\bar C_i$ does not contain the image of the corresponding component 
of $\bar C_i$ for all $i$'s (this is the same condition required in 
Lemma \ref{ovvio}).

We construct an irreducible subscheme ${\rm Sl} _f (\bar C_0)$ of 
$\overline \cM _{0,0} \bigl( X , f_* [\bar C] \bigr) $, consisting 
of morphisms $g : \bar C' \rightarrow X$ with the following properties:
\begin{itemize}
\item there is a decomposition 
$\bar C' = \bar C_0' \cup \ldots \cup \bar C_\ell '$, where $\bar C_i'$ 
is a connected subcurve; 
\item there are isomorphisms $g|_{\bar C_i'} \simeq f|_{\bar C_i}$; 
\item there is a morphism $res : {\rm Sl} _f (\bar C_0) \rightarrow 
\overline \cM _{0,0} \bigl( X , f_* [\bar C_0] \bigr)$, which is 
surjective on the irreducible component containing $f|_{\bar C_0}$; 
\item there are morphisms $a_i : {\rm Sl} _f (\bar C_0) \rightarrow \bar C_i$, 
for $i \in \{ 1 , \ldots , \ell \}$.
\end{itemize}

Let $\bar p_i \in \bar C_0$ be the node between $\bar C_0$ and $\bar C_i$ 
and $f_i := f|_{\bar C_i}$; by Lemma \ref{ovvio} we may find an irreducible 
subscheme $V \subset 
\overline \cM _{0,\ell} \bigl( X , f _* [\bar C_0] \bigr) \times _{X ^\ell } 
\bigl( \bar C_1 \times \ldots \times \bar C_\ell \bigr)$ and a morphism 
$V \rightarrow \overline \cM _{0,0} \bigl( X , f _* [\bar C_0] \bigr)$ which is 
surjective onto the irreducible component containing $f_0$.

Identify $\bar C_i$ with 
$\overline \cM _{0,1} \bigl( \bar C_i , [\bar C_i] \bigr)$; thus we may write 
$$ V \subset 
\overline \cM _{0,\ell} \bigl( X , f _* [\bar C_0] \bigr) \times _{X ^\ell } 
\Bigl( \overline \cM _{0,1} \bigl( \bar C_1 , [\bar C_1] \bigr) \times \ldots 
\times \overline \cM _{0,1} \bigl( \bar C_\ell , [\bar C_\ell] \bigr) \Bigr) $$
Let $M_i \subset \bar C_0 \times P$ be the closed subscheme with closed points 
of the form 
$\bigl( \bar c_{0i} \,;\, [g \,;\, \bar c_{01} , \ldots , \bar c_{0\ell} ] 
\,;\, \bar c_1 , \ldots , \bar c_\ell \bigr)$.  
Let $N_i \subset \bar C_i \times P'$ be the closed subscheme with closed points 
of the form $\bigl( \bar c_i \,;\, [g \,;\, \bar c_{01} , \ldots , \bar c_{0\ell} ] 
\,;\, \bar c_1 , \ldots , \bar c_\ell \bigr)$.  
It is clear that projection onto the $P'$ factor induces isomorphisms 
$M_i \simeq P'$ and $N_i \simeq P'$, and that $M_i \cap M_j = \emptyset $ 
for all $i \neq j$.

Construct the scheme $\bar \cC $: glue to $\bar C_0 \times P'$ the schemes 
$\bar C_i \times P'$ along the subschemes $M_i \simeq N_i$, where the isomorphisms are the 
ones induced by projection onto the factor $P'$.  By construction, there is a morphism 
$\bar \cC \longrightarrow P'$, whose fiber over the point 
$\bar c = \bigl( [g \,;\, \bar c_{01} , \ldots , \bar c_{0\ell} ] \,;\, 
\bar c_1 , \ldots , \bar c_\ell \bigr)$ is the curve $\bar \cC _{\bar c}$ 
obtained by the nodal union of $\bar C_0$ and $\bar C_i$, for all $i$'s, where the nodes 
of $\bar \cC _{\bar c}$ are at the points $\bar c_{0i} \in \bar C_0$ and 
$\bar c_i \in \bar C_{i,1} \subset \bar C_i$.

The morphism $\bar \cC \rightarrow P'$ is flat on all irreducible components of $P'$ 
(remember that $P'$ is smooth) thanks to Theorem III.9.9 of \cite{Ha}, since all fibers 
$\bar \cC _{\bar c}$ have geometric genus zero.  Thus $\bar \cC \rightarrow P'$ is a family 
of connected nodal projective curves of arithmetic genus zero.

A typical application of this construction can be found in the proof of Theorem \ref{maschera} as 
well as in many of the later proofs.

\section{Divisors of Small Degree: the Picard Lattice}

\subsection{The Nef Cone}


We collect here some results on the nef cone of a del Pezzo surface.  
We prove a ``numerical'' decomposition of any nef divisor on a del 
Pezzo surface in Corollary \ref{maquale}.  In the later sections we 
will show how to realize geometrically this decomposition.

\begin{defi}
Let $X _\delta $ be a del Pezzo surface of degree $9-\delta$.  
Suppose that $X_\delta \neq \bP^1 \times \bP^1$.  We call an 
integral basis $\{ \ell , e_1 , \ldots , e_\delta \}$ of 
${\rm Pic} (X_\delta )$ a {\rm standard basis} if there is a 
presentation $b: X_\delta \rightarrow \bP^2$ of $X_\delta $ as 
the blow up of $\bP^2$ at $\delta $ points such that 
$\ell $ is the pull-back of the class of a line and the 
$e_i$'s are the exceptional divisors of $b$.
\end{defi}

\begin{lem} \label{dispari}
Let $C \subset X$ be an integral curve of canonical degree -1 on 
the smooth surface $X$.  Then $C^2$ is odd and it is at least -1.
\end{lem}
{\it Proof.}  This is immediate from the adjunction formula: 
$$ \begin{array} {c}
C^2 + K_X \cdot C = 2 p_a(C) - 2 \hspace{20pt} \Longrightarrow 
\hspace{20pt} C^2 = 2 p_a(C) - 1 \geq -1 
\end{array} $$
The lemma is proved.  \bo

\begin{lem} \label{meno}
Let $C \subset X$ be a curve of canonical degree -1 on a del Pezzo 
surface of degree $d$.  Either $C$ is a $(-1)-$curve, or $d=1$ and 
the divisor class of $C$ is $-K_X$.
\end{lem}
{\it Proof.}  Note that since $-K_X$ is ample, a curve of canonical degree 
$-1$ must be integral.  If $X \simeq \bP^1 \times \bP^1$, all divisor 
classes on $X$ have even canonical degree, thus we may exclude this 
case.  Let $\rho := C^2$ and $\delta = 9-d$; by the previous lemma 
we know that $\rho \geq -1$ and it is odd.  Moreover, if $\rho = -1$ 
then $C$ is a $(-1)-$curve; suppose therefore that $\rho \geq 1$.  
By \cite{Ma} Proposition IV.25.1 we may find a standard basis 
$\{ \ell , e_1 , \ldots , e_\delta \}$ of the Picard group of $X$.  
If we write $C = a\ell - b_1 e_1 - \ldots - b_\delta e_\delta $, 
we have 
$$ \left\{ \begin{array} {r@{ ~=~ }l} \displaystyle
3a - \sum _{i=1} ^\delta b_i & 1 \\[15pt] \displaystyle
a^2 - \sum _{i=1} ^\delta b_i ^2 & \rho 
\end{array} \right. $$
and these equations are easily seen to be equivalent to the 
following: 
$$ \left\{ \begin{array} {rcl} \displaystyle
3a - \sum _{i=1} ^8 b_i & = & 1 \\[7pt] \displaystyle
\sum _{i=1} ^8 \bigl( a - 2b_i - 1 \bigr) ^2 & = & 
4 \bigl( 1 - \rho \bigr) \\[7pt] 
\displaystyle b_i & = & 0 \hspace{20pt} i \geq \delta +1
\end{array} \right. $$

We deduce that $\rho \leq 1$, and hence $\rho = 1$.  We conclude 
that $a - 2b_i - 1 = 0$ for all $i$'s and hence 
$\bigl( a \,;\, b_1 , \ldots , b_8 \bigr) = 
\bigl( 2b+1 \,;\, b , \ldots , b \bigr)$ and $3a - \sum b_i = 1$.  
Therefore $b=1$, $\delta = 8$ and the divisor class of $C$ is 
$\bigl( 3 \ell - e_1 - \ldots - e_8 \bigr) = -K_X$.  \bo

We need a criterion to determine which classes are nef on any del 
Pezzo surface $X$.  This is immediate in the cases of del Pezzo 
surfaces of degree 9 and 8.  If the degree is 9, then $X$ is 
isomorphic to $\bP^2$.  The non-negative multiples of the 
class of a line are the only nef divisors, and the only ample 
divisors are the positive such multiples.  If the degree of the 
del Pezzo surface is 8, then there are two cases: either $X$ is 
isomorphic to $\bP^1 \times \bP^1$ or $X$ is isomorphic to the 
blow-up of $\bP^2$ at one point.

If $X \simeq \bP^1 \times \bP^1$, then any divisor class $C$ on 
$X$ is of the form $a_1 F_1 + a_2 F_2$, where $F_1$ and $F_2$ are 
the two divisor classes of $\{p\} \times \bP^1$ and 
$\bP^1 \times \{p\}$ and $a_1$ and $a_2$ are integers.  Then $C$ 
is nef if and only if $a_1 , a_2 \geq 0$, while $C$ is ample if and 
only if $a_1 , a_2 > 0$.

If $X \simeq Bl _p (\bP^2)$, then any divisor class $C$ on $X$ 
is of the form $a \ell - b e$, where $\ell $ is the pull-back of 
the divisor class of a line in $\bP^2$, while $e$ is the exceptional 
divisor.  The divisor class $C$ is nef if and only if $a\geq b \geq 0$, 
while $C$ is ample if and only if $a > b > 0$.

The remaining cases are dealt with in the next Proposition.

\begin{prop} \label{clane}
Let $X$ be a del Pezzo surface of degree $d \leq 7$.  A divisor 
class $C \in {\rm Pic} (X)$ is nef (respectively ample) if and 
only if $C \cdot L \geq 0$ (respectively $C \cdot L > 0$) for all 
$(-1)-$curves $L \subset X$.
\end{prop}
{\it Proof.}  The necessity of the conditions is obvious.  To establish 
the sufficiency, we only need to prove the result for nef classes, 
since the ample classes are precisely the ones in the interior of 
the nef cone.  Proceed by induction on $r := 9 - d$.

If $r=2$ write $C = a \ell - b_1 e_1 - b_2 e_2$, in some standard 
basis $\{\ell , e_1 , e_2 \}$.  By assumption we know that 
$b_i \geq 0$ and $a \geq b_1 + b_2$.  Thus we can write 
$$ C = \bigl( a - b_1 - b_2 \bigr) \ell + 
b_1 \bigl( \ell - e_1 \bigr) + b_2 \bigl( \ell - e_2 \bigr) $$
which shows that $C$ is a non-negative combination of nef classes.

Suppose $r>2$.  Let $n := \min \bigl\{ C \cdot L ~;~ 
L \subset X \text{ is a $(-1)-$curve} \bigr\}$; by assumption 
we know that $n \geq 0$.  Let $\tilde C := C + nK_X$; for any 
$(-1)-$curve $L \subset X$ we have 
$\tilde C \cdot L = C \cdot L -n \geq 0$, and there is a 
$(-1)-$curve $L'$ such that $\tilde C \cdot L' = 0$, by 
the definition of $n$.

Let $b : X \rightarrow X'$ be the contraction of the curve $L'$ 
and note that $X'$ is a del Pezzo surface of degree $9-(r-1)$.  
We have $\tilde C = b^* b_* \tilde C - r L'$ and 
$$0 = \tilde C \cdot L' = b^* b_* \tilde C \cdot L' - r L' \cdot L' 
= b_* \tilde C \cdot b_* L' + r = r $$
and therefore $\tilde C = b^*b_*\tilde C$ is the pull-back of the 
divisor class $C' := b_* \tilde C$ on $X'$.  Since all 
$(-1)-$curves on $X'$ are images of $(-1)-$curves on $X$, by 
induction we know that $C'$ is nef, and thus $\tilde C$ is 
nef.  Hence $C = \tilde C + n (-K_X)$ is a non-negative linear 
combination of nef divisors, and thus $C$ is nef.  \bo

From this Proposition we deduce immediately the following Corollary.

\begin{cor} \label{maquale}
Let $X_\delta $ be a del Pezzo surface of degree $9 - \delta \leq 8$.  
Let $D \in {\rm Pic} (X_\delta )$ be a nef divisor.  Then we can 
find 
\begin{itemize}
\item non-negative integers $n_2 , \ldots , n_\delta $;
\item a sequence of contraction of $(-1)-$curves 
$$ \xymatrix { X_\delta \ar[r] & X_{\delta -1} \ar[r] & {\cdots} \ar[r] & 
X_2 \ar[r] & X_1 } ; $$
\item a nef divisor $D' \in {\rm Pic} (X_1)$;
\end{itemize}
such that 
$$ D = n_\delta (-K_{X_\delta }) + n_{\delta -1} (-K_{X_{\delta -1}}) + 
\ldots + n_2 (-K_{X_2}) + D' $$
\end{cor}
{\it Proof.}  We proceed by induction on $\delta $.  If $\delta \leq 1$, 
there is nothing to prove.

Suppose that $\delta \geq 2$ and let 
\(n := \min \bigl\{ L \cdot D ~|~ 
L \subset X {\text{ a $(-1)-$curve} } \} \).  By assumption we have 
$n \geq 0$.  Let $\bar D := D + n K_{X_\delta}$; for every $(-1)-$curve 
$L \subset X_\delta $ we have 
$$ \bar D \cdot L = D \cdot L + n K_{X_\delta} \cdot L  \geq n - n = 0 $$

Thus thanks to the previous Proposition, $\bar D$ is nef.  By construction 
there is a $(-1)-$curve $L_0 \subset X$ such that $\bar D \cdot L_0 = 0$.  
Thus $\bar D$ is the pull-back of a nef divisor on the del Pezzo surface 
$X_{\delta -1}$ obtained by contracting $L_0$.  By induction, we have a 
sequence of contractions 
$$ \xymatrix { X_{\delta -1} \ar[r] & {\cdots} \ar[r] & X_2 \ar[r] & X_1 } , $$
non-negative integers $n_2$,\ldots , $n_{\delta - 1}$ and a nef divisor 
$D'$ on $X_1$ such that we may write 
$\bar D = n_{\delta - 1} (-K_{X_{\delta -1}}) + \ldots 
+ n_2 (-K_{X_2}) + D'$.  Let $n_\delta := n$; with this notation we have 
$$ D = n_\delta (-K_{X_\delta }) + \bar D' = n_\delta (-K_{X_\delta }) + \ldots 
+ n_2 (-K_{X_2}) + D' $$
and a sequence of contractions as in the statement of the corollary.  
This concludes the proof.  \bo

\subsection{First Cases of the Main Theorem}

\begin{thm} \label{pumba}
Let $X_\delta$ be a del Pezzo surface of degree $9-\delta$; then 
the linear system $|-K_{X_\delta}|$ has dimension $9-\delta$.  
If $\delta = 8$, then $|-K_{X_8}|$ has a unique base-point; if 
$\delta \leq 7$, then $|-K_{X_\delta}|$ is 
base-point free and if $\delta \leq 6$ it is very ample.
\end{thm}
{\it Proof.}  This result is well-known (cf. \cite{Ma}).  \bo

\begin{prop} \label{cane}
Let $X_\delta $ be a del Pezzo surface of degree $9-\delta \geq 3$.  
The scheme $\overline \cM _{bir} \bigl( X_\delta , -K_{X_\delta} \bigr)$ 
is birational to a $\bP^ {6-\delta }-$bundle over $X_\delta $; in 
particular, it is rational and irreducible.
\end{prop}
{\it Proof.}  The surface $X_\delta $ is embedded in $\bP ^{9-\delta}$ by 
the linear system $|-K_{X_\delta}|$.  A general point $[f : \bP^1 \rightarrow X_\delta]$ 
of $\overline \cM _{bir} \bigl( X_\delta , -K_{X_\delta} \bigr)$ corresponds to 
a morphism whose image has a unique singular point $p \in X_\delta $ and is 
uniquely determined by the hyperplane containing $f(\bP^1)$.  Such a hyperplane 
is tangent to $X_\delta $ at $p$.  The hyperplanes in $\bP^{9-\delta}$ intersecting 
$X_\delta $ in a curve with a singular point at $p$ are precisely the hyperplanes 
containing the tangent plane to $X_\delta $ at $p$.

We thus have a rational morphism 
$$\xymatrix {\pi : \overline \cM _{bir} \bigl( X_\delta , -K_{X_\delta} \bigr) \ar@{-->} [r] & X_\delta }$$
assigning to $[f : \bP^1 \rightarrow X_\delta]$ the unique singular point of $f(\bP^1)$.  
The general point of the fiber of $\pi $ over a general point $p \in X_\delta $ corresponds to 
a hyperplane containing the tangent plane to $X_\delta$ at $p$.  The space of such 
hyperplanes is isomorphic to $\bP ^{6-\delta}$.  Since $X_\delta $ is irreducible and the 
general fiber of $\pi$ is also, it follows that $\overline \cM _{bir} \bigl( X_\delta , -K_{X_\delta} \bigr)$ 
is irreducible.  \bo

\noindent
{\it Remark}.  The schemes $\overline \cM _{0,0} \bigl( X_\delta , -K_X \bigr)$ 
are not irreducible if $X_\delta$ is the blow-up of $\bP^2$ at $\delta = 1$ or 
2 points.  Indeed, let $X_1$ be the blow-up of $\bP^2$ at one point $p$; there 
are two morphisms 
$$ \xymatrix { & X_1 \ar[dl] _{\pi _1} \ar[dr] ^{\pi_2} 
{}\save[]+<33pt,2pt>*{\subset \bP^2 \times \bP^1} \restore \\ 
\bP^2 && \bP^1 } $$
and the divisor class of a fiber of $\pi _2$ is $\ell - e$, where 
$\ell $ is the pull-back of the class of a line in $\bP^2$ under 
$\pi_1$ and $e$ is the exceptional fiber of $\pi_1$.  It is clear 
that the space of morphisms from a curve with dual graph 
$$ \xygraph {[] !~:{@{=}} 
!{<0pt,0pt>;<20pt,0pt>:} 
{\bullet} [rr] {\bullet} 
*\cir<2pt>{}
!{\save +<0pt,8pt>*\txt{$\scriptstyle \bar C_2$}  \restore}
- [ll]
*\cir<2pt>{} 
!{\save +<0pt,8pt>*\txt{$\scriptstyle \bar C_1$}  \restore} } $$
where $\bar C_1$ is a (rational) triple cover of a fiber of $\pi_2$ 
and $\bar C_2$ is a double cover of the exceptional fiber of $\pi_1$ 
has dimension at least 7: there are 4 parameters for the triple 
cover of $\ell - e$, 1 for the choice of fiber of $\pi_2$ and 2 for 
the double cover of $e$.

Similarly, let $X_2$ be the blow-up of $\bP^2$ at two distinct 
points $p,q$ and let $\{ \ell , e_1 , e_2 \}$ be a standard 
basis.  It is clear that the space of morphisms from a curve 
with dual graph 
$$ \xygraph {[] !~:{@{=}} 
!{<0pt,0pt>;<20pt,0pt>:} 
{\bullet} [rr] {\bullet} [rr] {\bullet} 
*\cir<2pt>{}
!{\save +<0pt,8pt>*\txt{$\scriptstyle \bar C_2$}  \restore}
- [ll]
*\cir<2pt>{}
!{\save +<0pt,8pt>*\txt{$\scriptstyle \bar D$}  \restore}
- [ll]
*\cir<2pt>{} 
!{\save +<0pt,8pt>*\txt{$\scriptstyle \bar C_1$}  \restore} } $$
where $\bar C_i$ is a double cover of the $(-1)-$curve with divisor 
class $e_i$ and $\bar D$ is a triple cover of the $(-1)-$curve with 
divisor class $\ell - e_1 - e_2$ has dimension at least 8.

In both these cases it is easy to check (Proposition \ref{grafico}) 
that in fact the dimension of the components described is precisely 
the indicated lower bound.

It is also possible to show that these are the only irreducible 
components of $\overline \cM _{0,0} \bigl( X , -K_X \bigr)$ 
besides the closure of $\cM _{0,0} \bigl( X , -K_X \bigr)$, 
when $X$ is a del Pezzo surface.

\begin{prop} \label{cadute}
Let $X$ be a del Pezzo surface of degree two and let $K_X$ be 
the canonical divisor of $X$.  The scheme 
$\overline \cM _{bir} \bigl( X , -K_X \bigr)$ is isomorphic to a 
smooth plane quartic.
\end{prop}
{\it Proof.}  We know (Lemma \ref{pumba}) that there is a morphism 
$\kappa : X \rightarrow \bP^2$ associated to the anticanonical 
sheaf and since $(-K_X)^2 = 2$ (and $-K_X$ is ample), this morphism is 
finite of degree two.  Let $R \subset \bP^2$ be the branch curve, and let 
$2r$ be its degree; denote by $\bar R \subset X$ the ramification 
divisor.  Let $\cO _X (1) = \kappa ^* \cO _{\bP^2} (1) \simeq 
\cO _X (-K_X)$; then using the identity 
$K_X = \kappa ^* K_{\bP^2} + \bar R$, we have 
$\cO _X (-1) \simeq \cO _X (-3+r)$ and we deduce that $r=2$.  
Thus $R$ is a plane quartic.  It is smooth since the morphism 
$\kappa $ has degree two and $X$ is smooth.

The general point of every irreducible component of 
$\overline \cM _{bir} \bigl( X , -K_X \bigr)$ corresponds to 
a singular divisor in $|-K_X|$.  These in turn are parameterized 
by the tangent lines to the ramification curve $R$ of $\kappa $.  
Let $p \in R$ be a point and let $L_p$ be the tangent 
line to $R$ at $p$.  It is easy to convince oneself that by associating to 
each point $p$ in $R$ the morphism which is the normalization of 
$\kappa ^{-1} (L_p)$ at $\kappa ^{-1} (p)$, gives an isomorphism 
$R \simeq \overline \cM _{bir} \bigl( X , -K_X \bigr)$.  \bo

We now deal with the three spaces 
$\overline \cM _{bir} \bigl( X,-nK \bigr)$ for $n \in \{ 1,2,3 \}$, 
where $X$ is a del Pezzo surface of degree one.

\begin{prop}
Let $X$ be a del Pezzo surface of degree one and let $K_X$ be 
the canonical divisor of $X$.  The scheme 
$\overline \cM _{0,0} \bigl( X , -K_X \bigr)$ has dimension zero 
and length twelve.  \bo
\end{prop}

The next two results prove that 
$\overline \cM _{bir} \bigl( X,-2K \bigr)$ and 
$\overline \cM _{bir} \bigl( X,-3K \bigr)$ are irreducible assuming 
that the del Pezzo surface $X$ is general.

\begin{thm} \label{cabala}
Let $X$ be a general del Pezzo surface of degree one and let $C$ 
be the closure of the set of points of $|-2K_X|$ corresponding 
to reduced curves whose normalization is irreducible and of 
genus zero.  Then $C$ is a smooth irreducible curve.
\end{thm}
{\it Proof.}  The linear system associated to the line bundle $\cO_X(-2K_X)$ 
on the del Pezzo surface of degree one is base-point free and determines 
a finite morphism $\kappa $ of degree two to $\bP^3$, whose image is a 
quadric cone $Q$.  The image $R$ of the ramification divisor of $\kappa $ 
is a smooth canonically embedded curve of genus four which does not 
contain the vertex of the cone.  Clearly, the vertex $v$ of $Q$ is an 
isolated ramification point, since $X$ is smooth.

Let $C \subset |-2K_X| \simeq \bP^3$ be the closure of the set of all 
points corresponding to integral curves whose normalization is 
irreducible and of genus zero.  In order to prove that $C$ is smooth and 
irreducible, we will first prove it is connected, and then that it is 
smooth.

Since the arithmetic genus of a divisor $D$ in $|-2K_X|$ is two, in order 
for $D$ to be integral and have geometric genus zero, $D$ must be tangent 
to the ramification divisor $\bar R$ at two points.  If we translate this 
in terms of the image of the morphism $\kappa $, this implies that the plane 
$P$ corresponding to the divisor $D$ intersects $R$ along a divisor of the 
form $2(p) + 2(q) + (r) + (s)$, for some points $p,q,r,s \in R$.  The 
condition that $D$ should be integral translates to the requirement that 
the plane $P$ should not contain a line of $Q$.  If this happens, then we 
have $P \cap R = 2 \bigl( (p)+ (q) + (r) \bigr)$ and 
$2 \bigl( (p) + (q) + (r) \bigr)$ is the (scheme-theoretic) fiber of the 
projection $p|_R$ away from the vertex $v$ (alternatively, $(p) + (q) + (r)$ 
is the scheme-theoretic intersection of a line on $Q$ with $R$).

Consider the smooth surface $R \times R$ and the two projection morphisms 
$$ \xymatrix { & R \times R \ar[dl] _{\pi _1} \ar[dr] ^{\pi _2} \\
R && R } $$
where $\pi _1$ is the projection onto the first factor and $\pi _2$ onto 
the second.  Denote by $\Delta \subset R\times R$ the diagonal.  Let 
$\cF := \cO _{R \times R} \bigl( \pi _2 ^* K_R - 2 \Delta \bigr)$ be a sheaf 
on $R \times R$ and let $\cE := (\pi _1) _* \cF $ be a sheaf on $R$.

Clearly $\cF $ is invertible.  The sheaf $\cE $ is locally free of rank two.  
To prove this, we compute for any $p \in R$ 
$$ h^0 \bigl( p , \cF \bigr) := 
\dim {\rm H}^0 \Bigl( (\pi _1) ^{-1} (p) , \cF |_{(\pi _1) ^{-1} (p)} \Bigr) 
= \dim {\rm H}^0 \Bigl( R , \cO_R \bigl( K_R - 2(p) \bigr) \Bigr) $$
We know that the last dimension is at least two, since there is a pencil of 
planes in $\bP^3$ containing the tangent line to $R$ at $p$.  By Riemann-Roch 
it follows that the sheaf $\cO_R \bigl( K_R - 2(p) \bigr)$ has non-vanishing 
first cohomology group.  By Clifford's Theorem (\cite{Ha} Theorem IV.5.4) the 
dimension of ${\rm H}^0 \Bigl( R , \cO_R \bigl( K_R - 2(p) \bigr) \Bigr)$ is 
at most 3 and since $R$ is not hyperelliptic (because it is a canonical 
curve) and obviously $K_R - 2(p)$ is not 0 nor $K_R$, it follows that 
$h^0 \bigl( p , \cF \bigr) = 2$ for all $p \in R$.

We may now apply the first part of Grauert's Theorem (\cite{Ha} Corollary 
III.12.9) to conclude that $\cE = (\pi _1) _* \cF $ is locally free and 
the second part of the same theorem to conclude that the natural morphism 
of sheaves on $R \times R$ 
$$ \xymatrix { \pi _1 ^* \cE = \pi _1 ^* \bigl( (\pi _1 )_* \cF \bigr) 
\ar[r] & \cF } $$
is surjective.  In turn, this implies (\cite{Ha} Proposition II.7.12) that 
there is a commutative diagram 
$$ \xymatrix { R \times R \ar[r] ^{\varphi } \ar[d] _{\pi _1} & 
\bP \bigl( \cE \bigr) \ar[d] ^{\pi } \\ 
R  \ar[r] ^{id} & R } $$

The morphism $\varphi $ is finite of degree four.  Let 
$\bar C \subset R \times R$ be the ramification divisor of $\varphi $ 
and let $F \subset R \times R$ be the closure of the set of points 
$\bigl\{ (p,q) ~ \bigl| ~ p_v (p) = p_v (q) ~,~ p \neq q \bigr\}$ (remember 
that $p_v$ is the projection away from the cone vertex $v$ of $Q$).  
Note that $\bar C$ does not contain any fiber of $\pi $, since all the 
induced morphisms 
$\varphi _p : R_p := (\pi _1)^{-1} (p) \rightarrow \bP^1 _p := \pi ^{-1} (p)$ 
are ramified covers of degree 4.  Moreover we have 
$R_p \cdot \bar C = 14$, since for all $p$ such intersection represents the 
ramification divisor of the morphism $\varphi _p$ which has degree four, and 
we may therefore compute the intersection using the Hurwitz formula.

By definition, $\bar C$ is the set of pairs $(p,q)$ such that if we denote 
by $P _p ^q$ the plane containing $q$ and the tangent line to $R$ at $p$ 
(or the osculating plane to $R$ at $p$, if $p=q$), then we have 
$P_p ^q \cap R \geq 2\bigl( (p) + (q) \bigr)$.

We clearly have $F \subset \bar C$, since if $(p,q) \in F$ then 
$P_p ^q $ is in fact the tangent plane to $Q$ at $p$ and thus 
$P_p ^q \cap R = 2 \bigl( (p) + (q) + (r) \bigr) \geq 2 \bigl( (p) + (q) \bigr)$.

By definition we have $C \subset \bar C$ and no component of $C$ 
is also a component of $F$, since the plane corresponding to a point in 
$F$ intersects $Q$ in a non-reduced curve.  It is also immediate to 
check that in fact $C$ is the residual curve to $F$ in $\bar C$, that is 
we have $\bar C = C \cup F$.

We now prove that the residual curve $C$ to $F$ in 
$\bar C$ is connected and (for general $X$) smooth.

The connectedness of $C$ is a consequence of a theorem of Kouvidakis: 
the divisor class of $C$ in $R\times R$ is $4 (F_1 + F_2) - \Delta $, 
where $\Delta $ is the diagonal and the $F_i$'s are the fibers of the 
two projections to $R$.  Thanks to \cite{La} Example 1.5.13, we 
know that $C$ is an ample divisor.  In particular, $C$ is connected.

To prove the smoothness of $C$, we will show that for any point 
$(p,q) \in C$ the two numbers 
$mult _{(p,q)} \Bigl( \bigl( \pi _1 \bigl|_C \bigr) ^{-1} \bigl( p \bigr) \Bigr)$ 
and 
$mult _{(p,q)} \Bigl( \bigl( \pi _2 \bigl|_C \bigr) ^{-1} \bigl( q \bigr) \Bigr)$ 
cannot both be at least two.  Since this would be the case if $(p,q)$ 
were a singular point, the theorem follows.

Let $p \in R$ and let $(p')+(p'')$ be the divisor obtained by intersecting 
the curve $R$ with the line on $Q$ through $p$; we have 
\begin{eqnarray*}
\Bigl( \pi _1 \bigl|_{\bar C} \Bigr) ^{-1} \bigl( p \bigr) & \hspace{-5pt} = & \hspace{-7pt}
\sum _{q \in R_p} \Bigl( mult _q \bigl( P_p ^q \cap R \bigr) - 1 \Bigr) 
\bigl( p , q \bigr) - 2 \bigl( p , p \bigr) \\
\Bigl( \pi _1 \bigl|_F \Bigr) ^{-1} \bigl( p \bigr) & \hspace{-5pt} = & \hspace{-7pt}
\bigl( p , p' \bigr) + \bigl( p , p'' \bigr) \\
\Bigl( \pi _1 \bigl|_C \Bigr) ^{-1} \bigl( p \bigr) & \hspace{-5pt} = & \hspace{-7pt}
\sum _{q \in R_p} \Bigl( mult _q \bigl( P_p ^q \cap R \bigr) - 1 \Bigr) \bigl( p , q \bigr) 
- \bigl( p , p' \bigr) - \bigl( p , p'' \bigr) - 2 \bigl( p , p \bigr) 
\end{eqnarray*}
and thus we deduce that 
\begin{eqnarray*}
Ram \Bigl( \pi _1 \bigl|_C \Bigr) & = & 
Ram \Bigl( \pi _1 \bigl|_{\bar C} \Bigr) 
- Ram \Bigl( \pi _1 \bigl|_F \Bigr) = \\[7pt]
& = & \sum _{p \in R} \left( \sum _{q \in R_p \cap \bar C} 
\Bigl( mult _q \bigl( P_p ^q \cap R \bigr) - 2 \Bigr) \bigl( p,q \bigr) \right) - \\[7pt]
& & ~~~- 
\Bigl( \bigl( p_1 , p_1' \bigr) + \ldots + \bigl( p_{12} , p_{12}' \bigr) \Bigr) 
\end{eqnarray*}
where $R_{p_i} \cap F = 2 \bigl( p_i , p_i' \bigr) + \bigl( p_i , p_i \bigr)$ 
(equivalently, the line $L_i$ on $Q$ containing $p_i$ is tangent to the image of 
the ramification divisor of $\kappa $ at $p_i' \neq p_i$).

We conclude that $\bigl( p , q \bigr) \in C$ is a ramification point for 
$\pi _1 \bigl| _C$ if and only if $P_p^q \cap R = 2 (p) + 3(q) + (r)$, for 
some $r \in R$.  In view of this asymmetry between $p$ and $q$, we deduce that 
$\bigl( p , q \bigr)$ can be a ramification for both projections 
$\pi _1 \bigl| _C$ and $\pi _2 \bigl| _C$ if and only if 
$P_p^q \cap R = 3 (p) + 3(q)$.  If $p$ and $q$ are on the same line on $Q$, 
then the inverse image under $\kappa $ of that line would be a cuspidal 
divisor in $|-K_X|$, which we are excluding.  We will now prove that the 
dimension of the space of smooth canonically embedded curves $R$ of 
arithmetic genus four lying on a singular quadric and having a plane $P$ 
transverse to the quadric cone and intersecting $R$ along a divisor of the 
form $3 \bigl( (p) + (q) \bigr)$ is at most seven, and thus for the general 
del Pezzo surface of degree one, this configuration does not happen.  This 
will conclude the proof.

This is simply a dimension count: using automorphisms of $\bP^3$ we may 
assume that the plane $P$ has equation $X_3 = 0$ and that the quadric 
cone has equations $X_0 X_1 = X_2 ^2$.  We may also assume that $p$ and 
$q$ have coordinates $[1,0,0,0]$ and $[0,1,0,0]$ respectively.  Note that 
we still have a two-dimensional group of automorphisms (with one generator 
corresponding to rescaling the coordinate $X_3$, and the other corresponding 
to multiplying the coordinate $X_0$ by a non-zero scalar and the coordinate 
$X_1$ by its inverse).  With these choices, the quadric cone is completely 
determined, as well as the plane $P$.  We still need to compute how many 
parameters are accounted for by the cubic intersecting the cone in $R$.

For this last computation, we consider the short exact sequences of sheaves 
$$ \xymatrix @R=10pt { 
0 \ar[r] & \cO_{\bP^3} \ar[r] & \cO_{\bP^3} (2) \ar[r] & \cO_Q (2) \ar[r] & 0 \\
0 \ar[r] & \cO_Q (2) \ar[r] & \cO_Q (3) \ar[r] & \cO_R (3) \ar[r] & 0 } $$

The first sequence implies that the cohomology groups 
${\rm H} ^i \bigl( Q , \cO_Q (2) \bigr)$ are zero for $i \geq 1$; therefore, 
from the second sequence we deduce that the dimension of the space of cubics 
vanishing on $R$ is nine.  Subtracting the two-dimensional automorphism group 
leaves us with a family of dimension seven.  Since there is a family of dimension 
eight of del Pezzo surfaces of degree one, we conclude.  \bo

In order to prove a similar result for the divisor class $-3K_X$, we first 
establish a lemma.

\begin{lem} \label{pizzica}
Let $X$ be a smooth projective surface and let $K_1$, $K_2$ and 
$K_3$ be three distinct nodal rational divisors of anticanonical degree one meeting 
at a point $p \in X$.  Suppose that two of the components meet transversely at $p$.  
Let $f : \bar C := \bar K_1 \cup \bar K_2 \cup \bar K_3 \cup \bar E 
\longrightarrow X$ be the stable map of genus zero, such that 
\begin{itemize}
\item the morphism $f_i := f|_{\bar K_i}$ is the normalization of $K_i$ followed 
by the inclusion in $X$;
\item the component $\bar E$ is contracted to the point $p \in K_1 \cap K_2 \cap K_3$;
\item the dual graph of the morphism $f$ is
$$ \xygraph {[] !~:{@{=}} 
!{<0pt,0pt>;<20pt,0pt>:} 
{\bullet} [rr] {\bullet} [ur] {\bullet} [dd] {\bullet}
*\cir<2pt>{}
!{\save +<8pt,0pt>*\txt{$\scriptstyle \bar K_3$}  \restore}
- [ul]
*\cir<2pt>{}
!{\save +<-1pt,8pt>*\txt{$\scriptstyle \bar E$}  \restore}
- [ur]
*\cir<2pt>{} 
!{\save +<8pt,0pt>*\txt{$\scriptstyle \bar K_2$}  \restore}
 [dl] - [ll]
*\cir<2pt>{}
!{\save +<0pt,8pt>*\txt{$\scriptstyle \bar K_1$}  \restore} } $$
\vglue0pt{ \centerline{Dual graph of $f$}}
\end{itemize}

Then the point represented by the morphism $f$ in 
$\overline { \cM } _{0,0} \bigl( X , -3K_X \bigr)$ lies in a unique 
irreducible component.
\end{lem}
{\it Proof.}  The expected dimension of $\overline { \cM } _{0,0} \bigl( X , -3K_X \bigr)$ 
is $-K_X \cdot (K_1 + K_2 + K_3) - 1 = 2$.  The first step of the proof 
consists of proving that the embedding dimension of 
$\overline { \cM } _{0,0} \bigl( X , K_1 + K_2 + K_3 \bigr)$ 
at $[f]$ is at most three.  To prove this, it suffices to prove that 
${\rm H}^1 \bigl( \bar C , f^* \cT_X \bigr)$ is one-dimensional.  This in 
turn will follow from the fact that 
${\rm H}^0 \bigl( \bar C , f^* \cT_X \bigr)$ has dimension six.  On each 
irreducible component $\bar K_i$ we have 
$f_i ^* \cT_X \simeq \cO_{\bar K_i} (2) \oplus \cO_{\bar K_i} (-1)$, where the 
$\cO_{\bar K_i} (2)$ summand is the tangent sheaf of $\bar K_i$.  Denote by 
$f_E$ the restriction of $f$ to the component $\bar E$; we have 
$f_E ^* \cT_X \simeq \cO_{\bar E} \oplus \cO_{\bar E}$.  Consider the sequence 
\begin{equation} \label{paganini}
\xymatrix @C=15pt { 0 \ar[r] & f^* \cT _X \ar[r] &  
f_1^* \cT _X \oplus f_2^* \cT _X \oplus f_3^* \cT _X \oplus f_E^* \cT _X \ar[r] & 
\cT _{X,p} \oplus \cT _{X,p} \oplus \cT _{X,p} \ar[r] & 0 }
\end{equation}

Since two if the $K_i$'s meet transversely at $p$, it follows that in order for 
the global sections on the irreducible components of $\bar C$ to glue together, it 
is necessary that the sections on the $\bar K_i$'s vanish at the node with $\bar E$.  
Moreover, if such a condition is satisfied, clearly the sections on the components 
$\bar K_i$ together with the zero section on $\bar E$ glue to give a global section 
of $f^* \cT_X$.  We deduce that ${\rm H}^0 \bigl( \bar C , f^* \cT_X \bigr)$ has 
dimension six, and it is isomorphic to 
${\rm H} ^0 \bigl( \bar K_1 , \cT_{\bar K_1} (-\bar p_1) \bigr) \oplus 
 {\rm H} ^0 \bigl( \bar K_2 , \cT_{\bar K_2} (-\bar p_2) \bigr) \oplus 
 {\rm H} ^0 \bigl( \bar K_3 , \cT_{\bar K_3} (-\bar p_3) \bigr) \bigr)$, where 
$\bar p_i \in \bar K_i$ is the node with $\bar E$.  From the exact sequence 
(\ref{paganini}) and the fact that 
${\rm H} ^1 \bigl( \bar C , f_1^* \cT_X \oplus f_2^* \cT_X \oplus f_3^* \cT_X \oplus 
f_E^* \cT_X \bigr) = 0$, we deduce that 
$$ h^1 \bigl( \bar C , f^* \cT _X \bigr) = 6 - \chi \bigl( \bar C , f^* \cT _X \bigr) = 
6 - (3+3+3+2) + 6 = 1 $$

Thus the embedding dimension of $\overline { \cM } _{0,0} \bigl( X , K_1 + K_2 + K_3 \bigr)$ 
at $[f]$ is at most three, as stated above.  It follows that we may write 
$$ \hat {\cO} _{[f]} \overline { \cM } _{0,0} \bigl( X , K_1 + K_2 + K_3 \bigr) 
\simeq k [\![ t_1 , t_2 , t_3 ]\!] / (g) $$

We thus deduce that all the components through $[f]$ have dimension equal to two, 
since there is a component of dimension two through $[f]$ and if there were also 
a component of dimension three or more containing $[f]$, then the embedding 
dimension would be more than three.  
Moreover, if there are two components containing $[f]$, then the singular points 
of $\overline { \cM } _{0,0} \bigl( X , K_1 + K_2 + K_3 \bigr)$ near $[f]$ must 
have dimension equal to one.  We prove that $[f]$ is an isolated singular point, 
and thus we conclude that there is a unique component containing $[f]$.

Let $U \subset \overline { \cM } _{0,0} \bigl( X , K_1 + K_2 + K_3 \bigr)$ be the 
open subset of morphisms $g : \bar D \rightarrow X$ which are immersions and 
birational to their image.

The subset $U$ is contained in the smooth locus of 
$\overline { \cM } _{0,0} \bigl( X , K_1 + K_2 + K_3 \bigr)$ thanks to Proposition 
\ref{grafico}.  Moreover $U \cup \bigl\{ [f] \bigr\}$ is a neighbourhood of $[f]$: 
all the morphisms in a neighbourhood of $[f]$ must have image consisting of at most 
two components, since the morphisms $f_i$ have no infinitesimal deformations.  It 
follows that there are neighbourhoods of $[f]$ such that $[f]$ is the only morphism 
with a contracted component.  Since the image of $f$ has no cusps and any two 
components meet transversely, the same statement holds for all the morphisms in a 
neighbourhood of $[f]$.  It follows that $U \cup \bigl\{ [f] \bigr\}$ is a 
neighbourhood of $[f]$.  Thus $[f]$ is an isolated singular point (possibly a smooth 
point) and since the embedding dimension of 
$\overline { \cM } _{0,0} \bigl( X , K_1 + K_2 + K_3 \bigr)$ at $[f]$ is at most three 
it follows that $\overline { \cM } _{0,0} \bigl( X , K_1 + K_2 + K_3 \bigr)$ is 
locally irreducible near $[f]$, thus concluding the proof of the lemma.  \bo

\begin{thm} \label{maschera}
Let $X$ be a del Pezzo surface of degree one such that the space 
$\overline { \cM } _{bir} \bigl( X , -2K_X \bigr)$ is irreducible 
and all the rational divisors in $|-K_X|$ are nodal.  Let $S$ be 
the closure of the set of points of $|-3K_X|$ corresponding to 
reduced curves whose normalization is irreducible and of genus 
zero.  Then $S$ is an irreducible surface.
\end{thm}
{\it Proof.}  Let $f : \bP^1 \rightarrow X$ be a morphism in 
$\overline { \cM } _{bir} \bigl( X , -3K_X \bigr)$.  Thanks to Proposition 
\ref{rpunti} and Lemma \ref{immersione}, we may assume that $f$ is an 
immersion and that its image contains a general point $p$ of $X$.  In 
particular it follows that $[f]$ represents a smooth point of 
$\overline { \cM } _{bir} \bigl( X , -3K_X \bigr)$.  We choose the 
point $p$ to be an independent point (Definition \ref{indip}).

Consider the space of morphisms of 
$\overline { \cM } _{bir} \bigl( X , -3K_X \bigr)$ in the same irreducible 
component as $[f]$ which contain the point $p$ in their image, denote 
this space by $\overline \cM _{bir} (p)$.  It follows immediately from 
the dimension estimates (\ref{dimdibarbi}) that 
$\dim _{[f]} \overline \cM _{bir} (p) = 1$ and that $[f]$ is a smooth 
point of $\overline \cM _{bir} (p)$.  We may therefore find a smooth 
irreducible projective curve $B$, a normal surface 
$\pi : S \rightarrow B$ and a morphism $F : S \rightarrow X$ such that 
the induced morphism $B \rightarrow \overline \cM _{bir} (p)$ is 
surjective onto the component containing $[f]$.  From \cite{Ko} 
Corollary II.3.5.4, it follows immediately that the morphism $F$ is 
dominant.  We want to show that there are fibers of $\pi$ that are 
reducible.  This is clear, since the morphism 
$F^* : {\rm Pic} (X) \rightarrow {\rm Num} (S)$ is injective, and 
${\rm Pic} (X)$ has rank nine, while if every fiber of $\pi$ were a smooth 
rational curve, it would follow that $\pi : S \rightarrow B$ is a ruled surface 
(\cite{Ha} V.2) and thus that the rank of ${\rm Num} (S)$ is two.

This implies that there must be a morphism $f_0 : \bar C \rightarrow X$ 
with reducible domain in the family of 
stable maps parametrized by $B$, and since all such morphisms contain the 
general point $p$ in their image, the same is true of the morphism $f_0$.  
In particular, since the point $p$ does not lie on any rational curve of 
anticanonical degree 1, it follows that $\bar C$ consists of exactly two 
components $\bar C_1$ and $\bar C_2$, where each $\bar C_i$ is irreducible 
and we may assume that $f_0 (\bar C_1)$ has anticanonical degree one and 
$f_0 (\bar C_2)$ has anticanonical degree two.  Denote by $C_i$ the image of 
$\bar C_i$.  It also follows from the definition of an independent point and 
Proposition \ref{grafico} that $f_0$ represents a smooth point of 
$\overline { \cM } _{bir} \bigl( X , -3K_X \bigr)$.

There are two possibilities for the divisor $C_1$: it is either a 
$(-1)-$curve (there are 240 such divisors on $X$), or it is rational curve 
in the anticanonical divisor class (there are 12 such divisors on $X$).  
We will prove that we may assume that $C_1$ is a rational divisor in the 
anticanonical linear system.

Suppose that $C_1$ is a $(-1)-$curve and let $C_1' \subset X$ be the 
$(-1)-$curve such that $C_1 + C_1' = -2K_X$.  The curve $C_2$ is thus an 
integral curve in the linear system $-3K_X - C_1 = -K_X - C_1'$.  It 
follows that $C_2$ is in the anticanonical linear system on the del Pezzo 
surface of degree two obtained by contracting $C_1'$.

The morphism $\varphi : X \rightarrow \bP^2$ 
associated to the divisor $C_2$ is the contraction of the $(-1)-$curve $C_1'$ 
followed by the degree two morphism to $\bP^2$ induced by the anticanonical 
divisor on the resulting surface $X'$.  In the plane $\bP^2$ we therefore have 
\begin{itemize}
\item the image of the ramification curve $R$, which is a smooth plane quartic;
\item the image of $C_1'$, which is a point $q$;
\item the image of $C_1$, which is a plane quartic with a triple point at $q$ 
and is everywhere tangent to $\varphi (R)$;
\item the image of $C_2$, which is a tangent line to $\varphi (R)$.
\end{itemize}
To be precise, the ramification divisor of $\varphi$ consists of two disjoint 
components, one is the $(-1)-$curve $C_1'$, whose image is the point $q$, and 
the other is a curve whose image is a smooth plane quartic.

Consider the morphism 
$$ \xymatrix { {\rm Sl}_{f_0} (\bar C_2) \ar[r] ^{\hspace{10pt}a} & \bar C_1 } $$
and let $\bar p \in \bar C_1$ be one of the (three) points mapping to the 
intersection $C_1 \cap C_1'$ (and in particular, $\varphi \bigl( f_0 (\bar p) \bigr) = q$).  
Let $f_1$ be a morphism in the fiber of $a$ above the point $\bar p$.  The image of 
$f_1$ consists of the divisor $\varphi (C_1)$ together with one of the tangent lines 
$L$ to $\varphi (R)$ containing the point $q$.

The domain curve of $f_1$ consists of possibly a contracted component and three 
more non-contracted components $\bar C_1$ mapped to $C_1$, $\bar L$ mapped to 
the closure of $\varphi ^{-1} \bigl( L \bigr) \setminus C_1'$ and finally 
$\bar C_1'$ mapped to the $(-1)-$curve $C_1'$.  The possible dual graphs of 
$f_1$ are 
$$ \xygraph {[] !~:{@{=}} 
!{<0pt,0pt>;<20pt,0pt>:} 
{\bullet} [rr] {\bullet} [rr] {\bullet} 
*\cir<2pt>{}
!{\save +<0pt,8pt>*\txt{$\scriptstyle \bar C_1$}  \restore}
- [ll]
*\cir<2pt>{}
!{\save +<0pt,8pt>*\txt{$\scriptstyle \bar C_1'$}  \restore}
- [ll]
*\cir<2pt>{} 
!{\save +<0pt,8pt>*\txt{$\scriptstyle \bar L$}  \restore} } 
\hspace{20pt}
\xygraph {[] !~:{@{=}} 
!{<0pt,0pt>;<20pt,0pt>:} 
{\bullet} [rr] {\bullet} [ur] {\bullet} [dd] {\bullet}
*\cir<2pt>{}
!{\save +<8pt,0pt>*\txt{$\scriptstyle \bar C_1'$}  \restore}
- [ul] -[ur] 
*\cir<2pt>{}
!{\save +<0pt,8pt>*\txt{$\scriptstyle \bar C_1$}  \restore}
 [dl] 
*\cir<2pt>{}
!{\save +<0pt,8pt>*\txt{$\scriptstyle \bar E$}  \restore}
- [ll]
*\cir<2pt>{} 
!{\save +<0pt,8pt>*\txt{$\scriptstyle \bar L$}  \restore} } $$
\vglue0pt{ \centerline{Possible dual graphs of $f_1$} \vspace{5pt}}

Note that in the first case $f_1$ represents a smooth point of 
the space $\overline \cM _{0,0} \bigl( X,-3K_X \bigr)$; in the second case, we 
may apply Lemma \ref{pizzica} to conclude that even if $[f_1]$ is not a 
smooth point, deforming it produces morphisms in the same irreducible 
component as $f_1$.  Smoothing out the components $\bar C_1 \cup \bar C_1'$ 
(or $\bar C_1 \cup \bar E \cup \bar C_1'$ if there is a contracted component) 
we obtain a morphism which has one component (the one obtained by smoothing) 
mapped birationally to a rational curve in $|-2K_X|$ and another component 
(the component $\bar L$, with notation as above) mapped birationally to a 
rational divisor in $|-K_X|$.

Thus we may deform the original morphism $f$ to a morphism 
$f_0 : \bar C_1 \cup \bar C_2 \rightarrow X$ such that $\bar C_1$ is mapped 
birationally to a rational curve in the anticanonical linear system and 
$\bar C_2$ is mapped birationally to a rational curve in $|-2K_X|$.

Choose three nodal rational curves $K_1$, $K_2$ and $K_3$ in the linear 
system $-K_X$.  We prove now that we may deform $f_0$ without changing 
the irreducible component of $\overline \cM _{0,0} \bigl( X,-3K_X \bigr)$ 
to the morphism 
$g : \bar K_1 \cup \bar K_2 \cup \bar K_3 \cup \bar E \longrightarrow X$ 
such that $\bar K_i$ is the normalization of $K_i$, $\bar E$ is contracted to 
the point in the intersection $K_1 \cap K_2 \cap K_3$ and the dual graph 
of $g$ is 
$$ \xygraph {[] !~:{@{=}} 
!{<0pt,0pt>;<20pt,0pt>:} 
{\bullet} [rr] {\bullet} [ur] {\bullet} [dd] {\bullet}
*\cir<2pt>{}
!{\save +<8pt,0pt>*\txt{$\scriptstyle \bar K_3$}  \restore}
- [ul]
*\cir<2pt>{}
!{\save +<-1pt,8pt>*\txt{$\scriptstyle \bar E$}  \restore}
- [ur]
*\cir<2pt>{} 
!{\save +<8pt,0pt>*\txt{$\scriptstyle \bar K_2$}  \restore}
 [dl] - [ll]
*\cir<2pt>{}
!{\save +<0pt,8pt>*\txt{$\scriptstyle \bar K_1$}  \restore} } $$
\vglue0pt{ \centerline{Dual graph of $g$} \vspace{5pt}}
\noindent
It follows from this and Lemma \ref{pizzica} that 
$\overline \cM _{bir} \bigl( X,-3K_X \bigr)$ is irreducible.

To achieve the required deformation, we consider the morphism 
$$ \xymatrix @C=40pt { {\rm Sl}_{f_0} (\bar C_2) \ar[r] ^{\pi \hspace{25pt}} & 
\overline \cM _{bir} \bigl( X,-2K_X \bigr) } $$
and note that $\pi $ is surjective since 
$\overline \cM _{bir} \bigl( X,-2K_X \bigr)$ is irreducible by assumption.

Relabeling $K_1$, $K_2$ and $K_3$, we may suppose that $C_1 \neq K_2, K_3$.  
Thus we may specialize $f_0$ to a morphism 
$f_1 : \bar C_1 \cup \bar K_2 \cup \bar K_3 \cup \bar E \longrightarrow X$ 
such that $f_1 (\bar K_i) = K_i$, $\bar E$ is contracted by $f_1$ and 
the dual graph of $f_1$ is 
$$ \xygraph {[] !~:{@{=}} 
!{<0pt,0pt>;<20pt,0pt>:} 
{\bullet} [rr] {\bullet} [ur] {\bullet} [dd] {\bullet}
*\cir<2pt>{}
!{\save +<8pt,0pt>*\txt{$\scriptstyle \bar K_3$}  \restore}
- [ul]
*\cir<2pt>{}
!{\save +<-1pt,8pt>*\txt{$\scriptstyle \bar E$}  \restore}
- [ur]
*\cir<2pt>{} 
!{\save +<8pt,0pt>*\txt{$\scriptstyle \bar K_2$}  \restore}
 [dl] - [ll]
*\cir<2pt>{}
!{\save +<0pt,8pt>*\txt{$\scriptstyle \bar C_1$}  \restore} } $$
\vglue0pt{ \centerline{Dual graph of $f_1$} \vspace{5pt}}
\noindent
Thanks to Lemma \ref{pizzica} any deformation of such morphism is in the same 
irreducible component of $\overline \cM _{0,0} \bigl( X,-3K_X \bigr)$ as 
$f_0$ and hence in the same irreducible component as the morphism $f$.

We may now smooth the components $\bar C_1 \cup \bar K_2 \cup \bar E$ to 
a single irreducible component mapped birationally to the divisor class 
$-2K_X$ and then we may use irreducibility of 
$\overline \cM _{bir} \bigl( X,-2K_X \bigr)$ again to prove that we may 
specialize the component thus obtained to break as $\bar K_1 \cup \bar K_2$.  
The morphism $g$ thus obtained is the one we were looking for, and the theorem 
is proved.  \bo

\noindent
{\it Remark.}  Thanks to Theorem \ref{cabala}, the space 
$\overline \cM _{bir} \bigl( X,-2K_X \bigr)$ is irreducible for the 
general del Pezzo surface of degree one.  Thus it follows from Theorem 
\ref{maschera} that also the space 
$\overline \cM _{bir} \bigl( X,-3K_X \bigr)$ is irreducible for the 
general del Pezzo surface of degree one.

\subsection{The Picard Group and the Orbits of the Weyl Group} \label{sporco}


In this section we prove some results on the divisor classes 
of the blow-up of $\bP^2$ at eight or fewer general points.  In 
particular we analyze several questions regarding the divisor 
classes of the conics and their orbits under the Weyl group.

Let $X_\delta $ be the blow-up of $\bP^2$ at $\delta \leq 8$ points 
such that no three are on a line, no six of them are on a conic and 
there is no cubic through seven of them with a node at the eighth.

\begin{defi}
A divisor $C$ on $X_\delta $ is called a {\rm conic} if 
$-K_{X_\delta } \cdot C = 2$ and $C^2 = 0$.
\end{defi}

Suppose that $\{ \ell , e_1 , \ldots , e_\delta \}$ is a standard basis of 
${\rm Pic} (X_\delta )$.  If 
$C = a\ell - b_1 e_1 - \ldots - b_\delta e_\delta $ is a divisor class 
on $X_\delta$, then to simplify the notation we simply write it as 
$(a \,;\, b_1 , \ldots , b_\delta )$.

\begin{prop}
The conics on $X_8$ are given, up to permutation of the $e_i$'s, 
by the following table:
\begin{equation} \label{soluco} 
\begin{array} {@{ \vline~~ }c@{ ~~\vline~~ }c@{ ~~\vline~~ }c@{ ~~\vline~~ }c@{ ~~\vline~~ }
c@{ ~~\vline~~ }c@{ ~~\vline~~ }c@{ ~~\vline~~ }c@{ ~~\vline~~ }c@{ ~~\vline~~ }c@{ ~~\vline }}
\hline Type & \ell & e_1 & e_2 & e_3 & e_4 & e_5 & e_6 & e_7 & e_8 \vphantom{\Bigl|} \\[-1pt]
\hline \vphantom{\Bigl|} 
A & 1 & 1 & 0 & 0 & 0 & 0 & 0 & 0 & 0 \\[-3pt]
B & 2 & 1 & 1 & 1 & 1 & 0 & 0 & 0 & 0 \\[-3pt]
C & 3 & 2 & 1 & 1 & 1 & 1 & 1 & 0 & 0 \\[-3pt]
D & 4 & 2 & 2 & 2 & 1 & 1 & 1 & 1 & 0 \\[-3pt]
E & 5 & 2 & 2 & 2 & 2 & 2 & 2 & 1 & 0 \\[-3pt]
D'& 4 & 3 & 1 & 1 & 1 & 1 & 1 & 1 & 1 \\[-3pt]
F & 5 & 3 & 2 & 2 & 2 & 1 & 1 & 1 & 1 \\[-3pt]
G & 6 & 3 & 3 & 2 & 2 & 2 & 2 & 1 & 1 \\[-3pt]
H & 7 & 3 & 3 & 3 & 3 & 2 & 2 & 2 & 1 \\[-3pt]
H'& 7 & 4 & 3 & 2 & 2 & 2 & 2 & 2 & 2 \\[-3pt]
I & 8 & 4 & 3 & 3 & 3 & 3 & 2 & 2 & 2 \\[-3pt]
I'& 8 & 3 & 3 & 3 & 3 & 3 & 3 & 3 & 1 \\[-3pt]
J & 9 & 4 & 4 & 3 & 3 & 3 & 3 & 3 & 2 \\[-3pt]
K & 10& 4 & 4 & 4 & 4 & 3 & 3 & 3 & 3 \\[-3pt]
L & 11& 4 & 4 & 4 & 4 & 4 & 4 & 4 & 3 \\ \hline 
\end{array} 
\end{equation}
Their numbers are given by the table:
$$ \begin{array} {|c|c|c|c|c|c|c|c|}
\hline \delta &    8 &   7 &  6 &  5 & 4 & 3 & 2 \\ \hline 
conics & 2160 & 126 & 27 & 10 & 5 & 3 & 2 \\ \hline 
\end{array} $$
\end{prop}
{\it Proof.}  We proceed just like in \cite{Ma} IV, \S 25.  
The condition of being a conic translates to the equations 
$$ \left\{ \begin{array} {rcl} \displaystyle 
a^2 - \sum _{i=1} ^8 b_i ^2 & = & 0 \\[15pt] \displaystyle 
3a - \sum _{i=1} ^8 b_i & = & 2
\end{array} \right. $$
and we may equivalently rewrite these as 
$$ \left\{ \begin{array} {rcl} \displaystyle 
\sum _{i=1} ^8 \bigl( a - 2 b_i - 2 \bigr) ^2 & = & 16 \\[15pt] \displaystyle 
3a - \sum _{i=1} ^8 b_i & = & 2
\end{array} \right. $$
It is now easy (but somewhat long) to check that (\ref{soluco}) is the 
complete list of solutions up to permutations.  \bo

\noindent
{\it Remark}.  The classes of conics on $X_\delta $ for $\delta \leq 7$ 
are obtained from the ones in list (\ref{soluco}) by erasing $8-\delta$ 
zeros and permuting the remaining coordinates.  Thus (up to permutations) 
the first five rows and seven columns describe conics on $X_7$, the 
first three rows and six columns are the conics on $X_6$ and so on.

We introduce the following notation (which luckily won't be extremely 
useful, but allows us to name conics!) for the classes of the conics 
on $X_\delta $, $\delta \leq 8$ (we set also $\bar E := e_1 + \ldots + e_8$): 
\begin{equation} \label{nomico}
\hspace{-11.3pt}
\left\{ \begin{array} {l@{ \,=\, }l}
A_i & \ell - e_i \\ \vphantom{\vdots}
B_{ijkl} & 2 \ell - e_i - e_j - e_k - e_l \\ \vphantom{\vdots}
C_i ^{jk} & 3 \ell -\bar E - e_i + e_j + e_k \\ \vphantom{\vdots}
D_{ijk} ^l & 4 \ell - \bar E - e_i - e_j - e_k + e_l \\ \vphantom{\vdots}
E_i ^j & 5 \ell - 2\bar E + e_i + 2e_j \\ \vphantom{\vdots}
D' _i & 4 \ell - \bar E - 2e_i \\ \vphantom{\vdots}
F_i ^{jkl} & 5 \ell - \bar E -2e_i - e_j - e_k - e_l \\ \vphantom{\vdots}
G_{ij} ^{kl} & 6 \ell - 2 \bar E - e_i - e_j + e_k + e_l
\end{array} \right. \hspace{-5pt}
\left\{ \begin{array} {l@{ \,=\, }l}
H_{ijk} ^l & 7 \ell - 3 \bar E + e_i + e_j + e_k + 2e_l \\ \vphantom{\vdots}
(H')_i ^l & 7 \ell - 2 \bar E - 2e_i - e_j \\ \vphantom{\vdots}
I_i ^{ijk} & 8 \ell - 3 \bar E - e_i + e_j + e_k + e_l \\ \vphantom{\vdots}
I' _i & 8 \ell - 3 \bar E + 2e_i \\ \vphantom{\vdots}
J_{ij} ^k & 9 \ell - 3 \bar E - e_i - e_j + e_k \\ \vphantom{\vdots}
K_{ijkl} & 10 \ell - 3 \bar E - e_i - e_j - e_k - e_l \\ \vphantom{\vdots}
L_i & 11 \ell - 4 \bar E + e_i 
\end{array} \right. 
\end{equation}

Denote by $\, \cdot \, $ the intersection form on the lattice 
${\rm Pic} (X_\delta )$.  From now on by an automorphism of 
${\rm Pic} (X_\delta )$ we will always mean a group automorphism of 
the lattice which preserves the intersection form and the canonical 
class; we let $W_\delta := {\rm Aut } \bigl( {\rm Pic } (X_\delta ) , 
K_{X_\delta } \, , \, \, \cdot \, \bigr)$, and we refer to $W_\delta$ 
as the Weyl group.  It will be useful later to know what are the 
orbits of pairs of conics under the automorphism group $W_\delta $ of 
${\rm Pic} (X_\delta )$.

\begin{lem} \label{trave}
The group $W_\delta $, $2 \leq \delta \leq 8$, acts transitively on the conics.
\end{lem}
{\it Proof.}  We only prove this in the case $\delta = 8$ and it will be 
clear from the proof that the same argument applies to the other cases.

Choose a standard basis $\{ \ell , e_1 , \ldots , e_8 \}$ of 
${\rm Pic} (X)$; it is enough to prove that the elements in the list 
(\ref{soluco}) are in the same orbit, since any permutation of the 
indices is an element of $W_8$.

Introduce the following automorphism of ${\rm Pic} (X_8)$: 
$$ T_{123} : \left\{ \begin{array} {l} 
\begin{array} {rcl@{\hspace{30pt}}rcl} 
\ell & \longmapsto & 2 \ell - e_1 - e_2 - e_3 \\ 
e _1 & \longmapsto & \ell - e _2 - e _3 \\
e _2 & \longmapsto & \ell - e _1 - e _3 \\
e _3 & \longmapsto & \ell - e _1 - e _2 \\
e _\alpha & \longmapsto & e_\alpha \end{array} \\
\text{ ~~\small $4 \leq \alpha \leq 8$} 
\end{array} \right. $$
and note that applying $T_{123}$ to an element 
$\bigl( a \,;\, b_1 , \ldots , b_8 \bigr)$ transforms it to 
$$\bigl( a \,;\, b_1 , \ldots , b_8 \bigr) \stackrel {T_{123}} 
{\longrightarrow} \bigl( 2a-b_1-b_2-b_3 \,;\, 
a-b_2-b_3 , a-b_1-b_3, a-b_1-b_2 , b_4 , \ldots , b_8 \bigr)$$
By inspection, the quantity $2a - b_1 - b_2 - b_3$ for elements 
in list (\ref{soluco}) is always strictly smaller than the 
initial value of $a$ unless $a=1$.  Permuting the indices so 
that $b_1,b_2,b_3$ are the three largest coefficients among the 
$b_i$'s and iterating this strategy finishes the argument.  
Note that we are always ``climbing up'' list (\ref{soluco}) and 
the conics on $X_7$ are the ones above line 5, and are hence 
preserved by the automorphism $T_{123}$ and the permutations 
needed.  Similar remarks are valid for $X_\delta $, with 
$3 \leq \delta \leq 6$, and the result is obvious for $X_2$, 
where the automorphism $T_{123}$ is not defined.  \bo

\noindent
{\it Remark}.  It is known (\cite{Ma} Theorem IV.23.9) that the group 
$W_\delta $ is generated by the permutations of the indices of the 
$e_i$'s together with the transformation $T_{123}$.

Suppose now we consider the action of the Weyl group on ordered pairs 
of conics $\bigl( Q_1 , Q_2 \bigr)$.  Clearly the number $Q_1 \cdot Q_2$ 
is an invariant of this action, and by looking at the list (\ref{soluco}) 
it is easy to convince oneself that 
$$ \begin{array} {@{ \vline~~ }c@{ ~~\vline~~ }c@{ ~~\vline~~ }c@{ ~~\vline~~ }c
@{ ~~\vline~~ }c@{ ~~\vline~~ }c@{ ~~\vline~~ }c@{ ~~\vline~~ }c@{ ~~\vline }} \hline
\delta =           & 8 & 7 & 6 & 5 & 4 & 3 & 2 \\ \hline
Q_1 \cdot Q_2 \leq & 8 & 4 & 2 & 2 & 1 & 1 & 1 \\ \hline
\end{array} $$
and that all the possible values between 0 and the number given above are 
attained.

Thus, for example, we know that the action of $W_8$ on pairs of conics 
has at least 9 orbits.

If $\delta = 8$, there is one more ``invariant'' under $W_8$ of pairs of 
conics: define a pair $\bigl( Q_1 , Q_2 \bigr)$ to be {\it ample} if 
$Q_1 + Q_2$ is an ample divisor on $X_8$.  Since the property of being 
ample is a numerical property, it follows that it is a property of the 
$W_8 -$orbit of the pair.

The next proposition proves that the lower bounds on the number of orbits 
obtained by considering the intersection product and ampleness (in case 
$\delta = 8$) of the pair are in fact the correct number of orbits.  Indeed, 
unless $\delta = 8$ it is enough to consider the intersection product, while 
if $\delta = 8$, there are two orbits with $Q_1 \cdot Q_2 = 4$, and only one 
of the two consists of ample pairs.

\begin{prop} \label{codico}
Let $Q_1$ and $Q_2$ be two conics in $X_\delta $, $2 \leq \delta \leq 8$.  
The intersection product $Q_1 \cdot Q_2$ determines uniquely the orbit 
of the (ordered) pair $\bigl( Q_1 , Q_2 \bigl)$ under $W_\delta $ with the 
only exception of $\delta = 8$ and $Q_1 \cdot Q_2 = 4$ which has exactly 
two orbits.
\end{prop}
{\it Proof.}  As for the previous lemma, we will only prove this proposition in 
the case $\delta = 8$; for the remaining cases simply ignore the inexistent 
indices.

Thanks to the previous lemma, we already know that we may assume 
$Q_1 = \ell - e_1$ which is the conic labeled $A_1$ in (\ref{nomico}).

The strategy is very simple: we again climb up the list (\ref{soluco}) 
using the automorphism $T_{123}$ followed by a permutation of the indices 
$\{ 2 , \ldots , 8 \}$ so that the resulting $b_2$ and $b_3$ are the two largest 
$b_i$'s, with $i \geq 2$.  Note that the elements of $W _\delta $ described above 
do indeed fix $Q_1$.

The case $Q_2 = 11 \ell - 4 \bar E + e_1$ is easily seen to be fixed 
by all permutations of $\{2, \ldots , 8 \}$ and by the automorphism 
$T_{123}$.  

In the following diagrams we write all possible conics with the 
given intersection product with $A_1$, sorting the entries 
$b_2, \ldots , b_8$ in non-increasing order.  An arrow going up 
means: apply $T_{123}$ and permute the indices different from 1 
so that the entries under $e_2 , \ldots , e_8$ are in 
non-increasing order.

\hspace{-45pt} \includegraphics{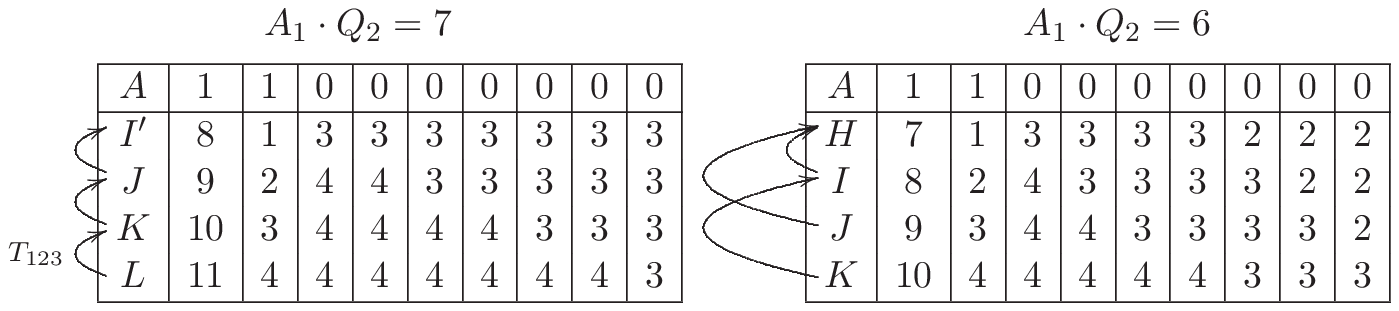}

\includegraphics{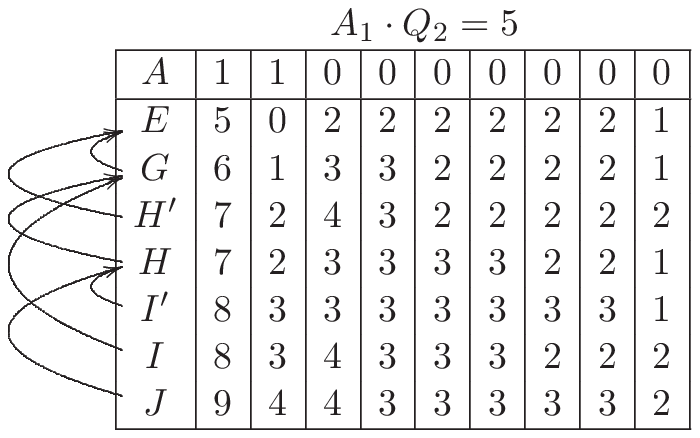}

Next is the case in which there is the exception.  Note that if $\delta = 7$, 
the possible intersection numbers $A_1 \cdot Q_2$ are at most 4, and 
$A_1 \cdot Q_2 = 4$ only if 
\(Q_2 = 5 \ell - e_1 - 2e_2 - 2e_3 - 2e_4 - 2e_5 - 2e_6 - 2e_7 \); thus the 
``top orbit'' of the next diagram does not appear for $\delta \leq 7$.

\hspace{-50pt} \includegraphics{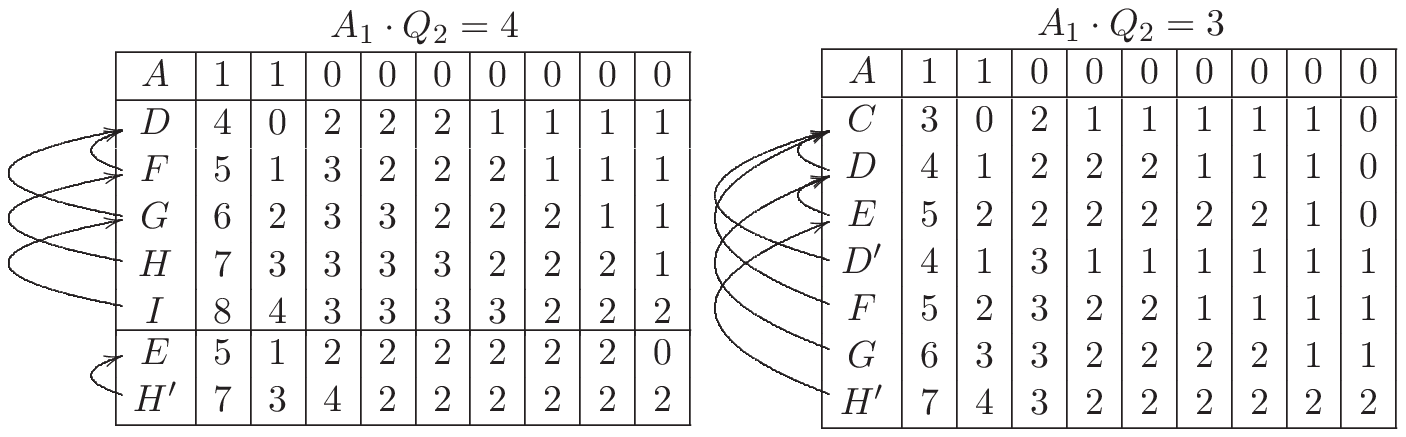}

\hspace{-50pt} \includegraphics{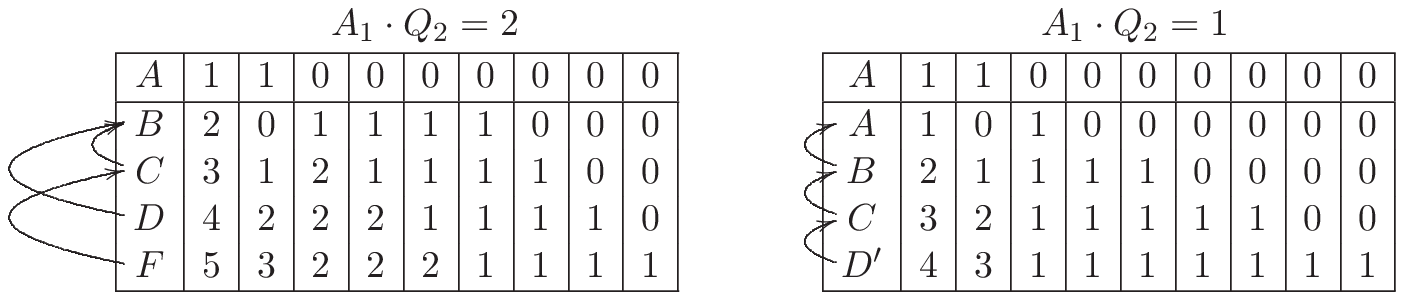}

Finally, note that $A_1 + D_{234} ^1 = -K_{X_8} + B_{234}$ is ample (being 
the sum of an ample divisor and a nef divisor), while 
$\bigl( A_1 + E_1 ^8 \bigr) \cdot e_8 = 0$.  Thus the pair 
$\bigl( A_1 , D_{234} ^1 \bigr)$ is ample, while the pair 
$\bigl( A_1 , E_1 ^8 \bigr)$ is not ample and therefore they cannot lie in 
the same orbit under the Weyl group.  This concludes the proof.  \bo

\noindent
{\it Remark}.  The same statement of Proposition \ref{codico} is clearly 
true if we are only interested in unordered pairs of conics.  This is 
obvious because the invariants we needed to detect all the orbits are 
invariants of unordered pairs, rather than ordered pairs.


The next two lemmas deal with a del Pezzo surface of degree one.

\begin{lem} \label{orco}
Let $X$ be a del Pezzo surface of degree one and let $L \subset X$ be a 
$(-1)-$curve.  If $L_1, L_2 \subset X$ are $(-1)-$curves such that 
$L_1 \cdot L = L_2 \cdot L$, then $L_1$ and $L_2$ are in the same orbit 
of the stabilizer of $L$ in ${\rm Aut} \bigl( {\rm Pic} (X) \bigr)$.
\end{lem}
\noindent
{\it Remark}.  The possible intersection numbers between any two 
$(-1)-$curves on a del Pezzo surface of degree one are -1, 0, 1, 2 and 3.  
Moreover, the group $W_8 := {\rm Aut} \bigl( {\rm Pic} (X) \bigr)$ acts 
transitively on $(-1)-$curves (\cite{Ma} Corollary 25.1.1).  
Thus as a consequence of this fact and the lemma we conclude that the 
stabilizer in the group $W_8$ of a $(-1)-$curve has exactly five orbits 
on the set of $(-1)-$curves.

\noindent
{\it Proof.}  We may choose a standard basis 
$\bigl\{ \ell , e_1 , \ldots , e_8 \bigr\}$ of ${\rm Pic} (X)$ such that 
$L = e_8$.  Given any divisor class $D \in {\rm Pic} (X)$, we write 
$D = a \ell - b_1 e_1 - \ldots - b_8 e_8$.  With these conventions, the 
classes of the $(-1)-$curves up to permutations of the indices 
$1, \ldots , 8$ are (\cite{Ma} Table IV.8)
\begin{equation} \label{radic8} 
\begin{array} {@{ \vline~~~ }l@{ ~~\vline~~ }c@{ ~~\vline~~ }c@{ ~~\vline~~ }
c@{ ~~\vline~~ }c@{ ~~\vline~~ }c@{ ~~\vline~~ }c@{ ~~\vline~~ }c@{ ~~\vline~~ }c@{ ~~\vline }}
\hline a & b_1 & b_2 & b_3 & b_4 & b_5 & b_6 & b_7 & b_8 \vphantom{\Bigl|} \\ 
\hline \vphantom{\Bigl|} 
0 &-1 & 0 & 0 & 0 & 0 & 0 & 0 & 0 \\[-2pt] 
1 & 1 & 1 & 0 & 0 & 0 & 0 & 0 & 0 \\[-2pt] 
2 & 1 & 1 & 1 & 1 & 1 & 0 & 0 & 0 \\[-2pt] 
3 & 2 & 1 & 1 & 1 & 1 & 1 & 1 & 0 \\[-2pt] 
4 & 2 & 2 & 2 & 1 & 1 & 1 & 1 & 1 \\[-2pt] 
5 & 2 & 2 & 2 & 2 & 2 & 2 & 1 & 1 \\[-2pt] 
6 & 3 & 2 & 2 & 2 & 2 & 2 & 2 & 2 \\ \hline 
\end{array} 
\end{equation}
and the stabilizer of $e_8$ contains the group $S$ generated by 
all permutations of $1, \ldots , 7$ and the automorphism $ T_{123}$ 
considered above.  In fact the stabilizer of $e_8$ is equal to the group $S$ 
just described, but we do not need this fact.

The proof consists simply in fixing one value for the coordinate 
$b_8$ and checking that all vectors with that last coordinate 
are in the same orbit of the group $S$.

\begin{description}
\item{$b_8 = 3$}.  There is only one vector in the list (\ref{radic8}) 
with an entry 3 in one of the $b_i$ columns and there is nothing to 
prove in this case.
\item{$b_8 = 2$}.  We have
\begin{eqnarray*}
T_{123} \bigl( 6 \,;\, 3 , 2 , 2 , 2 , 2 , 2 , 2 , 2 \bigr) & = & 
\bigl( 5 \,;\, 2 , 1 , 1 , 2 , 2 , 2 , 2 , 2 \bigr) \\
T_{145} \bigl( 5 \,;\, 2 , 1 , 1 , 2 , 2 , 2 , 2 , 2 \bigr) & = & 
\bigl( 4 \,;\, 1 , 1 , 1 , 1 , 1 , 2 , 2 , 2 \bigr) \\
T_{167} \bigl( 4 \,;\, 1 , 1 , 1 , 1 , 1 , 2 , 2 , 2 \bigr) & = & 
\bigl( 3 \,;\, 0 , 1 , 1 , 1 , 1 , 1 , 1 , 2 \bigr)
\end{eqnarray*}
and using permutations of $1, \ldots , 7$ we conclude.
\item{$b_8 = 1$}.  We have
\begin{eqnarray*}
T_{123} \bigl( 5 \,;\, 2 , 2 , 2 , 2 , 2 , 2 , 1 , 1 \bigr) & = & 
\bigl( 4 \,;\, 1 , 1 , 1 , 2 , 2 , 2 , 1 , 1 \bigr) \\
T_{145} \bigl( 4 \,;\, 1 , 1 , 1 , 2 , 2 , 2 , 1 , 1 \bigr) & = & 
\bigl( 3 \,;\, 0 , 1 , 1 , 1 , 1 , 2 , 1 , 1 \bigr) \\
T_{456} \bigl( 3 \,;\, 0 , 1 , 1 , 1 , 1 , 2 , 1 , 1 \bigr) & = & 
\bigl( 2 \,;\, 0 , 1 , 1 , 0 , 0 , 1 , 1 , 1 \bigr) \\
T_{236} \bigl( 2 \,;\, 0 , 1 , 1 , 0 , 0 , 1 , 1 , 1 \bigr) & = & 
\bigl( 1 \,;\, 0 , 0 , 0 , 0 , 0 , 0 , 1 , 1 \bigr) 
\end{eqnarray*}
\item{$b_8 = 0$}.  We have
\begin{eqnarray*}
T_{123} \bigl( 3 \,;\, 2 , 1 , 1 , 1 , 1 , 1 , 1 , 0 \bigr) & = & 
\bigl( 2 \,;\, 1 , 0 , 0 , 1 , 1 , 1 , 1 , 0 \bigr) \\
T_{145} \bigl( 2 \,;\, 1 , 0 , 0 , 1 , 1 , 1 , 1 , 0 \bigr) & = & 
\bigl( 1 \,;\, 0 , 0 , 0 , 0 , 0 , 1 , 1 , 0 \bigr) \\
T_{167} \bigl( 1 \,;\, 0 , 0 , 0 , 0 , 0 , 1 , 1 , 0 \bigr) & = & 
\bigl( 0 \,;\,-1 , 0 , 0 , 0 , 0 , 0 , 0 , 0 \bigr)
\end{eqnarray*}
\item{$b_8 =-1$}.  The only divisor class of a $(-1)-$curve with 
$b_8 =-1$ is $e_8$.
\end{description}

\noindent
This completes the cases we needed to check and the proof of the lemma.  \bo

The following is the last lemma of the section.

\begin{lem} \label{tremendi}
Let $L$ be the divisor class of a $(-1)-$curve on a del Pezzo surface $X$ 
of degree one, and let 
$$B := \left\{ \bigl\{ \lambda _1 , \lambda _2 , \lambda _3 \bigr\} ~ \Bigr| ~ 
\lambda _i {\text{ is a $(-1)-$curve, and }} 
\lambda _1 + \lambda _2 + \lambda _3 = -2K_X + L \right\} $$
The stabilizer in $W_8$ of $L$ has exactly four orbits on $B$.
\end{lem}
{\it Proof.}  Choose a standard basis of ${\rm Pic} (X)$ such that $L = e_8$.  
With this choice of basis, we have 
$$ \lambda _1 + \lambda _2 + \lambda _3 = \bigl( 6 \,;\, 2 , 2 , 2 , 2 , 2 , 2 , 2 , 1 \bigr) $$

Let $\beta _i$ be the coefficient of $-e_8$ in the chosen basis of $\lambda _i$.  
We deduce from above that 
\begin{eqnarray*}
& \beta _1 + \beta _2 + \beta _3 = 1 \\[5pt]
& -1 \leq \beta _i \leq 3 
\end{eqnarray*}
and thus, the solutions $\bigl\{ \beta _1 , \beta _2 , \beta _3 \bigr\}$ of 
the above system are $\bigl\{ 3 , -1 , -1 \bigr\}$, $\bigl\{ 2 , -1 , 0 \bigr\}$, 
$\bigl\{ 1 , 1 ,-1 \bigr\}$ and $\bigl\{ 1 , 0 , 0 \bigr\}$.

Permuting the $\lambda _i$'s we may assume that $\beta _1 \geq \beta _2 \geq \beta _3$ 
and using Lemma \ref{orco}, we may assume that the divisor class of $\lambda _1$ is 
$$ \begin{array} {c@{\text{ ~if~ }}l}
\bigl( 6 \,;\, 2 , 2 , 2 , 2 , 2 , 2 , 2 , 3 \bigr) & \beta _1 = 3, \\[5pt]
\bigl( 6 \,;\, 3 , 2 , 2 , 2 , 2 , 2 , 2 , 2 \bigr) & \beta _1 = 2, \\[5pt]
\bigl( 5 \,;\, 2 , 2 , 2 , 2 , 2 , 2 , 1 , 1 \bigr) & \beta _1 = 1.
\end{array} $$

It follows immediately that we must therefore have 
$$ \begin{array} {c}
\beta _1 = 3 ~:~ \left\{ \begin{array} {rcl}
\lambda _1 & = & \bigl( 6 \,;\, 2 , 2 , 2 , 2 , 2 , 2 , 2 , 3 \bigr) \\[5pt]
\lambda _2 & = & \bigl( 0 \,;\, 0 , 0 , 0 , 0 , 0 , 0 , 0 ,-1 \bigr) \\[5pt]
\lambda _3 & = & \bigl( 0 \,;\, 0 , 0 , 0 , 0 , 0 , 0 , 0 ,-1 \bigr) 
\end{array} \right. 
\end{array} $$ 
$$ \begin{array} {c} 
\beta _1 = 2 ~:~ \left\{ \begin{array} {rcl}
\lambda _1 & = & \bigl( 6 \,;\, 3 , 2 , 2 , 2 , 2 , 2 , 2 , 2 \bigr) \\[5pt]
\lambda _2 & = & \bigl( 0 \,;\, 0 , 0 , 0 , 0 , 0 , 0 , 0 ,-1 \bigr) \\[5pt]
\lambda _3 & = & \bigl( 0 \,;\,-1 , 0 , 0 , 0 , 0 , 0 , 0 , 0 \bigr) 
\end{array} \right. 
\end{array} $$

$$ \begin{array} {c} 
\begin{array} {c} 
\beta _1 = 1 \\[-2pt]
{\text{ and }} \\[-2pt]
\beta _2 = 1 
\end{array} ~:~ \left\{ \begin{array} {rcl}
\lambda _1 & = & \bigl( 5 \,;\, 2 , 2 , 2 , 2 , 2 , 2 , 1 , 1 \bigr) \\[5pt]
\lambda _2 & = & \bigl( 1 \,;\, 0 , 0 , 0 , 0 , 0 , 0 , 1 , 1 \bigr) \\[5pt]
\lambda _3 & = & \bigl( 0 \,;\, 0 , 0 , 0 , 0 , 0 , 0 , 0 ,-1 \bigr) 
\end{array} \right. 
\end{array} $$

$$ \begin{array} {c} 
\begin{array} {c} 
\beta _1 = 1 \\[-2pt]
{\text{ and }} \\[-2pt]
\beta _2 = 0 
\end{array} ~:~ \left\{ \begin{array} {rcl}
\lambda _1 & = & \bigl( 5 \,;\, 2 , 2 , 2 , 2 , 2 , 2 , 1 , 1 \bigr) \\[5pt]
\lambda _2 & = & \bigl( 0 \,;\,-1 , 0 , 0 , 0 , 0 , 0 , 0 , 0 \bigr) \\[5pt]
\lambda _3 & = & \bigl( 1 \,;\, 1 , 0 , 0 , 0 , 0 , 0 , 1 , 0 \bigr) 
\end{array} \right. 
\end{array} $$
thus proving the lemma.  \bo

\pagestyle{myheadings}
\markboth {3 \hspace{4pt} LARGE TO SMALL} {3 \hspace{4pt} LARGE TO SMALL}
\section{Realizing the Deformation: from Large to Small Degree}

\subsection{Breaking the Curve} \label{rompisezione}


In this section we construct deformations of a general point in 
every irreducible component of the space 
$\overline \cM _{bir} \bigl( X, \beta \bigr)$ to morphisms with 
image containing only curves of small anticanonical degree.

\begin{lem} \label{pezzenti}
Let $f : \bP^1 \rightarrow X$ be a free birational morphism to a del Pezzo 
surface.  In the same irreducible component of 
$\overline \cM _{bir} \bigl( X, f_*[\bP^1] \bigr)$ as $f$ there is a 
morphism $g : \bar C \rightarrow X$ birational to its image such that for 
every irreducible component $\bar C' \subset \bar C$, $g|_{\bar C'}$ is a 
free morphism whose image has anticanonical degree two or three.
\end{lem}
{\it Proof.}  We establish the lemma by induction on $d := -K_X \cdot f_*[\bP^1]$.  
There is nothing to prove if $d \leq 3$, since the image of a free morphism 
has anticanonical degree at least two (Lemma \ref{immersione}).

Suppose that $d \geq 4$.  Thanks to Proposition \ref{rpunti}, we may assume 
that the image of $f$ contains $d-2 \geq 2$ general points 
$p_1 , \ldots , p_{d-2}$ of $X$.  Denote by 
$\overline \cM _{bir} \bigl( p_1 , \ldots , p_{d-2} \bigr)$ the locus of 
morphisms of $\overline \cM _{bir} \bigl( X, f_*[\bP^1] \bigr)$ whose image 
contains the points $p_1 , \ldots , p_{d-2}$.  Using the dimension estimate 
(\ref{dimdibarbi}), we deduce that 
$\dim_{[f]} \overline \cM _{bir} \bigl( p_1 , \ldots , p_{d-2} \bigr) = 1$ 
and thus there is a one-parameter family of morphisms containing $f$ whose 
images contain the general points $p_1 , \ldots , p_{d-2}$.  Thanks to Lemmas 
\ref{nonnonred} and \ref{niette} we deduce that in the same irreducible 
component of $\overline \cM _{bir} \bigl( X, f_*[\bP^1] \bigr)$ as $f$ we can 
find a morphism $f_0 : \bar C_1 \cup \bar C_2 \rightarrow X$ such that $f_0$ 
is birational to its image, $\bar C_i \simeq \bP^1$ and $f_0 |_{\bar C_i}$ is 
a free morphism.  We also have $d_i := -K_X \cdot f_1(\bar C_i) \geq 2$, 
and thus by induction on $d$, we know that the irreducible component of 
$\overline \cM _{bir} \bigl( X, f_0 (\bar C_i) \bigr)$ containing 
$f_0|_{\bar C_i}$ contains a morphism 
$g_i : \bar C^i _1 \cup \ldots \cup \bar C^i _{r_i} \longrightarrow X$ with 
all the required properties.  Thus considering the morphism 
$$ \xymatrix @C=30pt { {\rm Sl}_{f_0} (\bar C_2) \ar[r] ^{\pi \hspace{20pt}} & 
\overline \cM _{bir} \bigl( X, f_0 (\bar C_2) \bigr) } $$
we deduce that we may find a morphism 
$f_1 : \bar C_1 \cup \bar C^2 _1 \cup \ldots \cup \bar C^2 _{r_2} \longrightarrow X$ 
with dual graph 
$$ \xygraph {[] !~:{@{.}} 
!{<0pt,0pt>;<20pt,0pt>:} 
{\bullet} [rr] {\bullet} [ur] {\bullet} [dd] {\bullet}
 [u] !{\save +<0pt,1pt>*{\txt{ {\begin{tabular} {c} 
\tiny dual \\[-7pt]
\tiny graph \\[-7pt]
\tiny of $g_2$ 
\end{tabular}}} }  \restore}
 [d] 
*\cir<2pt>{}
!{\save +<0pt,-8pt>*\txt{$\scriptstyle \bar C^2_c$}  \restore}
: [ul]
*\cir<2pt>{}
!{\save +<-3pt,8pt>*\txt{$\scriptstyle \bar C^2_a$}  \restore}
: [ur]
*\cir<2pt>{}
!{\save +<0pt,8pt>*\txt{$\scriptstyle \bar C^2_b$}  \restore}
 [dl] - [ll]
*\cir<2pt>{}
!{\save +<0pt,8pt>*\txt{$\scriptstyle \bar C_1$}  \restore} } $$
\vglue0pt {\centerline {Dual graph of $f_1$} \vspace{5pt}}

Similarly, considering the morphism 
$$ \xymatrix @C=30pt { {\rm Sl}_{f_1} (\bar C_1) \ar[r] ^{\pi \hspace{20pt}} & 
\overline \cM _{bir} \bigl( X, f_0 (\bar C_1) \bigr) } $$
we deduce that we may find a morphism 
$f_2 : \bar C^1 _1 \cup \ldots \cup \bar C^1 _{r_1} \cup 
\bar C^2 _1 \cup \ldots \cup \bar C^2 _{r_2} \longrightarrow X$ 
with dual graph 
$$ \xygraph {[] !~:{@{.}} 
!{<0pt,0pt>;<20pt,0pt>:} 
{\bullet} [d] !{\save +<0pt,1pt>*{\txt{ {\begin{tabular} {c} 
\tiny dual \\[-7pt]
\tiny graph \\[-7pt]
\tiny of $g_1$ 
\end{tabular}}} }  \restore} [d] 
{\bullet} [ur] {\bullet} [rr] {\bullet} [ur] {\bullet} [dd] {\bullet}
 [u] !{\save +<0pt,1pt>*{\txt{ {\begin{tabular} {c} 
\tiny dual \\[-7pt]
\tiny graph \\[-7pt]
\tiny of $g_2$ 
\end{tabular}}} }  \restore} [d] 
*\cir<2pt>{}
!{\save +<0pt,-8pt>*\txt{$\scriptstyle \bar C^2_c$}  \restore}
: [ul]
*\cir<2pt>{}
!{\save +<-3pt,8pt>*\txt{$\scriptstyle \bar C^2_a$}  \restore}
: [ur]
*\cir<2pt>{}
!{\save +<0pt,8pt>*\txt{$\scriptstyle \bar C^2_b$}  \restore}
 [dl] - [ll]
*\cir<2pt>{}
!{\save +<3pt,8pt>*\txt{$\scriptstyle \bar C^1_a$}  \restore} 
: [dl] 
*\cir<2pt>{}
!{\save +<0pt,-8pt>*\txt{$\scriptstyle \bar C^1_c$}  \restore} 
 [ur] : [ul] 
*\cir<2pt>{}
!{\save +<0pt,8pt>*\txt{$\scriptstyle \bar C^1_b$}  \restore} } $$
\vglue0pt {\centerline {Dual graph of $f_2$} \vspace{5pt}}

To conclude, we need to show that the images $C^1_a$ and $C^2_a$ of 
$\bar C^1_a$ and $\bar C^2_a$ respectively can be assumed to be distinct.

Suppose $C^1_a = C^2_a$.  If the anticanonical degree of $C^1_a$ is at 
least three, then we may deform one of them, keeping the image of the node 
between $\bar C^1_a$ and $C^2_a$ fixed and conclude.  Suppose therefore 
that $-K_X \cdot C^1_a = 2$.  Let $\varphi : \bar C^1_a \rightarrow \bar C^2_a$ 
be the morphism $(f_2)^{-1} \circ (f_2)|_{\bar C^1_a}$ and let 
$\bar p_i \in \bar C^i_a$ be the point in the intersection 
$\bar C^1_a \cap \bar C^2_a$.  There are two possibilities: either
$\varphi (\bar p_1) \neq \bar p_2$, or $\varphi (\bar p_1) = \bar p_2$.  
In the first case, the deformations of the morphism 
$f_2|_{\bar C^1_a \cup \bar C^2_a}$ fixing the component $\bar C^1_a$
actually change the image of the other component, allowing us to
conclude.  In the second case, there is a one-dimensional space of 
deformations of the stable map obtained by ``sliding the point $\bar
p_i$ along $\bar C^i_a$.''  Moreover, there must be components in the
image of $f_2$ different from $C^i_a$, since otherwise the morphism
$f$ could not have been birational to its image.  Thus we may assume
that $\bar C^2_a$ is adjacent to a curve mapped to a curve different
from $C^1_a = C^2_a$, call this curve $\bar D$ (remember that $g_i$ is 
birational to its image).  Let $\bar q \in \bar C^2_a$ be the node 
between $\bar C^2_a$ and $\bar D$.  We may slide the node $\bar p_i$ 
until it reaches the point $\bar q$ to obtain a morphism $f_3$ with 
dual graph 
$$ \xygraph {[] !~:{@{.}} 
!{<0pt,0pt>;<20pt,0pt>:} 
{\bullet} [d] !{\save +<0pt,1pt>*{\txt{ {\begin{tabular} {c} 
\tiny dual \\[-7pt]
\tiny graph \\[-7pt]
\tiny of $g_1$ 
\end{tabular}}} }  \restore} [d] 
{\bullet} [ur] {\bullet} [rr] {\bullet} [u] {\bullet} [drr] 
{\bullet} [ur] {\bullet} [dd] {\bullet}
 [u] !{\save +<0pt,1pt>*{\txt{ {\begin{tabular} {c} 
\tiny dual \\[-7pt]
\tiny graph \\[-7pt]
\tiny of $g_2$ 
\end{tabular}}} }  \restore} [d] 
*\cir<2pt>{}
!{\save +<0pt,-8pt>*\txt{$\scriptstyle \bar C^2_c$}  \restore}
: [ul]
*\cir<2pt>{}
!{\save +<-3pt,8pt>*\txt{$\scriptstyle \bar C^2_a$}  \restore}
: [ur]
*\cir<2pt>{}
!{\save +<0pt,8pt>*\txt{$\scriptstyle \bar C^2_b$}  \restore}
 [dl] - [ll]
*\cir<2pt>{}
!{\save +<0pt,-8pt>*\txt{$\scriptstyle \bar E$}  \restore} 
- [u] 
*\cir<2pt>{}
!{\save +<0pt,8pt>*\txt{$\scriptstyle \bar D$}  \restore} 
: [ur] [dl] : [ul] [ddr] - [ll]
*\cir<2pt>{}
!{\save +<3pt,8pt>*\txt{$\scriptstyle \bar C^1_a$}  \restore} 
: [dl] 
*\cir<2pt>{}
!{\save +<0pt,-8pt>*\txt{$\scriptstyle \bar C^1_c$}  \restore} 
 [ur] : [ul] 
*\cir<2pt>{}
!{\save +<0pt,8pt>*\txt{$\scriptstyle \bar C^1_b$}  \restore} } $$
\vglue0pt {\centerline {Dual graph of $f_3$} \vspace{5pt}}
\noindent
where the component labeled $\bar E$ is contracted to the point 
$f_2(\bar q)$.  Since the sheaf $f_3 ^* \cT_X$ is globally generated 
on each component of the domain of $f_3$ it follows that $f_3$ is a
smooth point of $\overline \cM _{0,0} \bigl( X, f_*[\bP^1] \bigr)$.

Clearly the morphism $f_3$ is also a limit of morphisms $f_4$ with
dual graphs 
$$ \xygraph {[] !~:{@{.}} 
!{<0pt,0pt>;<20pt,0pt>:} 
{\bullet} [d] !{\save +<0pt,1pt>*{\txt{ {\begin{tabular} {c} 
\tiny dual \\[-7pt]
\tiny graph \\[-7pt]
\tiny of $g_1$ 
\end{tabular}}} }  \restore} [d] 
{\bullet} [ur] {\bullet} [rr] {\bullet} [rr] {\bullet} [ur] 
{\bullet} [dd] {\bullet}
 [u] !{\save +<0pt,1pt>*{\txt{ {\begin{tabular} {c} 
\tiny dual \\[-7pt]
\tiny graph \\[-7pt]
\tiny of $g_2$ 
\end{tabular}}} }  \restore} [d] 
*\cir<2pt>{}
!{\save +<0pt,-8pt>*\txt{$\scriptstyle \bar C^2_c$}  \restore}
: [ul]
*\cir<2pt>{}
!{\save +<-3pt,8pt>*\txt{$\scriptstyle \bar C^2_a$}  \restore}
: [ur]
*\cir<2pt>{}
!{\save +<0pt,8pt>*\txt{$\scriptstyle \bar C^2_b$}  \restore}
 [dl] - [ll]
*\cir<2pt>{}
!{\save +<0pt,8pt>*\txt{$\scriptstyle \bar D$}  \restore} 
: [ur] [dl] : [ul] [dr] - [ll]
*\cir<2pt>{}
!{\save +<3pt,8pt>*\txt{$\scriptstyle \bar {C^1_a}'$}  \restore} 
: [dl] 
*\cir<2pt>{}
!{\save +<0pt,-8pt>*\txt{$\scriptstyle \bar C^1_c$}  \restore} 
 [ur] : [ul] 
*\cir<2pt>{}
!{\save +<0pt,8pt>*\txt{$\scriptstyle \bar C^1_b$}  \restore} } $$
\vglue0pt {\centerline {Dual graphs of the morphisms $f_4$} \vspace{5pt}}
\noindent
where $\bar C^1_a{}'$ is mapped to a general divisor linearly
equivalent to $C^1_a$ and transverse to it.  This concludes the proof
of the lemma.  \bo

\begin{lem} \label{piuma}
Let $f : \bar C := \bar C_1 \cup \ldots \cup \bar C_r \longrightarrow X$ be a
stable map of genus zero and suppose that $f_i := f|_{\bar C_i}$ is a
free morphism.  If $f(\bar C_1) \cdot f(\bar C_2) > 0$, then in the
same irreducible component of $\overline \cM _{0,0} 
\bigl( X, f_*[\bar C] \bigr)$ containing 
$[f]$ there is a morphism $g : \bar D_1 \cup \ldots \cup \bar D_r
\longrightarrow X$ such that $\bar D_1$ and $\bar D_2$ are adjacent, 
$g|_{\bar D_i}$ is a free morphism and $f_*[\bar C_i] = g_*[\bar D_i]$ 
for all $i$'s.
\end{lem}
{\it Proof.}  Renumbering the components of the domain of $f$, we may assume
that the curve 
$\bar C_{12} := \bar C_3 \cup \bar C_4 \ldots \cup \bar C_s$ is
the connected component of $\bar C_3 \cup \bar C_4 \ldots \cup \bar
C_r$ which has a point in common with both $\bar C_1$ and $\bar C_2$.
Moreover, we may also assume that no component of $\bar C_{12}$ is
mapped to a curve in the same divisor class as $\bar C_1$ or $\bar
C_2$.

Since all the morphisms $f|_{\bar C_i}$ are free, we may deform
$f|_{\bar C_{12}}$ to a free morphism with irreducible domain 
$\bar C_{12}'$.  Consider the morphism 
$$ \xymatrix @C=35pt 
{{\rm Sl}_f (\bar C_{12}) \ar[r]^{\pi \hspace{20pt}} & 
\overline \cM _{0,0} \bigl( X, f_*[\bar C_{12}] \bigr) } $$
and note that it is dominant on the component of $\overline \cM
_{0,0} \bigl( X, \bar C_{12} \bigr)$ containing $f|_{\bar C_{12}}$.
Thus we can find a morphism 
$$ f_1 : \bar C_1 \cup \bar C_{12}' \cup \bar C_2 \cup \bar C_{s+1} \cup 
\ldots \cup \bar C_r \longrightarrow X $$
with dual graph 
$$ \xygraph {[] !~:{@{.}} 
!{<0pt,0pt>;<20pt,0pt>:} 
{\bullet} [dd] {\bullet} [ur] {\bullet} [rr] {\bullet} [rr] {\bullet}
[ur] {\bullet} [dd] {\bullet} *\cir<2pt>{} : [ul] : [ur] [dl] 
*\cir<2pt>{}
!{\save +<-3pt,8pt>*\txt{$\scriptstyle \bar C_2$}  \restore}
- [ll]
*\cir<2pt>{}
!{\save +<0pt,8pt>*\txt{$\scriptstyle \bar C_{12}'$}  \restore}
- [ll]
*\cir<2pt>{}
!{\save +<3pt,8pt>*\txt{$\scriptstyle \bar C_1$}  \restore}
: [dl] [ur] : [ul] } $$
\vglue0pt {\centerline {Dual graph of $f_1$} \vspace{5pt}}

We want to deform $f_1$ to a morphism $f_2$ with dual graph 
$$ \xygraph {[] !~:{@{.}} 
!{<0pt,0pt>;<20pt,0pt>:} 
{\bullet} [dd] {\bullet} [ur] {\bullet} [rr] {\bullet} [u] {\bullet}
[drr] {\bullet} [ur] {\bullet} [dd] {\bullet} *\cir<2pt>{} 
: [ul] : [ur] [dl] 
*\cir<2pt>{}
!{\save +<-3pt,8pt>*\txt{$\scriptstyle \bar C_2$}  \restore}
- [ll]
*\cir<2pt>{}
!{\save +<0pt,-8pt>*\txt{$\scriptstyle \bar E$}  \restore}
- [u]
*\cir<2pt>{}
!{\save +<0pt,8pt>*\txt{$\scriptstyle \bar C_{12}'$}  \restore}
 [d] - [ll]
*\cir<2pt>{}
!{\save +<3pt,8pt>*\txt{$\scriptstyle \bar C_1'$}  \restore}
: [dl] [ur] : [ul] } $$
\vglue0pt {\centerline {Dual graph of $f_2$} \vspace{5pt}}
\noindent
where $\bar E$ is a contracted component.  This is immediate
considering the morphism
$$ \xymatrix @C=35pt 
{{\rm Sl}_{f_1} (\bar C_1) \ar[r]^{\hspace{10pt} a} & 
\bar C_{12}' } $$
and noting that it is dominant.

It is clear that we may similarly deform $f_2$ to a morphism 
$f_3$ obtained by sliding $\bar C_1'$ along $\bar C_2$ away from 
the component $\bar C_{12}'$.  The dual graph of the morphism $f_3$ is 
$$ \xygraph {[] !~:{@{.}} 
!{<0pt,0pt>;<20pt,0pt>:} 
{\bullet} [dd] {\bullet} [ur] {\bullet} [d] {\bullet}
[urr] {\bullet} [ur] {\bullet} [dd] {\bullet} *\cir<2pt>{} 
: [ul] : [ur] [dl] 
*\cir<2pt>{}
!{\save +<-3pt,8pt>*\txt{$\scriptstyle \bar C_2$}  \restore}
- [ll]
*\cir<2pt>{}
!{\save +<3pt,8pt>*\txt{$\scriptstyle \bar C_1''$}  \restore}
- [d]
*\cir<2pt>{}
!{\save +<10pt,0pt>*\txt{$\scriptstyle \bar C_{12}'$}  \restore}
 [u] : [dl] [ur] : [ul] } $$
\vglue0pt {\centerline {Dual graph of $f_3$} \vspace{5pt}}

To conclude we consider the morphism 
$$ \xymatrix @C=35pt 
{{\rm Sl}_{f_3} (\bar C_{12}') \ar[r]^{\pi \hspace{20pt}} & 
\overline \cM _{0,0} \bigl( X, f_*[\bar C_{12}] \bigr) } $$
to deform $f_3|_{\bar C_{12}'} \simeq f_1|_{\bar C_{12}'}$ back to
$f|_{\bar C_{12}}$ and conclude the proof of the lemma.  \bo

\subsection{Easy Cases: $\bP^2$, $\bP^1 \times \bP^1$ and $Bl_p (\bP^2)$}


This section proves the irreducibility of the spaces 
$\overline { \cM } _{bir} \bigl( X , \alpha \bigr)$ where $X$ 
is a del Pezzo surface of degree eight or nine.  Of course, in 
the case of $\bP^2$ this result is obvious: for a given degree 
$d$ of the image,  the space 
${\rm Hom }_{d} \bigl( \bP^1 , \bP^2 \bigr)$ of maps with image 
of degree $d$ is birational to the set of triples of homogeneous 
polynomials of degree $d$ up to scaling.  Since the space 
${\rm Hom }_{d} \bigl( \bP^1 , \bP^2 \bigr)$ dominates 
$\overline { \cM } _{0,0} \bigl( \bP^2 , d[line] \bigr)$, we 
deduce the stated irreducibility.  Similar considerations apply 
to $\bP^1 \times \bP^1$.  The result is less obvious for 
$Bl _p (\bP^2)$.  We prove the result for $Bl _p (\bP^2)$, but 
similar techniques would also apply to the other two cases.

Note that the same result for the cases $\bP^2$ and 
$\bP^1 \times \bP^1$ follows also from \cite{KP}.

\begin{thm} \label{barbapapa}
The spaces of stable maps $\overline { \cM } _{bir} \bigl( \bP^2 , \alpha \bigr)$, 
$\overline { \cM } _{bir} \bigl( \bP^1 \times \bP^1 , \beta \bigr)$ and 
$\overline { \cM } _{bir} \bigl( Bl_p(\bP^2) , \gamma \bigr)$ 
are irreducible for all divisor classes $\alpha$, $\beta $ and $\gamma $.
\end{thm}
{\it Proof.}  As remarked above, we only treat the case of $Bl _p (\bP^2)$.  
To simplify the notation, let $\bP$ denote $Bl _p (\bP^2)$.  Let 
$f: \bP^1 \rightarrow Bl _p (\bP^2)$ be a general morphism in 
an irreducible component of $\overline { \cM } _{bir} \bigl( \bP , \delta \bigr)$.

We first examine the cases $-K_\bP \cdot f_*[\bP^1] \leq 3$ separately.

If $K_\bP \cdot f_*[\bP^1] = -1$, then $f_* [\bP^1]$ is the unique 
$(-1)-$curve $E$.  In this case we clearly have 
$\overline { \cM } _{0,0} \bigl( Bl_p(\bP^2) , E \bigr) = 
\bigl\{ [f] \bigr\}$.

From now on, we may assume that $f(\bP^1)$ is not a $(-1)-$curve, 
and thus $f_* [\bP^1]$ is a nef divisor, since the only integral 
curve on $\bP$ having negative square is the exceptional divisor.

Suppose that $-K_\bP \cdot f_* [\bP^1] = 2$.  It follows easily that 
$f_*(\bP^1)$ is the strict transform of a line in $\bP^2$.  In this case the 
result is evidently true, since the mapping space is identified with the linear 
$|f_*(\bP^1)|$.

Suppose that $-K_\bP \cdot f_* [\bP^1] = 3$.  
Let $E \subset Bl_p (\bP^2)$ be the exceptional divisor and let 
$L \subset Bl_p (\bP^2)$ be an irreducible divisor representing the 
class obtained by pulling-back the divisor class of a line in $\bP^2$.  
The divisor classes of $L$ and $E$ generate the Picard group of 
$Bl_p (\bP^2)$.  The canonical divisor class on $Bl_p (\bP^2)$ is 
$K = -3 [L] + [E]$.

We have $f_* [\bP^1] = a \cdot [L] - b \cdot [E]$, with $a \geq b \geq 0$, since 
$f_* [\bP^1]$ is a nef divisor.  Moreover we know that $3a - b = 3$, and 
thus we see that we necessarily have $a=1$ and $b=0$.  Thus we deduce 
that 
$\overline { \cM } _{bir} \bigl( \bP , f_* [\bP^1] \bigr) \simeq 
\overline { \cM } _{0,0} \bigl( \bP^2 , [L] \bigr) \simeq \bP^2$, in the 
case of a divisor of anticanonical degree three.

Suppose that $-K_\bP \cdot f_*[\bP^1] \geq 4$.  We may use Lemma 
\ref{pezzenti} to deform $f$ to a morphism $f' : \bar C \rightarrow \bP$ 
where $\bar C = \bar C_1 \cup \ldots \cup \bar C_\ell $ are the 
irreducible components, all immersed by $f'$ and each having 
anticanonical degree two or three.

We have a morphism $f' : \bar C \rightarrow Bl_p (\bP^2)$, 
birational to its image, such that each component of $\bar C$ is mapped 
to a curve of anticanonical degree two or three.  We already saw that 
this means that each component represents one of the two divisor classes 
$[L] - [E]$ or $[L]$.

Since $f(\bar C)$ is connected, if all the components of $\bar C$ were 
mapped to curves whose divisor class is $[L] - [E]$ (there are at least 
two such components because we are assuming the anticanonical degree of 
the image is at least four), then they would 
all have the same image, which is ruled out by the fact that $f'$ is 
birational to its image.  It follows that at least one component of 
$\bar C$, say $\bar C_1$, is mapped to the divisor class $[L]$.

Using Lemma \ref{piuma} we may slide all the components of 
$\bar C$ mapped to the divisor class $[L] - [E]$ to be adjacent to the 
component $\bar C_1$.  After having done this, let $\bar F_1$, \ldots , 
$\bar F_l$ denote the components mapped to the divisor class 
$[L] - [E]$, and let $\bar C_1'$, \ldots , $\bar C_k'$ denote the 
components mapped to the divisor class $[L]$, where $\bar C_1'$ is the 
only component adjacent to all the components $\bar F_j$ and no other 
component $\bar C_r'$ is adjacent to any $\bar F_j$.

Consider the subgraph of the dual graph spanned by the components 
$\bar C_r'$; this is clearly a tree.  Suppose that one of the 
components adjacent to $\bar C_1'$ is $\bar C_2'$.  Using Lemma 
\ref{piuma}, we may slide all the components adjacent to $\bar C_1'$ 
(and mapped to $[L]$) 
to be adjacent to $\bar C_2'$, making $\bar C_1'$ a leaf of the resulting 
tree.  Similarly, considering the subgraph spanned by the components 
mapped to the divisor class $[L]$ different from $\bar C_1'$, we may 
again assume that $\bar C_2'$ is a leaf, and so on.  Eventually we 
end up with a morphism $g : \bar D \rightarrow Bl_p (\bP^2)$, where the 
components of $\bar D$ mapped to $[L] - [E]$ are $\bar F_1$, \ldots , 
$\bar F_l$ and the components mapped to $[L]$ are $\bar H_1$, \ldots , 
$\bar H_k$ and the dual graph of $g$ is 
\begin{equation} \label{grablo}
\xygraph {[] !~:{@{=}} 
!{<0pt,0pt>;<20pt,0pt>:} 
{\bullet} [ll] {\bullet} [ll] {\bullet} [ll] {\bullet} [ul] {\bullet} 
[dd] {\bullet} 
*\cir<2pt>{}
!{\save +<-8pt,0pt>*\txt{$\scriptstyle \bar F_l$}  \restore}
 [u] {\mathop {\vdots} \limits _{\vphantom {a}}} [u]
*\cir<2pt>{}
!{\save +<-8pt,0pt>*\txt{$\scriptstyle \bar F_1$}  \restore}
- [dr]
*\cir<2pt>{}
!{\save +<3pt,8pt>*\txt{$\scriptstyle \bar H_1 $}  \restore}
- [dl]
*\cir<2pt>{}
 [ur]
- [rr]
*\cir<2pt>{}
!{\save +<0pt,8pt>*\txt{$\scriptstyle \bar H_2 $}  \restore}
- [r] *\txt{\dots}
- [r]
*\cir<2pt>{}
!{\save +<0pt,8pt>*\txt{$\scriptstyle \bar H_{k-1} $}  \restore}
- [rr]
!{\save +<0pt,8pt>*\txt{$\scriptstyle \bar H_k $}  \restore} }
\end{equation}

Note that there are isomorphisms 
\begin{eqnarray*}
& \overline \cM _{0,0} \bigl( Bl_p (\bP^2) , [L] \bigr) \simeq 
\overline \cM _{0,0} \bigl( \bP^2 , [line] \bigr) \simeq 
\bigl( \bP^2 \bigr) ^\vee \\ & {\rm and} \\
& \overline \cM _{0,0} \bigl( Bl_p (\bP^2) , [L] - [E] \bigr) \simeq \bP^1
\end{eqnarray*}

Thus, since $[L] \cdot [L] = 1$ and 
$[L] \cdot ([L]-[E]) = 1$, we deduce that the space of all morphisms 
with dual graph (\ref{grablo}) is birational to 
$\bigl( \bP^2 \bigr) ^k \times \bigl( \bP^1 \bigr) ^l$, and in 
particular it is irreducible.  Since all the components of the 
morphisms with dual graph (\ref{grablo}) are free smooth rational 
curves, it follows that this locus contains smooth points of 
$\overline \cM _{bir} \bigl( Bl_p (\bP^2) , f_*[\bP^1] \bigr)$, and 
therefore we deduce that the space 
$\overline \cM _{bir} \bigl( Bl_p (\bP^2) , f_*[\bP^1] \bigr)$ is 
irreducible.  \bo

\section{Realizing the Deformation: from Small to Large Degree}
\pagestyle{myheadings}
\markboth{4 \hspace{4pt} SMALL TO LARGE}{4 \hspace{4pt} SMALL TO LARGE}

\subsection{Growing from the Conics} \label{conichette}


In this section we prove some results that allow us to deform 
unions of conics to divisors which are the anticanonical divisor 
on a del Pezzo surface dominated by $X$.  These results will be 
the main building blocks in the proof of Theorem \ref{passo}.

\begin{prop} \label{struttura}
Let $X_\delta $ be a del Pezzo surface of degree $9-\delta $ such 
that the spaces $\overline \cM _{bir} \bigl( X_\delta , \beta \bigr)$ 
are irreducible or empty if $-K_{X_\delta} \cdot \beta = 2, 3$.  In 
the case $\delta = 8$, or equivalently if the degree of $X_\delta $ is 
one, suppose also that all the rational divisors in the anticanonical 
linear system are nodal.  Let $f: \bar Q \rightarrow X_\delta $ be a 
morphism from a connected, projective, nodal curve of arithmetic 
genus zero.  Suppose that $\bar Q_1$ and $\bar Q_2$ are the irreducible 
components of $\bar Q$ and that $f_*[\bar Q_1]$ and $f_*[\bar Q_2]$ 
are conics.  If $f(\bar Q_1) \cdot f(\bar Q_2) \geq 2$, then in the 
irreducible component of 
$\overline \cM _{0,0} \bigl( X_\delta , f_*[\bar Q] \bigr)$ 
containing $[f]$ there is a morphism $g: \bar C \rightarrow X_\delta $ 
such that 
\begin{itemize}
\item all the irreducible components of $\bar C$ are immersed and 
represent nef divisor classes; 
\item there is a component $\bar C_1 \subset \bar C$ and a standard 
basis $\bigl\{ \ell , e_1 , \ldots , e_\delta \bigr\}$ of 
${\rm Pic } (X_\delta )$ with 
$$g_*[\bar C_1] = 3 \ell - e_1 - \ldots - e_\alpha $$
for some $\alpha \leq \delta $; 
\item if $g_*[\bar C_1] = -K_{X_8}$, then we may choose which of 
the twelve rational divisors in $|-K_{X_8}|$ the image of 
$\bar C_1$ is;
\item the point $[g]$ is smooth.
\end{itemize}
\end{prop}
{\it Proof.}  Observe that $f$ represents a smooth point of 
$\overline \cM _{0,0} \bigl( X_\delta , f_*[\bar Q] \bigr)$, 
since $f^*\cT_X$ is globally generated on both components of 
$\bar Q$.  Note also that by considering 
${\rm Sl }_f (\bar Q_i)$ we may assume that $Q_i:= f(\bar Q_i)$ 
misses any preassigned subscheme of $X$ of codimension 2.  In 
particular, we may suppose that $Q_i$ does not contain the 
intersection points between any two $(-1)-$curves.

We first take care of the case $Q_1 \cdot Q_2 = 2$: we may 
assume by Proposition \ref{codico} that $Q_1 = \ell - e_1$ and 
$Q_2 = 2 \ell - e_2 - e_3 - e_4 - e_5$.  It is therefore 
enough to smooth $\bar Q_1 \cup \bar Q_2$ to prove the 
proposition.

This concludes the proof if $\delta \leq 6$ since on a del Pezzo 
surface of degree at least three there do not exist conics $Q_1$ 
and $Q_2$ such that  $Q_1 \cdot Q_2 \geq 3$.

Suppose that $Q_1 \cdot Q_2 \geq 3$.  Our first step is to write 
$Q_2$ as a sum of two $(-1)-$curves $M_1$ and $M_2$ so that 
in some standard basis $\{ \ell' , e_1' , \ldots , e_\delta ' \}$ 
we have 
$$Q_1 + M_1 = \bigl( 3 \ell' - e_1' - \ldots - e_\alpha' \bigr) + N$$
where $N$ is a nef divisor.  We assume $Q_1 = A_1 = \ell - e_1$ 
(Lemma \ref{trave}).  Here is the explicit decomposition 
$Q_2 = M_1 + M_2$ in all the needed cases (Proposition \ref{codico}):

$$ \begin{array} {c} \\ Q_1 \cdot Q_2 = 3 \\ \fbox{\( \begin{array} {c} 
Q_2 = \bigl( 5 \,;\, 2,2,2,2,2,2,1,0 \bigr) \\[3pt]
M_1 = \bigl( 3 \,;\, 1,2,1,1,1,1,1,0 \bigr) \\
M_2 = \bigl( 2 \,;\, 1,0,1,1,1,1,0,0 \bigr) 
\end{array} \)} \end{array}
\hspace{11pt}
\begin{array} {c} Q_1 \cdot Q_2 = 4 \\ 
(Q_1+Q_2) \cdot e_8 = 0 \\ \fbox{\( \begin{array} {c} 
Q_2 = \bigl( 5 \,;\, 1,2,2,2,2,2,2,0 \bigr) \\[3pt]
M_1 = \bigl( 3 \,;\, 1,2,1,1,1,1,1,0 \bigr) \\
M_2 = \bigl( 2 \,;\, 0,0,1,1,1,1,1,0 \bigr) 
\end{array} \)} \end{array} $$
$$ \begin{array} {c} Q_1 \cdot Q_2 = 4 \\ 
Q_1+Q_2 \text{ ample (on $X_8$)}\\ \fbox{\( \begin{array} {c} 
Q_2 = \bigl( 4 \,;\, 0,2,2,2,1,1,1,1 \bigr) \\[3pt]
M_1 = \bigl( 3 \,;\, 0,2,1,1,1,1,1,1 \bigr) \\
M_2 = \bigl( 1 \,;\, 0,0,1,1,0,0,0,0 \bigr) 
\end{array} \)} \end{array} 
\hspace{11pt}
\begin{array} {c} \\ Q_1 \cdot Q_2 = 5 \\ \fbox{\( \begin{array} {c} 
Q_2 = \bigl( 5 \,;\, 0,2,2,2,2,2,2,1 \bigr) \\[3pt]
M_1 = \bigl( 3 \,;\, 0,2,1,1,1,1,1,1 \bigr) \\
M_2 = \bigl( 2 \,;\, 0,0,1,1,1,1,1,0 \bigr) 
\end{array} \)} \end{array} $$
$$ \begin{array} {c} Q_1 \cdot Q_2 = 6 \\ \fbox{\( \begin{array} {c} 
Q_2 = \bigl( 7 \,;\, 1,3,3,3,3,2,2,2 \bigr) \\[3pt]
M_1 = \bigl( 3 \,;\, 0,2,1,1,1,1,1,1 \bigr) \\
M_2 = \bigl( 4 \,;\, 1,1,2,2,2,1,1,1 \bigr) 
\end{array} \)} \end{array} 
\hspace{11pt}
\begin{array} {c} Q_1 \cdot Q_2 = 7 \\ \fbox{\( \begin{array} {c} 
Q_2 = \bigl( 8 \,;\, 1,3,3,3,3,3,3,3 \bigr) \\[3pt]
M_1 = \bigl( 3 \,;\, 0,2,1,1,1,1,1,1 \bigr) \\
M_2 = \bigl( 5 \,;\, 1,1,2,2,2,2,2,2 \bigr) 
\end{array} \)} \end{array} $$
$$ \hspace{32pt} \begin{array} {c} Q_1 \cdot Q_2 = 8 \\ \fbox{\( \begin{array} {r} 
Q_2 = \bigl( 11 \,;\, 3,4,4,4,4,4,4,4 \bigr) \\[3pt]
M_1 = \bigl(  5 \,;\, 1,2,2,2,2,2,2,1 \bigr) \\
M_2 = \bigl(  6 \,;\, 2,2,2,2,2,2,2,3 \bigr) 
\end{array} \)} \end{array} $$
\vspace{3pt}

Let us check that the previous decomposition has the required property: 
$$ \begin{array} {lr@{\,=\,}l}
Q_1 \cdot Q_2 = 3 & 
Q_1 + M_1 & T_{127} \bigl( 3 \ell - e_1 - \ldots - e_6 \bigr) \\[7pt]
Q_1 \cdot Q_2 = 4 & Q_1 + M_1 & \left\{ \hspace{-5pt} \begin{array} {ll}
T_{127} \bigl( 3 \ell - e_1 - \ldots - e_6 \bigr) & \hspace{-7pt} {\text{ if }}(Q_1+Q_2) \cdot e_8 = 0 \\[6pt]
-K_{X_8} + \bigl( \ell - e_2 \bigr) & \hspace{-7pt} {\text{ if }}Q_1+Q_2 \text{ is ample} 
\end{array} \right. \\[15pt]
Q_1 \cdot Q_2 = 5 & 
Q_1 + M_1 & -K_{X_8} + \bigl( \ell - e_2 \bigr) \\[7pt]
Q_1 \cdot Q_2 = 6 & 
Q_1 + M_1 & -K_{X_8} + \bigl( \ell - e_2 \bigr) \\[7pt]
Q_1 \cdot Q_2 = 7 & 
Q_1 + M_1 & -K_{X_8} + \bigl( \ell - e_2 \bigr) \\[7pt]
Q_1 \cdot Q_2 = 8 & 
Q_1 + M_1 & -K_{X_8} + \bigl(  3 \ell - e_1 - \ldots - e_7 \bigr)
\end{array} $$

Next we show that we can deform $f$ so that the dual graph of 
the resulting morphism $f_1$ is 
\begin{equation} \label{quasi}
\xygraph {[] !~:{@{=}} 
!{<0pt,0pt>;<20pt,0pt>:} 
{\bullet} [rr] {\bullet} [rr] {\bullet} 
*\cir<2pt>{}
!{\save +<0pt,8pt>*\txt{$\scriptstyle \bar M_2$}  \restore}
- [ll]
*\cir<2pt>{}
!{\save +<0pt,8pt>*\txt{$\scriptstyle \bar M_1$}  \restore}
- [ll]
*\cir<2pt>{}
!{\save +<0pt,8pt>*\txt{$\scriptstyle \bar Q _1$}  \restore} }
\end{equation}
\vglue0pt {\centerline {Dual graph of $f_1$} \vspace{5pt}}
\noindent
where of course $\bar M_i$ maps to the $(-1)-$curve with divisor 
class $M_i$.  To achieve this, consider 
$$\xymatrix { a : {\rm Sl}_f (\bar Q_2) \ar[r] & \bar Q_1 } $$

The morphism $a$ is not constant because $f|_{\bar Q_2}$ is free, 
and hence it is surjective.  We denote with the symbols $M_1$ and 
$M_2$ both the divisor classes and the $(-1)-$curves on $X$ with 
the same divisor class.  Let $\bar p \in \bar Q_1$ be a point 
such that $f(\bar p) =: p \in M_1$; such a point exists, since 
$Q_1 \cdot M_1 \geq 2$ by inspection.

Thanks to the surjectivity of $a$, we may find 
$f_1 : \bar Q_1 \cup \bar Q_2' \rightarrow X$ such that 
$a(f_1) = \bar p$, and in particular, the node between $\bar Q_1$ 
and $\bar Q_2'$ maps to $p \in M_1$.  Since $Q_2 \cdot M_1 = 0$ and 
since $f_1(\bar Q_2') \cap M_1 \ni p$, it follows that 
$f_1 (\bar Q_2') \supset M_1$.  Thus we have that 
$\bar Q_2' = \bar M_1 \cup \bar M_2$, where $\bar M_2$ maps to the 
$(-1)-$curve $M_2 \subset X$; the dual graph of $f_1$ is the one in 
(\ref{quasi}): by construction $\bar Q_1$ and $\bar M_1$ are adjacent, 
and by connectedness of $\bar Q_2$ it follows that $\bar M_1$ and 
$\bar M_2$ are adjacent; the assumption that $Q_1$ does not 
contain the intersections of two $(-1)-$curves shows that there 
cannot be contracted components.  Note that the node between 
$\bar M_1$ and $\bar M_2$ maps to a node, since the intersection 
number $M_1 \cdot M_2$ equals one.

Let us check that $f_1$ represents a smooth point of its moduli 
space.  Thanks to Proposition \ref{grafico}, we have that the sheaf 
$\cC_1 := \cC _{f_1} \otimes \omega _{\bar Q_1 \cup \bar Q_2'}$, 
whose global sections represent the obstructions, has degrees 
given by the following diagram:
$$ \xygraph{[] !~:{@{.}} 
!{<0pt,0pt>;<20pt,0pt>:} 
{\bullet} [rr] {\bullet} [rr] {\bullet} 
*\cir<2pt>{}
!{\save +<0pt,8pt>*\txt{$\scriptstyle \bar M_2$}  \restore}
!{\save +<0pt,-8pt>*\txt{$\scriptstyle -1$}  \restore}
- [ll]
*\cir<2pt>{}
!{\save +<0pt,8pt>*\txt{$\scriptstyle \bar M_1$}  \restore}
!{\save +<0pt,-8pt>*\txt{$\scriptstyle \leq 0$}  \restore}
: [ll]
!{\save +<0pt,8pt>*\txt{$\scriptstyle \bar Q_1$}  \restore}
!{\save +<0pt,-8pt>*\txt{$\scriptstyle \leq -1$}  \restore} } $$
\vglue0pt {\centerline {Multi-degree of $\cC_1$} \vspace{5pt}}
\noindent
A solid edge means that the sheaf $\cC_1$ is locally free at 
the corresponding node, while a dotted edge means that the sheaf 
$\cC_1$ need not be locally free at that node (we could make 
sure that the sheaf $\cC_1$ is locally free by reducing to the 
case in which $Q_1$ intersects transversely $M_1$, but this is 
not needed).  It is now clear that $\cC_1$ has no global sections, 
and thus the point $f_1$ is smooth.

We smooth the components $\bar Q_1 \cup \bar M_1$ to a single 
irreducible component $\bar Q_1'$.  We obtain a morphism 
$g' : \bar Q_1' \cup \bar M_2 \rightarrow X$, such that in some 
standard basis $\{ \ell' , e_1 , \ldots , e_\delta' \}$ we have 
$g'_* [\bar Q_1'] = \bigl( 3 \ell' - e_1' - \ldots - e_\alpha' \bigr) + N$ 
where $\alpha \geq 6$ and $N$ is a nef divisor.  By construction 
the anticanonical degree of $g'_* [\bar Q_1']$ is three.

In the first two cases above, that is if $Q_2$ equals 
$\bigl( 5\,;\, 2,2,2,2,2,2,1,0 \bigr)$ or 
$\bigl( 5\,;\, 1,2,2,2,2,2,2,0 \bigr)$, the divisor $N$ above 
is zero, but in both cases we may write 
$$ Q_1' = \bigl( 3\,;\, 1,1,1,1,1,1,1,0 \bigr) + \bigl( 1\,;\, 1,1,0,0,0,0,0,0 \bigr) $$
We let $C_2$ be the $(-1)-$curve with divisor class 
$\bigl( 1\,;\, 1,1,0,0,0,0,0,0 \bigr)$.  By inspection we see that 
$C_2 \cdot M_2 \geq 1$, and therefore we may find a point $\bar c$ of 
$\bar Q_2$ such that $g' (\bar p) \in C_2$.  Considering the morphism 
$$\xymatrix { a : {\rm Sl}_{g'} (\bar Q_1') \ar[r] & \bar M_2 } $$
we let $g_1 : \bar C_1 \cup \bar C_2 \cup \bar M_2 \longrightarrow X$ 
be a morphism such that $a(g_1) = \bar c$, where we denote by $\bar C_2$ 
the component mapped to $C_2$ and by $\bar C_1$ the component mapped to 
$g'_*[\bar Q_1'] - C_2 = \bigl( 3\,;\, 1,1,1,1,1,1,1,0 \bigr)$.  By 
construction, the dual graph of $g_1$ is 
$$ \xygraph {[] !~:{@{.}} !{<0pt,0pt>;<20pt,0pt>:} 
{\bullet} [rr] {\bullet} [rr] {\bullet} 
*\cir<2pt>{}
!{\save +<0pt,8pt>*\txt{$\scriptstyle \bar M_2$}  \restore}
- [ll]
*\cir<2pt>{}
!{\save +<0pt,8pt>*\txt{$\scriptstyle \bar C_2$}  \restore}
- [ll]
!{\save +<0pt,8pt>*\txt{$\scriptstyle \bar C_1$}  \restore} } $$
\vglue0pt {\centerline {Dual graph of $g_1$} \vspace{5pt}}
\noindent
Smoothing the components $\bar C_2 \cup \bar M_2$ we conclude the proof 
of the proposition in these cases.

In the remaining cases (the ones for which $Q_1+Q_2$ is ample on a del 
Pezzo surface of degree one) we write $g'_*[\bar Q_1'] = -K_{X_8} + N$, 
where $N$ is $\bigl( \ell - e_2 \bigr)$, if $Q_1 \cdot Q_2 \leq 7$ and 
$N$ is $\bigl(  3 \ell - e_1 - \ldots - e_7 \bigr)$, if 
$Q_1 \cdot Q_2 = 8$.

By assumption the space 
$\overline \cM _{bir} \bigl( X_\delta , g'_*[\bar Q_1'] \bigr)$ 
is irreducible.  We may therefore deform the morphism $g'$ to a 
morphism $g_1 : \bar K \cup \bar N \cup \bar M_2 \longrightarrow X$, 
such that $\bar K$ is mapped to any preassigned rational divisor in 
$|-K_{X_8}|$ and $\bar N$ is mapped to a general divisor in $|N|$.  
The possible dual graphs for $g_1$ are 
$$ \xygraph {[] !~:{@{.}} !{<0pt,0pt>;<20pt,0pt>:} 
{\bullet} [rr] {\bullet} [rr] {\bullet} 
*\cir<2pt>{}
!{\save +<0pt,8pt>*\txt{$\scriptstyle \bar M_2$}  \restore}
- [ll]
*\cir<2pt>{}
!{\save +<0pt,8pt>*\txt{$\scriptstyle \bar N$}  \restore}
- [ll]
!{\save +<0pt,8pt>*\txt{$\scriptstyle \bar K$}  \restore} } 
\hspace{30pt}
\xygraph {[] !~:{@{.}} !{<0pt,0pt>;<20pt,0pt>:} 
{\bullet} [rr] {\bullet} [rr] {\bullet} 
*\cir<2pt>{}
!{\save +<0pt,8pt>*\txt{$\scriptstyle \bar M_2$}  \restore}
- [ll]
*\cir<2pt>{}
!{\save +<0pt,8pt>*\txt{$\scriptstyle \bar K$}  \restore}
- [ll]
!{\save +<0pt,8pt>*\txt{$\scriptstyle \bar N$}  \restore} } $$
\vglue0pt {\centerline {Possible dual graphs of $g_1$} \vspace{5pt}}

We smooth the two components $\bar M_2$ and the one adjacent 
to it.  In either case, the proposition is proved: this is 
obvious if $\bar N$ is adjacent to $\bar M_2$; if $\bar K$ is 
adjacent to $\bar M_2$, note that $-K_{X_8} + M_2$ is the 
pull-back of the anticanonical divisor on the del Pezzo surface 
obtained by contracting $M_2$.  This concludes the proof of the 
proposition.  \bo

{\it Remark.}  The proof above only requires the existence of 
one nodal rational divisor in $|-K_{X_8}|$.

\begin{prop} \label{infradito}
Let $X_\delta $ be a del Pezzo surface of degree $9-\delta $.  
Let $f: \bar Q \rightarrow X_\delta $ be a morphism from a connected, 
projective, nodal curve of arithmetic genus zero.  Suppose that 
$\bar Q_1$, $\bar Q_2$, $\bar Q_3$ are the irreducible components of 
$\bar Q$ and that $f_*[\bar Q_i]$ is a conic, for all $i$.  If 
$f(\bar Q_i) \cdot f(\bar Q_j) = 1$ for all $i \neq j$, then in the 
irreducible component of 
$\overline \cM _{0,0} \bigl( X_\delta , f_*[\bar Q] \bigr)$ 
containing $[f]$ there is an immersion 
$g: \bar C \rightarrow X_\delta $ such that $\bar C$ is irreducible 
and $g_* [\bar C] = 3 \ell - e_1 - e_2 - e_3$, for some choice of 
standard basis $\{ \ell , e_1 , \ldots , e_\delta \}$.
\end{prop}
{\it Proof.}  It is enough to show that we may find a standard basis such 
that $f_*[\bar Q_i] = \ell - e_i$, for $i \in \{ 1,2,3 \}$, since then 
smoothing out all the components we conclude.  Denote by $Q_i$ the 
image of $\bar Q_i$.  Thanks to Proposition \ref{codico} we may 
assume that $Q_1 = \ell - e_1$ and $Q_2 = \ell - e_2$.  Looking at 
the list (\ref{soluco}) we easily see that either we may assume 
that $Q_3 = \ell - e_3$ and we are done, or 
$Q_3 = 2 \ell - e_1 - e_2 - e_3 - e_4$, up to permutations of the 
coordinates.  In this last case, we apply $T_{124}$ to all three 
divisor classes.  Both $Q_1$ and $Q_2$ are fixed by $T_{124}$, 
while $T_{124} \bigl( Q_3 \bigr) = \ell - e_3$.  \bo

\subsection{Reduction of the Problem to Finitely Many Cases} \label{pochipochi}


This section gathers the information obtained in the previous sections 
to prove that the irreducibility of 
$\overline \cM _{bir} \bigl( X , \beta \bigr)$ for all $\beta $ can 
be checked by examining only finitely many cases.

First we prove two simple results.

\begin{lem}
Let $X$ be a smooth projective surface and let $D \in {\rm Pic} (X)$ 
be a base-point free nef divisor such that $D^2 > 0$.  If $N$ is a 
nef divisor such that $D \cdot N = 0$, then $N \equiv 0$.
\end{lem}
{\it Proof.}  Let $\varphi : X \rightarrow \bP ^n$ be the morphism induced by 
the linear system $|D|$ and denote by $X'$ the image of $\varphi $.  We 
clearly have that $D' := \varphi _* [D]$ is an ample divisor on $X'$.

The push-forward of a nef divisor $N$ on $X$ is a nef on $X'$: let 
$C \subset X'$ be an effective curve; we have 
$\varphi _*N \cdot C = N \cdot \varphi ^* C \geq 0$, since $\varphi ^* C$ 
is an effective curve.

Let $N$ be a nef divisor on $X$ such that $N \cdot D = 0$.  We have 
$\varphi _*N \cdot D' = N \cdot \varphi ^* D' = N \cdot D = 0$, and 
therefore by the Hodge Index Theorem we deduce that either 
$\varphi _*N \equiv 0$ or $(\varphi _*N) ^2 < 0$.  Since $\varphi _*N$ 
is nef, it is a limit of ample divisors and it follows that 
$(\varphi _*N) ^2 \geq 0$.  We deduce that $\varphi _*N \equiv 0$ and 
thus that $N$ is numerically equivalent to a linear combination of curves 
contracted by $\varphi $.  Since the intersection form on the span of the 
contracted curves is negative definite and $N$ is nef, we deduce that 
$N \equiv 0$.  \bo

\begin{cor} \label{benaco}
Let $X$ be a del Pezzo surface and let $D$ be a nef divisor on $X$ 
which is not a multiple of a conic.  If a nef divisor $N$ on $X$ 
is such that $D \cdot N = 0$, then $N=0$.
\end{cor}
{\it Proof.}  The result is obvious in the case $X = \bP^2$.  
Thanks to the previous lemma and the fact that numerical 
equivalence is the same as equality of divisor classes on a del Pezzo 
surface, it is enough to check that a multiple of a nef divisor class 
$D$ on $X$ is base-point free and has positive square, unless $D$ is the 
divisor class of a conic.

Write $D = n_\delta (-K_{X_\delta}) + \ldots + n_2 (-K_2) + D'$ as in 
Corollary \ref{maquale}.  It is immediate to check that $2D$ is 
base-point free (in fact, unless $D = -K_X$ and $X$ has degree one, then 
$D$ itself is base-point free).  If one of the $n_\alpha $'s is non-zero, 
then clearly the square of $D$ is positive (note that all the divisors 
appearing in the above expression of $D$ are nef and thus effective since 
$X$ is a del Pezzo surface).  If all the $n_\alpha $'s are zero, then 
$D = D'$ is a nef divisor on a del Pezzo surface of degree eight.

If $X = \bP^1 \times \bP^1$, let $\ell _1$ and $\ell _2$ be the two 
divisor classes $\{ p \} \times \bP^1$ and $\bP^1 \times \{ p \}$ 
respectively.  Any nef divisor class is a non-negative linear combination 
of $\ell _1$ and $\ell _2$; thus we may write $D = a_1 \ell _1 + a_2 \ell _2$, 
with $a_1, a_2 \geq 0$.  Moreover, if one of the $a_i$'s were zero, then $D$ 
would be a multiple of a conic: we deduce that $a_i > 0$.  Thus we compute 
$D^2 = 2 a_1 a_2 > 0$.

If $X = Bl_p (\bP^2 )$, let $\ell $ and $e$ be the pull-back of the divisor 
class of a line and the exceptional divisor under the blow down morphism to 
$\bP^2$ respectively.  Any nef divisor class is a non-negative linear 
combination of $\ell $ and $\ell - e$; thus we may write 
$D = a \ell + b (\ell - e)$, with $a, b \geq 0$.  Moreover, if $a = 0$, then 
$D$ is a multiple of a conic: we deduce that $a > 0$.  Thus we compute 
$D^2 = a (a + 2 b) > 0$ and the proof is complete.  \bo

We are now ready to prove the main result of the section.  The 
proof involves several steps and is quite long.

\begin{thm} \label{passo}
Let $X$ be a del Pezzo surface such that the spaces 
$\overline \cM _{bir} \bigl( X, \beta \bigr)$ are irreducible (or empty) 
for all nef divisors $\beta $ such that $2 \leq - K_X \cdot \beta \leq 3$.  
In the case $\deg X = 1$, suppose that all the rational divisors in the 
anticanonical system are nodal.  Then, for any nef divisor $D \subset X$ 
such that $-K_X \cdot D \geq 2$, the space 
$\overline \cM _{bir} \bigl( X, D \bigr)$ is irreducible or empty.
\end{thm}
{\it Proof.}  We establish the theorem by induction on $d := -K_X \cdot D$.  By 
hypothesis, the theorem is true if $d \leq 3$.

Suppose that $d \geq 4$.  Let $f : \bP^1 \rightarrow X$ be a general morphism 
in an irreducible component of $\overline \cM _{bir} \bigl( X, D \bigr)$.  
Since the morphism $f$ is a general point on an irreducible component of 
$\overline \cM _{bir} \bigl( X, D \bigr)$ and $d \geq 2$, it follows that $f$ 
is an immersion and that it is a free morphism.

If there is a $(-1)-$curve $L \subset X$ such that $L \cdot D = 0$, then 
let $b : X \rightarrow X'$ be the contraction of $L$.  We have 
$\overline \cM _{bir} \bigl( X, D \bigr) \simeq 
\overline \cM _{bir} \bigl( X', b_* D \bigr)$, and thus we reduce to the case 
in which the divisor $D$ intersects strictly positively every $(-1)-$curve.  
By Theorem \ref{barbapapa} we may also assume that the degree of $X$ is at 
most seven.  Thus Proposition \ref{clane} implies that $D$ is an ample divisor.

Thanks to Lemma \ref{pezzenti} we may deform $f$ to a morphism 
$g : \bar C \rightarrow X$ such that each component $\bar C_0 \subset
\bar C$ is immersed to a curve of anticanonical degree two or three.  We want 
to show that we may specialize $g$ to a morphism in which one component is 
mapped to a multiple of the divisor class $-K_X$.  We will prove this in a 
series of steps.

{\bf Step 1.}  There is a standard basis $\{ \ell , e_1 , \ldots , e_\delta \}$ 
of ${\rm Pic} (X)$ and a component $\bar C_1$ of $\bar C$ mapped birationally 
either to the divisor class $3 \ell - e_1 - \ldots - e_{\alpha}$, for 
$\alpha \in \{ 1 , \ldots , 7 \}$, or to $-rK_{X_8}$, for $r \in \{ 1,2,3 \}$.  
If the image of $\bar C_1$ represents $-K_{X_8}$, then we can choose to which 
of the twelve rational divisors in $|-K_{X_8}|$ the component $\bar C_1$ maps.  
The morphism is free on all the components of $\bar C$, except on $\bar C_1$ 
if it represents $-K_{X_8}$.

The divisors of anticanonical degree two on $X$ are 
\begin{itemize}
\item the divisor $-2K_X$, if $\deg X = 1$;
\item the divisor $-K_X$, if $\deg X = 2$;
\item the divisor class of a conic.
\end{itemize}
The divisors of anticanonical degree three on $X$ are 
\begin{itemize}
\item the divisor $-3K_X$, if $\deg X = 1$;
\item the divisor $-K_X-K_{X'}$, if $\deg X = 1$ and $X'$ is obtained from $X$ 
by contracting a $(-1)-$curve;
\item the divisor $-K_X + C$, if $\deg X = 1$ and $C$ is the class of a conic;
\item the divisor $-K_X$, if $\deg X = 3$;
\item the divisor $\ell $, for some standard basis $\{ \ell , e_1 , \ldots , e_\delta \}$.
\end{itemize}

Thanks to the irreducibility assumption on the spaces 
$\overline \cM _{bir} \bigl( X, \beta \bigr)$, for $2 \leq -K_X \cdot \beta \leq 3$, 
we reduce to the case in which all components of $\bar C$ are mapped to 
either the divisor class of a conic or the divisor class $\ell $, for some 
choice of standard basis.

We reduce further to the following case: 

\vspace{17pt}
\noindent
\begin{tabular} {l@{\hspace{7pt}}c}
($\star$) & 
\begin{tabular} {c} There is a standard basis $\{ \ell , e_1 , \ldots , e_\delta \}$ of 
${\rm Pic} (X)$ such that all \\
curves of degree three in the image of $g$ 
have divisor class $\ell $.
\end{tabular} \end{tabular}
\vspace{17pt}

This is easily accomplished.  Suppose that $\bar C_1$ and $\bar C_2$ are components 
of $\bar C$ such that $g_*[\bar C_1] = \ell _1$ and $g_*[\bar C_2] = \ell _2$, 
where $\{ \ell _i , e^i_1 , \ldots , e^i_\delta \}$ are two standard basis of 
${\rm Pic} (X)$ and $\ell _1 \neq \ell _2$.  We may first of all apply Lemma 
\ref{piuma} to reduce to the case in which $\bar C_1$ and $\bar C_2$ are adjacent 
in the dual graph of $g$.  If $\ell _2$ were orthogonal to 
$e^1_1 , \ldots , e^1_\delta $, then $\ell _2$ would be proportional and hence equal 
to $\ell_1$.  It follows that $\ell _2$ is not orthogonal to all the $e^1_j$'s.  By 
permuting the indices if necessary, we may assume that $\ell _2 \cdot e^1_1 > 0$.  
Since $g|_{\bar C_2}$ is free, we may assume that $g(\bar C_2)$ and $E := E^1_1$, 
the $(-1)-$curve whose divisor class is $e^1_1$, meet transversely.  Denote by 
$\bar p \in \bar C_2$ a point such that $g(\bar p) \in E$.  Consider the morphism 
$$ \xymatrix @C=35pt 
{{\rm Sl}_{g} (\bar C_1) \ar[r]^{\hspace{10pt} a} & \bar C_2 } $$
and note that it is dominant, since $g|_{\bar C_1}$ is free.  It follows that 
we may find a morphism 
$g_1 : \bar C_1' \cup \bar C_2 \cup \ldots \cup \bar C_r \longrightarrow X$ such 
that $a(g_1) = \bar p$.  We deduce that $g_1(\bar C_1') \ni p$ and 
$(g_1)_* [\bar C_1'] = \ell _1$.  Since $\ell _1 \cdot e^1_1 = 0$, we conclude that 
$g_1(\bar C_1')$ contains $E$ and another (irreducible) component whose divisor class 
is $\ell _1 - e^1_1$.  Finally, the subgraph of the dual graph of $g_1$ spanned by 
$\bar C_1'$ and $\bar C_2$ is 
$$ \xygraph {[] !~:{@{.}} !{<0pt,0pt>;<20pt,0pt>:} 
{\bullet} [rr] {\bullet} [rr] {\bullet} 
*\cir<2pt>{}
!{\save +<0pt,8pt>*\txt{$\scriptstyle \bar C_2$}  \restore}
- [ll]
*\cir<2pt>{}
!{\save +<0pt,8pt>*\txt{$\scriptstyle \bar E$}  \restore}
- [ll]
*\cir<2pt>{}
!{\save +<0pt,8pt>*\txt{$\scriptstyle \bar C_1''$}  \restore} } $$
\vglue0pt {\centerline {Subgraph of the dual graph of $g_1$} \vspace{5pt}}
\noindent
where $(g_1)_* [\bar C_1''] = \ell _1 - e^1_1$ and $(g_1)_*[\bar E] = e^1_1$.  
We may now smooth $\bar E \cup \bar C_2$ to a single irreducible component 
$\bar C_2'$ mapped to a curve of anticanonical degree four.  With the usual 
argument of fixing two general points on the image of $\bar C_2'$ we may 
deform the morphism so that $\bar C_2'$ breaks in two components each mapped 
to a curve of anticanonical 
degree two.  The result of this deformation was to replace two components 
mapped to the different divisor classes $\ell _1$ and $\ell _2$ by three 
components mapped to divisor classes of curves of anticanonical degree two.

Iterating the previous argument we deduce that we may assume that condition 
($\star$) holds.

We first treat the case in which there are no components mapped to $\ell$.  
If all the irreducible components of the domain of $g$ are mapped to conics, 
and two of them have images with intersection number at least two, then 
we may use Proposition \ref{struttura} to conclude.  If all the conics in 
the image of $g$ have intersection number at most one, then there must be 
at least three having pairwise intersection products exactly one.  
Otherwise we would be able to find a standard basis 
$\{ \ell , e_1 , \ldots , e_\delta \}$ such that the divisor classes of 
the components of the image of $g$ are in the span of $\ell - e_1$ and 
$\ell - e_2$.  Clearly, no linear combination of these divisors is ample, 
since 
$$ \bigl( \ell - e_1 - e_2 \bigr) \cdot 
\Bigl( a_1 \bigl( \ell - e_1 \bigr) + a_2 \bigl( \ell - e_2 \bigr) \Bigr) = 0 $$
Thus there must be at least three components mapped to divisor classes of 
conics with pairwise intersection products exactly one.  Lemma \ref{piuma} 
allows us to assume that three of these components are adjacent and using 
Proposition \ref{infradito} we conclude.

Suppose now that there is a component mapped to a curve with divisor class 
$\ell$.  Since on a del Pezzo surface of degree at most seven no multiple 
of $\ell $ is ample, it follows that there must be components of the domain 
of $g$ mapped to divisor classes of conics.

Suppose that $\bar C_1$ is mapped to the divisor class $\ell $ and $\bar C_2$ 
is mapped to the divisor class of a conic $Q$.  Thanks to Lemma \ref{piuma} 
we may assume that $\bar C_1$ and $\bar C_2$ are adjacent.  Permuting the 
indices $1, \ldots , \delta $ if necessary, we may assume that the component 
$\bar C_2$ mapped to the conic $Q = a \ell - b_1 e_1 - \ldots - b_\delta e_\delta $ 
has largest possible $b_1$.  Looking at the table (\ref{soluco}), it is easy 
to check that 
\begin{itemize}
\item if $Q$ is of type $H',I,J,K,L$, then $Q+e_1 = -K_{X_8} - K_{X'}$, 
where $X'$ is obtained by contracting a $(-1)-$curve on $X$;
\item if $Q$ is of type $D',F,G,H,I'$, then $Q+e_1 = -K_{X_8} + Q'$, 
where $Q'$ is a conic;
\item if $Q$ is of type $C,D,E$, then $Q+e_1$ is already of the required 
form (for a different choice of standard basis, when $Q=D$ or $E$);
\item if $Q$ is of type $B$, then $Q+\ell $ is of the required form.
\end{itemize}

If $Q$ is of type $B$, then we smooth $\bar C_1 \cup \bar C_2$ to a single 
irreducible component to conclude.  Otherwise, we deform the morphism so that 
$\bar C_1$ breaks into a component mapped to $\ell - e_1$ and a component $e_1$ 
adjacent to $\bar C_2$.  To achieve this splitting, consider the morphism 
$$ \xymatrix @C=35pt 
{{\rm Sl}_{g} (\bar C_1) \ar[r]^{\hspace{10pt} a} & \bar C_2 } $$
and let $\bar p \in \bar C_2$ be a point mapped to the $(-1)-$curve whose 
divisor class is $e_1$ (note that $Q \cdot e_1 \geq 1$).  Since the 
restriction of $g$ to each irreducible component of its domain free, we may 
assume that $g(\bar p)$, as well as the image of the node between $\bar C_1$ 
and $\bar C_2$ are general points of $X$.  Since the morphism $a$ is 
dominant, we may find a deformation $g'$ of $g$ such that $a(g') = \bar p$.  
This means that the ``limiting component'' of $\bar C_1$ breaks in the 
desired way.  Smoothing the union of $\bar C_2$ and the component mapped to 
$e_1$, we obtain a morphism $\bar g' : \bar C_1' \cup \bar C_2' \rightarrow X$ 
where $\bar g'_*[\bar C_1'] = \ell - e_1$ and $\bar g'_*[\bar C_2'] = Q+ e_1$.  
The hypotheses of the theorem imply that 
$\overline \cM _{bir} \bigl( X, \bar g'_*[\bar C_2'] \bigr)$ is irreducible 
since $-K_X \cdot \bigl( Q + e_1 \bigr) = 3$.  Thanks to the previous analysis 
of the divisor class $Q+e_1$, we conclude considering the dominant morphism 
$$ \xymatrix @C=35pt 
{{\rm Sl}_{\bar g'} (\bar C_2') \ar[r]^{\pi \hspace{20pt}} & 
\overline \cM _{bir} \bigl( X, \bar g'_*[\bar C_2'] \bigr) } $$
that we may assume that there is a component of $g$ mapped to the divisor 
class $3\ell - e_1 - \ldots - e_\alpha $, for some $\alpha \leq 8$.

The only remaining case is the one in which the conic $Q$ is of type $A$.  
We may therefore suppose that if a component of the domain of $g$ is mapped 
to the divisor class of a conic, then the divisor class of the image is 
$\ell - e_j$ for some $j$.  Since the image of $g$ is an ample divisor, it 
follows that there must be components $\bar Q_1, \ldots , \bar Q_\delta $ 
mapped to $\ell - e_1, \ldots , \ell - e_\delta $ respectively.

Repeated application of Lemma \ref{piuma} allows us to assume that the 
component mapped to $\ell $ and the two components $\bar Q_1$ $\bar Q_2$ 
are adjacent.  Smoothing the union of these three components to a single 
irreducible free morphism, concludes the first step of the proof.

{\bf Step 2.}  There is a component $\bar C_1$ of $\bar C$ mapped to the 
divisor class $-K_{X_\delta}$, if $\delta \leq 7$.  If $\delta = 8$ (that 
is, the degree of $X = X_8$ is one), then $\bar C_1$ mapped to 
$-rK_{X_8}$, for $r \in \{1,2,3 \}$.  If $r = 1$, then we may choose to 
which rational divisor in $|-K_{X_8}|$ the component $\bar C_1$ maps.

If the component $\bar C_1$ of {\bf Step 1} is mapped to $-rK_{X_8}$, 
there is nothing to prove.

Let $\bar C_1$ be the component of $g$ mapped to the divisor class 
$3 \ell - e_1 - \ldots - e_\alpha $.  If $\alpha = \delta $, then again 
there is nothing to prove.  Suppose therefore that $\alpha < \delta $.  
There is a component of $\bar C$, say $\bar C_2$, such that 
$g_*[\bar C_2] \cdot e_{\alpha +1} \geq 1$, since the image of $g$ is 
an ample divisor; let $C_2 := g_*[\bar C_2]$  Moreover, $C_1 := 
g_*[\bar C_1] = 3\ell - e_1 - \ldots - e_\alpha $ intersects positively 
every non-zero nef divisor, thanks to Corollary \ref{benaco}.  Thus 
$C_1 \cdot C_2 > 0$ and thanks to Lemma \ref{piuma} we may assume that 
$\bar C_1$ and $\bar C_2$ are adjacent in the dual graph of $g$.  Since 
the morphism $g|_{\bar C_2}$ is free, we assume also that $C_2$ meets 
transversely the $(-1)-$curve $E_{\alpha + 1}$ whose divisor class is 
$e_{\alpha + 1}$.  Let $\bar p \in \bar C_2$ be a point such that 
$p := g(\bar p) \in E_{\alpha + 1}$.  Consider the morphism 
$$ \xymatrix @C=35pt 
{{\rm Sl}_{g} (\bar C_1) \ar[r]^{\hspace{10pt} a} & \bar C_2 } $$
and note that it is dominant, since $g|_{\bar C_1}$ is free.  It follows that 
we may find a morphism 
$g_1 : \bar C_1' \cup \bar C_2 \cup \ldots \cup \bar C_r \longrightarrow X$ such 
that $a(g_1) = \bar p$.  We deduce that $g_1(\bar C_1') \ni p$ and 
$(g_1)_* [\bar C_1'] = 3 \ell - e_1 - \ldots - e_\alpha $.  Since 
$\bigl( 3 \ell - e_1 - \ldots - e_\alpha \bigr) \cdot e_{\alpha + 1} = 0$, we 
conclude that $g_1(\bar C_1')$ contains $E_{\alpha + 1}$ and another (irreducible) 
component whose divisor class is $3 \ell - e_1 - \ldots - e_{\alpha + 1}$.  Thus 
the subgraph of the dual graph of $g_1$ spanned by $\bar C_1'$ and 
$\bar C_2$ is 
$$ \xygraph {[] !~:{@{.}} !{<0pt,0pt>;<20pt,0pt>:} 
{\bullet} [rr] {\bullet} [rr] {\bullet} 
*\cir<2pt>{}
!{\save +<0pt,8pt>*\txt{$\scriptstyle \bar C_2$}  \restore}
- [ll]
*\cir<2pt>{}
!{\save +<0pt,8pt>*\txt{$\scriptstyle \bar E_{\alpha + 1}$}  \restore}
- [ll]
*\cir<2pt>{}
!{\save +<0pt,8pt>*\txt{$\scriptstyle \bar C_1''$}  \restore} } $$
\vglue0pt {\centerline {Subgraph of the dual graph of $g_1$} \vspace{5pt}}
\noindent
where $(g_1)_* [\bar C_1''] = 3 \ell - e_1 - \ldots - e_{\alpha + 1}$ and 
$(g_1)_*[\bar E_{\alpha + 1}] = e_{\alpha + 1}$.  We may now smooth 
$\bar E_{\alpha + 1} \cup \bar C_2$ to a single irreducible component 
$\bar C_2'$ mapped to a curve of anticanonical degree three or four.  If 
this new component has degree four, we break it into two components of 
anticanonical degree two.

If the degree of $X$ is at least two, iterating this procedure allows us 
to produce a component of $\bar C$ whose image represents the divisor 
class $-K_X$.  If the degree of $X$ is one, we may apply the same 
procedure to obtain a component $\bar C_1$ mapped to $-K_{X_8}$, but we 
still have to prove that we may choose which nodal rational divisor in 
$|-K_{X_8}|$ is in the image of the morphism.

If the component $\bar C_2$ adjacent to the component $\bar C_1$ has 
degree two, we smooth these two components to a single irreducible one 
and using the irreducibility of 
$\overline \cM _{bir} \bigl( X, \beta \bigr)$, for 
$-K_{X_8} \cdot \beta=3$ we conclude.

If the component $\bar C_2$ adjacent to the component $\bar C_1$ mapped 
to $-K_{X_8}$ has degree three, then it represents one of the five 
divisor classes $-3K_{X_8}$, $\bigl( -K_{X_8} - K_{X_7} \bigr)$, 
$\bigl( -K_{X_8} + Q \bigr)$, $-K_{X_6}$ or $\ell $, where $Q$ is a conic, 
$-K_{X_7}$ and $-K_{X_6}$ are del Pezzo surfaces of degree two and three 
respectively dominated by $X_8$.

In the first three cases, we deform the morphism so that the component 
$\bar C_2$ breaks into a component mapped to a preassigned nodal divisor 
in $|-K_{X_8}|$ and into a component where the morphism is free.  In these 
cases, smoothing the component $\bar C_1$ with the component adjacent to 
it into which $\bar C_2$ broke finishes the proof.

If $\bar C_2$ is mapped to $-K_{X_6}$ or $\ell $, then we may choose a 
standard basis so that $-K_{X_6} = 3\ell - e_1 - \ldots - e_6$.  The 
morphism $\varphi : X_8 \rightarrow \bP^2$ determined by the linear 
system $|\ell |$ is the contraction of the $(-1)-$curves with divisor 
classes $e_1$, \ldots , $e_8$ to the points $q_1, \ldots , q_8 \in \bP^2$.  
The image of $\bar C_1$ in $\bP^2$ is a nodal plane cubic through the 
eight points $q_1$, \ldots , $q_8$.  The image of $\bar C_2$ is either 
a rational cubic through $q_1$, \ldots , $q_6$ or a line.  We treat only 
the case in which the image of $\bar C_2$ is a nodal cubic, since the 
other one is simpler and the arguments are similar.

Deform the nodal cubic through $q_1$, \ldots , $q_6$ until it contains 
a general point $q \in \bP^2$.  We may now slide the node between $\bar C_1$ 
and $\bar C_2$ along $\bar C_1$ until it reaches a point on $\bar C_1$ mapped 
to the point $q_7 \in \bP^2$.  As we slide the node, we let the image of 
$\bar C_2$ always contain the general point $q$.  When the deformation is 
finished, the component $\bar C_2$ breaks as the $(-1)-$curve with divisor 
class $e_7$ and the divisor class $-K_{X_6} - e_7$.  Since the point $q$ is 
general, we know that there are only finitely many (in fact twelve) possible 
configurations for these limiting positions and we may assume that they are 
all transverse to the image of $\bar C_1$.  Thus the dual graph of the 
resulting morphism 
$g' : \bar C_1 \cup \bar C_2' \cup \bar E_7 \longrightarrow X_8$ 
is 
$$ \xygraph {[] !~:{@{.}} !{<0pt,0pt>;<20pt,0pt>:} 
{\bullet} [rr] {\bullet} [rr] {\bullet} 
*\cir<2pt>{}
!{\save +<0pt,8pt>*\txt{$\scriptstyle \bar C_2'$}  \restore}
- [ll]
*\cir<2pt>{}
!{\save +<0pt,8pt>*\txt{$\scriptstyle \bar E_7$}  \restore}
- [ll]
*\cir<2pt>{}
!{\save +<0pt,8pt>*\txt{$\scriptstyle \bar C_1$}  \restore} } $$
\vglue0pt {\centerline {Dual graph of $g'$} \vspace{5pt}}

We may smooth the components $\bar C_1 \cup \bar E_7$ to a unique irreducible 
component mapped to a curve of anticanonical degree two.  The assumptions of 
the theorem imply that $\overline \cM _{bir} \bigl( X_8, C_1 + e_7 \bigr)$ is 
irreducible and we may therefore deform the morphism so that its domain breaks 
as a preassigned rational nodal divisor $C''$ in $|-K_{X_8}|$ and a curve 
mapped to the divisor class $e_7$.  The dual graph of the resulting morphism 
$g''$ is of one of the following types:
$$ \xygraph {[] !~:{@{.}} !{<0pt,0pt>;<20pt,0pt>:} 
{\bullet} [rr] {\bullet} [rr] {\bullet} 
*\cir<2pt>{}
!{\save +<0pt,8pt>*\txt{$\scriptstyle \bar C_2'$}  \restore}
- [ll]
*\cir<2pt>{}
!{\save +<0pt,8pt>*\txt{$\scriptstyle \bar E_7$}  \restore}
- [ll]
*\cir<2pt>{}
!{\save +<0pt,8pt>*\txt{$\scriptstyle \bar C_1''$}  \restore} } 
\hspace{30pt}
\xygraph {[] !~:{@{.}} !{<0pt,0pt>;<20pt,0pt>:} 
{\bullet} [rr] {\bullet} [rr] {\bullet} 
*\cir<2pt>{}
!{\save +<0pt,8pt>*\txt{$\scriptstyle \bar C_2'$}  \restore}
- [ll]
*\cir<2pt>{}
!{\save +<0pt,8pt>*\txt{$\scriptstyle \bar C_1''$}  \restore}
- [ll]
*\cir<2pt>{}
!{\save +<0pt,8pt>*\txt{$\scriptstyle \bar E_7$}  \restore} } $$
\vglue0pt {\centerline {Possible dual graphs of $g''$} \vspace{5pt}}

In the first case we smooth $\bar E_7 \cup \bar C_2'$ to a single 
irreducible component and conclude.  In the second case, we slide the node 
between $\bar C_1''$ and $\bar C_2'$ until it reaches the node between 
$\bar C_1'$ and $\bar E_7$, in such a way that the limiting position of 
$\bar C_2'$ does not coincide with the image of $\bar C_1''$ (we can do this 
thanks to the irreducibility of 
$\overline \cM _{bir} \bigl( X_8, C_1'' + C_2' \bigr)$).  It follows that 
the dual graph of the morphism $\bar g$ thus obtained is 
$$ \xygraph {[] !~:{@{.}} !{<0pt,0pt>;<20pt,0pt>:} 
{\bullet} [rr] {\bullet} [ur] {\bullet} [dd] {\bullet} 
*\cir<2pt>{}
!{\save +<9pt,0pt>*\txt{$\scriptstyle \bar C_2''$}  \restore}
- [ul] - [ur]
*\cir<2pt>{}
!{\save +<8pt,0pt>*\txt{$\scriptstyle \bar E_7$}  \restore}
 [dl] 
!{\save +<-3pt,8pt>*\txt{$\scriptstyle \bar E$}  \restore} 
- [ll]
*\cir<2pt>{}
!{\save +<0pt,8pt>*\txt{$\scriptstyle \bar C_1''$}  \restore} }$$
\vglue0pt {\centerline {Dual graph of $\bar g$} \vspace{5pt}}

It is easy to check that $\bar g$ represents a smooth point of 
$\overline \cM _{0,0} \bigl( X_8, C_1'' + C_2'' + e_7 \bigr)$, since 
the sheaf $\bar g ^* \cT_{X_8}$ is globally generated on $\bar C_2''$ and 
has no first cohomology group on the remaining components, thanks to 
Lemmas \ref{blodo} and \ref{tretre}.  We may now smooth the components 
$\bar E_7 \cup \bar E \cup \bar C_2''$ to a single irreducible component 
on which the morphism is free to conclude.

A similar and simpler argument proves the same result if $\bar C_2$ has 
divisor class $\ell $.  This finishes the proof of this step.

{\bf Step 3.}  We may deform $g : \bar C \rightarrow X$ to a morphism 
$h : \bar D_1 \cup \bar D_2 \rightarrow X$ where $\bar D_1$ and $\bar D_2$ 
are irreducible, $h_*[\bar D_1] = -rK_X$ ($r \in \{1,2,3\}$), 
$h(\bar D_1) \neq h(\bar D_2)$ and $h|_{\bar D_2}$ is free.  If the 
degree of $X$ is at least two, then $r=1$.  If the degree of $X$ is one 
and $r = 1$, we may choose which rational divisor in $|-K_X|$ $h(\bar D_1)$ 
is.  Note that we are not requiring $h|_{\bar D_2}$ to be birational to its 
image.

Thanks to the previous steps, we may assume that $g_* [\bar C_1] = -rK_X$ 
(with the required restriction for $r$) and that all the components of 
$\bar C$ different from $\bar C_1$ are immersed to curves of anticanonical 
degree two or three.  Let $C_2$, \ldots , $C_r$ be the components of the 
image of $g$ different from $g(\bar C_1)$, and let $\bar C_i$ be the 
component of $\bar C$ whose image is $C_i$.

The divisor class $C_2 + \ldots + C_r$ is nef and if it is not a multiple 
of a conic, then it meets all nef curves positively, thanks to Corollary 
\ref{benaco}.  Thus, still assuming that $C_2 + \ldots + C_r$ is not a 
multiple of a conic, we may deform the morphism using Lemma \ref{piuma} 
and assume that the union of all the components of the domain of $g$ 
different from $\bar C_1$ is connected.  Smoothing the resulting union 
$\bar C_2 \cup \ldots \cup \bar C_r$ concludes the proof of this step in 
this case.

Suppose that $C_2 + \ldots + C_r$ is a multiple of a conic.  Then it follows 
that $C_2 = \ldots = C_r = C$ is a conic.  Since $g$ is birational to its 
image and two divisors linearly equivalent to the same conic are either 
disjoint or they coincide, it follows that the dual graph of $g$ must be 
$$ \xygraph {[] !~:{@{.}} !{<0pt,0pt>;<20pt,0pt>:} 
{\bullet} [ul] {\bullet} [rr] {\bullet} [d] {\vdots} [d] {\bullet} [ll] {\bullet} 
*\cir<2pt>{}
!{\save +<-8pt,0pt>*\txt{$\scriptstyle \bar C_r$}  \restore}
- [ur] - [dr] 
*\cir<2pt>{}
!{\save +<13pt,0pt>*\txt{$\scriptstyle \bar C_{r-1}$}  \restore}
 [uu] 
*\cir<2pt>{}
!{\save +<8pt,0pt>*\txt{$\scriptstyle \bar C_3$}  \restore}
- [dl] - [ul] 
*\cir<2pt>{}
!{\save +<-8pt,0pt>*\txt{$\scriptstyle \bar C_2$}  \restore} 
 [dr] 
!{\save +<0pt,9pt>*\txt{$\scriptstyle \bar C_1$}  \restore} } $$
\vglue0pt {\centerline {Dual graph of $g$} \vspace{5pt}}

By sliding each $\bar C_i$ along $\bar C_1$, we may assume that 
the images of all the components mapped to a conic coincide and 
that the nodes in the source curve all map to the same general 
point $p \in X$.  Thus the dual graph of the resulting morphism 
$g''$ is 
$$ \xygraph {[] !~:{@{.}} !{<0pt,0pt>;<20pt,0pt>:} 
{\bullet} [rr] {\bullet} [ul] {\bullet} [rr] {\bullet} [d] {\vdots} [d] 
{\bullet} [ll] {\bullet} 
*\cir<2pt>{}
!{\save +<8pt,0pt>*\txt{$\scriptstyle \bar C_r$}  \restore}
- [ur] - [dr] 
*\cir<2pt>{}
!{\save +<13pt,0pt>*\txt{$\scriptstyle \bar C_{r-1}$}  \restore}
 [uu] 
*\cir<2pt>{}
!{\save +<8pt,0pt>*\txt{$\scriptstyle \bar C_3$}  \restore}
- [dl] - [ul] 
*\cir<2pt>{}
!{\save +<8pt,0pt>*\txt{$\scriptstyle \bar C_2$}  \restore} 
 [dr] 
!{\save +<0pt,9pt>*\txt{$\scriptstyle \bar E$}  \restore} 
- [ll] 
!{\save +<0pt,9pt>*\txt{$\scriptstyle \bar C_1$}  \restore} } $$
\vglue0pt {\centerline {Dual graph of $g''$} \vspace{5pt}}
\noindent
where $\bar E$ is a contracted component whose image is $p$.  
The morphism $g''$ is a smooth point of 
$\overline \cM _{bir} \bigl( X, g_* [\bar C] \bigr)$ since the sheaf 
$(g'') ^* \cT_X$ is globally generated on each irreducible component 
of the domain curve of $g''$.

We may smooth all the components 
$\bar E \cup \bar C_2 \cup \ldots \cup \bar C_r$ to a single 
irreducible component which represents a multiple cover (in fact a 
degree $r-1$ cover) of its image.  The resulting morphism 
$h : \bar D_1 \cup \bar D_2 \rightarrow X$ is therefore such that 
the image of $\bar D_1$ is a rational divisor in the linear system 
$|-rK_X|$ ($r=1$ unless the degree of $X$ is one, in which case 
$r \leq 3$), which is an arbitrary preassigned one in case 
$\deg X = 1$ and $r=1$, and the morphism $h|_{\bar D_2}$ is a 
multiple cover of the divisor class of a conic.  This concludes the 
proof of the third step.

We may write $h_*[\bar D_2] = n_8 (-K_{X_8}) + \ldots + n_2 (-K_2) + 
D_2'$ as in Corollary \ref{maquale} (to simplify the notation we will 
assume that $\deg X = 1$; if this is not the case, simply set to zero 
all the coefficients $n_\alpha $, with $\alpha > 9 - \deg X$).  Let 
$n := [\frac {n_8} {2}]$, if $n_8 \neq 1$ and let $n = 1$, if $n_8 = 1$.  
Thus we have $n_8 = 2(n-1) + 3$, if $n_8$ is odd and at least three, 
and $n_8 = 2n$, if $n$ is even.

{\bf Step 4.}  Let $S \subset X$ be a nodal rational divisor in the 
linear system $|-K_X|$.  We may deform $h$ to a morphism 
$k : \bar K := \bar K_1 \cup \ldots \cup \bar K_\ell \longrightarrow X$ 
with the following properties:
\renewcommand{\labelenumi}{P\arabic{enumi})}
\begin{enumerate}
\item \label{aringoli} the morphism $k$ restricted to each irreducible 
component of $\bar K$ is free, except possibly on $\bar K_1$;
\item each irreducible component of $\bar K$ represents one of the 
divisor classes $-3K_{X_8}$, $-2K_{X_8}$, $-K_{X_7}$, \ldots , 
$-K_{X_2}$, $D_2'$, except $\bar K_1$, whose image may also be $S$;
\item the dual graph of $k$ is 
$$ \xygraph {[] !~:{@{.}} !{<0pt,0pt>;<20pt,0pt>:} 
{\bullet} [rr] {\bullet} [rrrr] {\bullet} [rr] {\bullet} 
*\cir<2pt>{}
!{\save +<0pt,8pt>*\txt{$\scriptstyle \bar K_\ell $}  \restore}
- [ll] 
*\cir<2pt>{}
!{\save +<0pt,8pt>*\txt{$\scriptstyle \bar K_{\ell -1}$}  \restore}
- [ll] {\cdots} -[ll]
*\cir<2pt>{}
!{\save +<0pt,8pt>*\txt{$\scriptstyle \bar K_2$}  \restore}
- [ll] 
*\cir<2pt>{}
!{\save +<0pt,8pt>*\txt{$\scriptstyle \bar K_1$}  \restore} } $$
\vglue0pt {\centerline {Dual graph of $k$} \vspace{5pt}}
\item the component $\bar K_1$ is mapped to $-3K_{X_8}$ if $n_8$ is 
odd and at least three, to $S$ if $n_8 = 1$ and to $-2K_{X_8}$ if 
$n_8$ is even and bigger than zero;
\item the components 
$\bar K_2$, $\bar K_3$, \ldots , $\bar K_n $ are mapped to $-2K_{X_8}$;
\item let $N_\alpha := n+n_7 + \ldots + n_{\alpha + 1}$; the 
components $\bar K_{N_\alpha + 1}$, \ldots , $\bar K_{N_{\alpha - 1} }$ 
are mapped to $-K_{X_\alpha }$;
\item if $D_2' \neq 0$, then $\bar K_\ell $ is mapped to $D_2'$;
\item the morphism $k|_{\bar K_1 \cup \ldots \cup \bar K_{\ell -1}}$ 
is birational to its image. \label{testa}
\end{enumerate}
\renewcommand{\labelenumi}{\arabic{enumi})}

We call a morphism satisfying all the above properties a 
{\it morphism in standard form}.

By induction on the anticanonical degree of the divisor, we know that 
the space $\overline \cM _{bir} \bigl( X, h_* [\bar D_2] \bigr)$ is 
irreducible (or empty if $h_* [\bar D_2]$ is a multiple of a conic and 
it is clear how to proceed in this case; we will not mention this issue 
anymore).  We may therefore deform the morphism $h|_{\bar D_2}$ to a 
morphism $l : \bar E := \bar E_1 \cup \ldots \cup \bar E_t \longrightarrow X$ 
in standard form.  Considering the morphism 
$$ \xymatrix @C=40pt 
{{\rm Sl}_{h} (\bar D_2) \ar[r]^{\pi \hspace{20pt}} & 
\overline \cM _{bir} \bigl( X, h_*[\bar D_2] \bigr) } $$
we may find a deformation $\widetilde l$ of $h$ such that 
$\pi (\widetilde l ) = l$.  The dual graph of $\widetilde l$ is 
$$ \xygraph {[] !~:{@{.}} !{<0pt,0pt>;<20pt,0pt>:} 
{\bullet} [rr] {\bullet} [rrrr] {\bullet} [d] {\bullet} [urrrr] {\bullet} [rr] {\bullet} 
*\cir<2pt>{}
!{\save +<0pt,8pt>*\txt{$\scriptstyle \bar E_t$}  \restore}
- [ll] 
*\cir<2pt>{}
!{\save +<0pt,8pt>*\txt{$\scriptstyle \bar E_{t-1}$}  \restore}
- [ll] {\cdots} -[ll]
*\cir<2pt>{}
!{\save +<0pt,8pt>*\txt{$\scriptstyle \bar E_j$}  \restore}
- [d] 
*\cir<2pt>{}
!{\save +<8pt,0pt>*\txt{$\scriptstyle \bar D_1$}  \restore}
 [u] - [ll] {\cdots} -[ll]
*\cir<2pt>{}
!{\save +<0pt,8pt>*\txt{$\scriptstyle \bar E_2$}  \restore}
- [ll] 
*\cir<2pt>{}
!{\save +<0pt,8pt>*\txt{$\scriptstyle \bar E_1$}  \restore} } $$
\vglue0pt {\centerline {Dual graph of $\widetilde l$} \vspace{5pt}}
\noindent
for some $ 1 \leq j \leq t$.  We want to show by induction on 
$j$ that we may assume that $j = 1$.  In the case $j=1$ there is 
nothing to prove.

Suppose $j \geq 2$ and consider the morphism 
$$ \xymatrix @C=35pt 
{{\rm Sl}_{\widetilde l} (\bar E_{j-1}) \ar[r]^{\hspace{15pt} a} & \bar E_j } $$

Unless $j=2$ and $\bar E_1$ represents $|-K_{X_8}|$, the morphism $a$ 
is dominant and we may find a morphism $l_1$ such that $a(l_1)$ is 
the node between $\bar D_1$ and $\bar E_j$.  The dual graph of the 
morphism $l_1$ is 
$$ \xygraph {[] !~:{@{.}} !{<0pt,0pt>;<20pt,0pt>:} 
{\bullet} [rr] {\bullet} [rrrr] {\bullet} [rr] {\bullet} [d] {\bullet} [urr] 
{\bullet} [rrrr] {\bullet} [rr] {\bullet} 
*\cir<2pt>{}
!{\save +<0pt,8pt>*\txt{$\scriptstyle \bar E_t$}  \restore}
- [ll] 
*\cir<2pt>{}
!{\save +<0pt,8pt>*\txt{$\scriptstyle \bar E_{t-1}$}  \restore}
- [ll] {\cdots} -[ll]
*\cir<2pt>{}
!{\save +<0pt,8pt>*\txt{$\scriptstyle \bar E_j$}  \restore}
- [ll] 
*\cir<2pt>{}
!{\save +<0pt,8pt>*\txt{$\scriptstyle \bar E$}  \restore}
- [d] 
*\cir<2pt>{}
!{\save +<8pt,0pt>*\txt{$\scriptstyle \bar D_1$}  \restore}
 [u] 
- [ll] 
*\cir<2pt>{}
!{\save +<6pt,6pt>*\txt{$\scriptstyle \bar E_{j-1}'$}  \restore}
- [ll] {\cdots} -[ll]
*\cir<2pt>{}
!{\save +<0pt,8pt>*\txt{$\scriptstyle \bar E_2$}  \restore}
- [ll] 
*\cir<2pt>{}
!{\save +<0pt,8pt>*\txt{$\scriptstyle \bar E_1$}  \restore} } $$
\vglue0pt {\centerline {Dual graph of $l_1$} \vspace{5pt}}
\noindent
and the component $\bar E$ is contracted.  The morphism $l_1$ represents 
a smooth point of $\overline \cM _{bir} \bigl( X, h_*[\bar C] \bigr)$, 
since $l_1^* \cT_X$ is globally generated on all components different 
from $\bar D_1$, and ${\rm H}^1 \bigl( \bar D_1 , l_1^* \cT_X \bigr) = 0$.

The morphism $l_1$ is also a limit of morphisms with dual graph 
$$ \xygraph {[] !~:{@{.}} !{<0pt,0pt>;<20pt,0pt>:} 
{\bullet} [rr] {\bullet} [rrrr] {\bullet} [rr] {\bullet} [d] {\bullet} [urr] 
{\bullet} [rrrr] {\bullet} [rr] {\bullet} 
*\cir<2pt>{}
!{\save +<0pt,8pt>*\txt{$\scriptstyle \bar E_t$}  \restore}
- [ll] 
*\cir<2pt>{}
!{\save +<0pt,8pt>*\txt{$\scriptstyle \bar E_{t-1}$}  \restore}
- [ll] {\cdots} -[ll]
*\cir<2pt>{}
!{\save +<0pt,8pt>*\txt{$\scriptstyle \bar E_j$}  \restore}
- [ll] 
*\cir<2pt>{}
!{\save +<6pt,6pt>*\txt{$\scriptstyle \bar E_{j-1}''$}  \restore}
- [d] 
*\cir<2pt>{}
!{\save +<8pt,0pt>*\txt{$\scriptstyle \bar D_1$}  \restore}
 [u] 
- [ll] 
*\cir<2pt>{}
!{\save +<6pt,6pt>*\txt{$\scriptstyle \bar E_{j-2}$}  \restore}
- [ll] {\cdots} -[ll]
*\cir<2pt>{}
!{\save +<0pt,8pt>*\txt{$\scriptstyle \bar E_2$}  \restore}
- [ll] 
*\cir<2pt>{}
!{\save +<0pt,8pt>*\txt{$\scriptstyle \bar E_1$}  \restore} } $$
\vglue0pt {\centerline {Dual graph of morphisms limiting to $l_1$} \vspace{5pt}}

We can apply the induction hypothesis to these last morphisms 
to conclude that we may deform the morphism $f$ to a morphism 
$m : \bar F := \bar D_1 \cup \bar F_1 \cup \ldots \cup \bar F_t \longrightarrow X$, 
where $m|_{\bar F_1 \cup \ldots \cup \bar F_t}$ is a morphism in 
standard form and $m|_{\bar D_1}$ is birational onto its image 
and the image is a rational divisor in $|-K_X|$, or $|-rK_X|$ 
if $\deg X = 1$.  We may specify which rational curve $m(\bar D_1)$ 
is if it represents $-K_{X_8}$, and we assume it is $S$ if and only 
if $\bar F_1$ is not mapped to $S$.  If $\bar F_1$ represents a 
divisor class different from $|-K_{X_8}|$, then we may assume that 
$\bar D_1$ is adjacent to $\bar F_1$ and conclude; if $\bar F_1$ 
represents the divisor class $|-K_{X_8}|$, then $\bar D_1$ is 
adjacent to $\bar F_2$ and $\bar D_1$ represents the divisor class 
$-rK_{X_8}$ for some $r \in \{ 1,2,3 \}$.  The divisor represented 
by $\bar F_2$ is either $-K_{X_\alpha}$, for some $\alpha \leq 7$ or 
it is a divisor on a del Pezzo surface of degree eight, i.e. 
$\bP^1 \times \bP^1$ or $Bl_p(\bP^2)$.

We are only going to treat the case $m_*[\bar F_2] = -K_{X_7}$, and 
$m_*[\bar D_1] = -K_{X_8}$.  The remaining cases are simpler and can 
be treated with similar techniques.

Let $J \subset \overline \cM _{bir} \bigl( X, h_*[\bar C] \bigr)$ 
be the closure of the set of morphisms 
$m' : \bar D_1' \cup \bar F_1' \cup \ldots \cup \bar F_t' \longrightarrow X_8$ 
such that $m'|_{\bar F_2' \cup \ldots \cup \bar F_t'} \simeq 
m|_{\bar F_2 \cup \ldots \cup \bar F_t}$, $m'|_{\bar D_1'} \simeq m|_{\bar D_1}$ 
the image of $\bar F_2'$ is a general rational divisor in $|-K_{X_7}|$ 
and the dual graphs of $m$ and $m'$ coincide.

We clearly have a morphism $J \rightarrow 
\overline \cM _{bir} \bigl( X, -K_{X_7} \bigr)$ obtained by 
``restricting a morphism in $J$ to its $\bar F_2'$ component.''  
Since the intersection number $m_*[\bar F_2] \cdot m_*[\bar D_1]$ 
equals two, it follows that $J$ has at most two irreducible 
components.  Moreover, even if the space $J$ is reducible, its 
components meet.  To see this, we construct a point in common to 
the two components.  Let $\psi : X_8 \longrightarrow \bP^2$ be 
the morphism induced by the linear system $|-K_{X_7}|$.  The 
morphism $\psi $ contracts a $(-1)-$curve and ramifies above a 
smooth plane quartic $R$.  The images of $\bar D_1$ and 
$\bar F_1$ are two distinct tangent lines in $\bP^2$ which 
contain the image $c$ of the contracted component and are tangent 
to $R$, but are not bitangent lines to $R$.  The images of 
the curves $\bar F_2'$ are tangent lines to $R$.  Since $R$ has 
degree four, it follows that there is a point $c' \in R$ where 
the image of $\bar F_1$ meets transversely $R$.  Through such a 
point $c'$, there are ten tangent lines to $R$, different from 
the tangent line to $R$ at $c'$.  Each of these lines corresponds 
to a point in the intersection of the two components of $J$.  
Moreover, these points in common to the components are easily seen 
to be smooth points of the mapping space, using Proposition 
\ref{grafico}.  We deduce that $J$ is connected, and thus in the 
same irreducible component of 
$\overline \cM _{bir} \bigl( X, h_*[\bar C] \bigr)$ containing $m$ 
there is a morphism $m_1$ with dual graph 
$$ \xygraph {[] !~:{@{.}} !{<0pt,0pt>;<20pt,0pt>:} 
{\bullet} [rr] {\bullet} [d] {\bullet} [urr] 
{\bullet} [d] {\bullet} [urr] {\bullet} [rrrr] {\bullet} [rr] {\bullet} 
*\cir<2pt>{}
!{\save +<0pt,8pt>*\txt{$\scriptstyle \bar F_t$}  \restore}
- [ll] 
*\cir<2pt>{}
!{\save +<0pt,8pt>*\txt{$\scriptstyle \bar F_{t-1}$}  \restore}
- [ll] {\cdots} -[ll]
*\cir<2pt>{}
!{\save +<0pt,8pt>*\txt{$\scriptstyle \bar F_3$}  \restore}
- [ll] 
*\cir<2pt>{}
!{\save +<0pt,8pt>*\txt{$\scriptstyle \bar F_2'$}  \restore}
- [d] 
*\cir<2pt>{}
!{\save +<8pt,0pt>*\txt{$\scriptstyle \bar F_2''$}  \restore}
 [u] 
- [ll] 
*\cir<2pt>{}
!{\save +<0pt,8pt>*\txt{$\scriptstyle \bar E$}  \restore}
- [d] 
*\cir<2pt>{}
!{\save +<8pt,0pt>*\txt{$\scriptstyle \bar D_1$}  \restore}
 [u] 
- [ll] 
*\cir<2pt>{}
!{\save +<0pt,8pt>*\txt{$\scriptstyle \bar F_1$}  \restore} } $$
\vglue0pt {\centerline {Dual graph of $m_1$} \vspace{5pt}}
\noindent
which agrees with $m$ on the components with the same label and such 
that $\bar F_2''$ is mapped to the divisor class contracted by $\psi$ 
and $\bar F_2'$ is mapped to a rational divisor in $|-K_{X_8}|$ 
distinct from the images of $\bar D_1$ and $\bar F_1$.  It follows 
from the computations in {\bf Step 6, Case 4} below that $m_1$ 
represents a smooth point of the mapping space.  We may now smooth 
$\bar D_1 \cup \bar E \cup \bar F_1$ to a single irreducible component 
$\bar G_1$ representing a nodal rational divisor in $|-2K_{X_8}|$.  
Similarly, we may smooth $\bar F_2' \cup \bar F_2''$ to an irreducible 
component $\bar G_2$ representing a nodal rational divisor in 
$|-K_{X_7}|$.  The resulting morphism $m_2$ has dual graph 
$$ \xygraph {[] !~:{@{.}} !{<0pt,0pt>;<20pt,0pt>:} 
{\bullet} [rr] {\bullet} [rr] {\bullet} [rrrr] {\bullet} [rr] {\bullet} 
*\cir<2pt>{}
!{\save +<0pt,8pt>*\txt{$\scriptstyle \bar F_t$}  \restore}
- [ll] 
*\cir<2pt>{}
!{\save +<0pt,8pt>*\txt{$\scriptstyle \bar F_{t-1}$}  \restore}
- [ll] {\cdots} -[ll]
*\cir<2pt>{}
!{\save +<0pt,8pt>*\txt{$\scriptstyle \bar F_3$}  \restore}
- [ll] 
*\cir<2pt>{}
!{\save +<0pt,8pt>*\txt{$\scriptstyle \bar G_2$}  \restore}
- [ll] 
*\cir<2pt>{}
!{\save +<0pt,8pt>*\txt{$\scriptstyle \bar G_1$}  \restore} } $$
\vglue0pt {\centerline {Dual graph of $m_2$} \vspace{5pt}}
\noindent
This concludes the proof of this step.

We now define a locally closed subset $K_\beta $ of 
$\overline \cM _{0,0} \bigl( X, \beta \bigr)$.  Let $\bar K_\beta $ be 
the closure of the set of morphisms in standard form.  The subspace 
$K_\beta \subset \bar K_\beta $ is the open subset of points lying in 
a unique irreducible component of 
$\overline \cM _{0,0} \bigl( X, \beta \bigr)$, or equivalently 
$K_\beta $ is the complement in $\bar K_\beta $ of the union of all 
the pairwise intersections of the irreducible components of 
$\overline \cM _{0,0} \bigl( X, \beta \bigr)$.  In particular, all the 
points of $\bar K_\beta $ that are smooth points in 
$\overline \cM _{0,0} \bigl( X, \beta \bigr)$ lie in $K_\beta $.

{\bf Step 5.}  The morphisms in standard form are contained in $K_\beta $.

It is enough to prove that a morphism in standard form is a smooth point 
of $\overline \cM _{0,0} \bigl( X, \beta \bigr)$.  Let 
$k : \bar K \rightarrow X$ be a morphism in standard form and let 
$\bar K_1$, \ldots , $\bar K_\ell $ be the components of $\bar K$.  We 
will always assume that the numbering of the components is the 
``standard'' one.  The morphism $k$ represents a smooth point of 
$\overline \cM _{0,0} \bigl( X, \beta \bigr)$ since 
$k^* \cT_X$ is globally generated on all components $\bar K_i$, for $i \geq 2$ 
and ${\rm H}^1 \bigl( \bar K_1 , k^* \cT_X \bigr) = 0$.

{\bf Step 6.}  The space $K_\beta $ is connected.

To prove connectedness of $K_\beta $, let $k : \bar K \rightarrow X$ be a 
morphism in standard form and suppose that all the nodes of $\bar K$ are 
mapped to points of $X$ not lying on $(-1)-$curves.  There are such 
morphisms in all the connected components of $K_\beta$ since 
$k|_{\bar K_i}$ is a free morphism, for $i \geq 2$ and $S$ is not a 
$(-1)-$curve.  Given any $k' : \bar K' \rightarrow X$, we construct a 
deformation from $k'$ to $k$ entirely contained in $K_\beta $.  This is 
clearly enough to prove the connectedness of $K_\beta$.

We are going to construct the deformation in stages.

We prove that in the same connected component of $K_\beta $ 
containing $k'$ there is a morphism 
$k_1 : \bar K_1^1 \cup \ldots \cup \bar K_\ell ^1 \longrightarrow X$ such 
that $k_1|_{\bar K_1^1} \simeq k|_{\bar K_1}$.

This is true by assumption if $k_*[\bar K_1] = -K_{X_8}$, since in this 
case $k(\bar K_1) = S = k'(\bar K_1')$ and $k$ and $k'$ are birational.  
Thus in this case we may choose $k_1 = k'$.

Suppose that $k_*[\bar K_1] \neq -K_{X_8}$.  
Since $k'|_{\bar K_i}$ is free for all $i$'s, we may assume that 
$k' (\bar K_i)$ is not contained in the image of $k$, for all $i$'s.  
Thanks to Theorem \ref{maschera}, Theorem \ref{cabala} and Theorem 
\ref{barbapapa}, we conclude that there is an irreducible 
curve $P \subset K_{k_*[\bar K_1]} \subset 
\overline \cM _{0,0} \bigl( X, k_*[\bar K_1] \bigr)$ 
containing $k'|_{\bar K_1'}$ and $k|_{\bar K_1}$.  Consider the morphism 
$$ \xymatrix @C=35pt { 
{\rm Sl}_{k'} \bigl( \bar K_1' \bigr) \ar[r] ^{\pi \hspace{20pt}} & 
\overline \cM _{0,0} \Bigl( X, k_* [\bar K_1] \Bigr) } $$
and let $\bar P \subset \pi ^{-1} (P)$ be an irreducible curve 
dominating $P$ and containing $k'$.  The curve $\bar P$ has finitely many 
points not lying in $K_\beta $: they are the points $\tilde k$ for which 
the image of $\pi (\tilde k)$ contains a component of 
$k'(\bar K_2' \cup \ldots \cup \bar K_\ell ')$.  By construction, $\bar P$ 
contains a morphism 
$k_1 : \bar K_1^1 \cup \ldots \cup \bar K_\ell ^1 \longrightarrow X$ such 
that $k_1|_{\bar K_1^1} \simeq k|_{\bar K_1}$ and 
$k_1|_{\bar K_2^1 \cup \ldots \cup \bar K_\ell ^1} \simeq 
k'|_{\bar K_2' \cup \ldots \cup \bar K_\ell '}$.  It follows that 
$\bar P \cap K_{\beta }$ is an irreducible curve contained in $K_{\beta }$ 
and containing $k'$ and $k_1$.  Therefore $k'$ and $k_1$ are in the 
same connected (in fact irreducible) component of $K_\beta $.

Thus to prove that $K_\beta $ is connected we may assume that 
$k'|_{\bar K_1'} \simeq k|_{\bar K_1}$.  Suppose that we found a morphism 
$k_{j-1}$ in the same connected component of $K_\beta $ as $k'$ such that 
$k_{j-1}|_{\bar K_1^{j-1} \cup \ldots \cup \bar K_{j-1}^{j-1}} \simeq 
k|_{\bar K_1 \cup \ldots \cup \bar K_{j-1}}$ for some $2 \leq j \leq \ell $.  
If we can find a connected subset of $K_\beta $ containing $k_{j-1}$ and a 
morphism 
$k_j : \bar K_1^j \cup \ldots \cup \bar K_\ell ^j \longrightarrow X$ such 
that $k_j|_{\bar K_1^j \cup \ldots \cup \bar K_j^j} \simeq 
k|_{\bar K_1 \cup \ldots \cup \bar K_j}$, then we may conclude by induction 
on $j$.

The remaining part of the proof will examine the several cases separately.  
To simplify the notation we assume that $j=2$, we write $\bar K_1$ 
also for $\bar K_1^1$ since $k|_{\bar K_1} \simeq k_2|_{\bar K_1^1}$, and 
we let 
$k = k|_{\bar K_1 \cup \bar K_2}$ and $k_1 = k_1|_{\bar K_1 \cup \bar K_2^1}$; 
to get the result for $k$ and $k_1$ simply consider the morphism 
$$ \xymatrix @C=35pt { 
{\rm Sl}_{k_1} \bigl( \bar K_2^1 \bigr) \ar[r] ^{\pi \hspace{20pt}} & 
\overline \cM _{0,0} \Bigl( X, k_* [\bar K_2] \Bigr) } $$
and lift the path in $\overline \cM _{0,0} \Bigl( X, k_* [\bar K_2] \Bigr)$ 
to a path in $\overline \cM _{0,0} \Bigl( X, k_* [\bar K] \Bigr)$, and note 
that the lift lies in the space $K_\beta $.

\noindent
{\bf \boldmath Case 1: $k_*[\bar K_2]$ is a multiple of a conic.}  We may 
assume that the node between $\bar K_1$ and $\bar K_2^1$ is mapped to the 
same point where the node between $\bar K_1$ and $\bar K_2$ is mapped; denote 
this point by $p_2$.  It follows that the image of $\bar K_2$ is uniquely 
determined.  From the irreducibility of the Hurwitz spaces (\cite{FuHu}) it 
follows that we may find an irreducible curve in $K_\beta $ containing $k_1$ 
and a morphism $k_2$ as above.

\noindent
{\bf \boldmath Case 2: $k_*[\bar K_2] = -K_{X_\alpha}$, for $1 \leq \alpha \leq 6$.}  
We may assume that the node between $\bar K_1$ and $\bar K_2^1$ is mapped to 
the same point where the node between $\bar K_1$ and $\bar K_2$ is mapped; 
denote this point by $p_2$.

Since the point $p_2$ does not lie on any $(-1)-$curve, the space of all 
rational divisors in $|-K_{X_\alpha}|$ containing the point $p_2$ is 
isomorphic to the space of all rational divisors in $|-K_{\tilde X_\alpha}|$, 
where $\tilde X_\alpha $ is the blow up of $X_\alpha $ at $p_2$.  It follows 
from the fact that $\tilde X_\alpha$ is a del Pezzo surface of degree at least 
two, that the space of rational curves in $|-K_{\tilde X_\alpha}|$ irreducible 
and thus we conclude also in this case.

\noindent
{\bf \boldmath Case 3: $k_*[\bar K_2] = -K_{X_7}$ and 
$k_*[\bar K_1] = -K_{X_7}$ or $-K_{X_8}$.}  The dual graphs of $k$ and $k_1$ are 
$$ \xygraph {[] !~:{@{.}} 
!{<0pt,0pt>;<20pt,0pt>:} 
{\bullet} [rr] {\bullet} 
*\cir<2pt>{}
!{\save +<0pt,8pt>*\txt{$\scriptstyle \bar K_2$}  \restore}
- [ll]
*\cir<2pt>{}
!{\save +<0pt,8pt>*\txt{$\scriptstyle \bar K_1$}  \restore} } 
\hspace{30pt}
\xygraph {[] !~:{@{.}} 
!{<0pt,0pt>;<20pt,0pt>:} 
{\bullet} [rr] {\bullet} 
*\cir<2pt>{}
!{\save +<0pt,8pt>*\txt{$\scriptstyle \bar K_2 ^1$}  \restore}
- [ll]
*\cir<2pt>{}
!{\save +<0pt,8pt>*\txt{$\scriptstyle \bar K_1$}  \restore} } $$
\vglue0pt {\centerline {Dual graphs of $k$ and $k_1$} \vspace{5pt}}
\noindent
and we have 
$ k_1 \bigl( \bar K_1 \bigr) \cdot k_1 \bigl( \bar K_2^1 \bigr) = 2 $.  
Consider the diagram 
$$ \xymatrix @C=40pt @R=35pt { S_{bir} \ar[r] \ar[dd]^F & S_2 \ar[r] \ar[d] & 
\bar K_1 \times \bar K_1 
\ar[d] ^{( k|_{\bar K_1} , k|_{\bar K_1} )} \\
& \overline \cM _{0,2} \bigl( X, -K_{X_7} \bigr) 
\ar[r]^{\hspace{20pt} \underline {ev}} \ar[d] & X \times X \\
\overline \cM _{bir} \bigl( X, -K_{X_7} \bigr) \ar[r] & 
\overline \cM _{0,0} \bigl( X, -K_{X_7} \bigr) } $$
where all squares are fiber products.  The morphism $F$ is dominant and its 
fiber over a stable map $f$ has length two, unless the image of $f$ contains 
$k(\bar K_1)$.  We denote by $S_{bir}' \subset S_{bir}$ the union of the 
components of $S_{bir}$ dominating 
$\overline \cM _{bir} \bigl( X, -K_{X_7} \bigr)$.  Clearly $k$ and $k_1$ both 
lie in $S_{bir}'$.  Moreover, since the fibers of $F':= F|_{S_{bir}'}$ have 
length two, it follows that $S_{bir}'$ has at most two components.  To prove the 
connectedness of $S_{bir}'$, we assume it is reducible and check that there is a 
point in $K_\beta $ common to the two components of $S_{bir}'$.  This will 
conclude the proof in this case.

Consider the morphism $\varphi : X \rightarrow \bP^2$ determined by $|-K_{X_7}|$.  
We can factor $\varphi $ as the contraction of the divisor class $e_8$ followed 
by the double cover of $\bP^2$ branched along a smooth plane quartic curve $R$.  
The image of $\bar K_2^1$ is a line tangent to $R$.  The image of the 
component $\bar K_1$ is herself a tangent line to the branch curve $R$.  Note 
that in case $k_*[\bar K_1] = -K_{X_7}$ we may assume that this tangent line is 
not a bitangent line nor a flex line.  In case $k_*[\bar K_1] = -K_{X_8}$, it 
follows from the fact that all the rational divisors in $|-K_{X_8}|$ are nodal 
that the image of $\bar K_1$ is not a flex line; the fact that it is not a 
bitangent line follows from the fact that $X$ is a del Pezzo surface.

Let $s \in \bP^2$ be one of the two points such that 
$s \in R \cap \varphi \bigl( k(\bar K_1) \bigr)$, but $s$ is not the point where 
$R$ and $\varphi \bigl( k(\bar K_1) \bigr)$ are tangent.  Through the point $s$ 
there are ten tangent lines to $R$ (counted with multiplicity, and not counting 
the tangent line to $R$ at $s$): tangent lines through $s$ correspond to 
ramification points of the morphism $R \rightarrow \bP^1$ obtained by projecting 
away from the point $s$.  Since $R$ has genus three and the morphism has degree 
three, by the Hurwitz formula we deduce that the degree of the ramification 
divisor is ten, as asserted above.  Let $L \subset \bP^2$ be one of the tangent 
lines to $R$ through $s$ different from $\varphi \bigl( k(\bar K_1) \bigr)$ and 
let $f : \bP^1 \rightarrow X$ be a morphism birational to its image and whose 
image is $\varphi ^{-1} \bigl( L \bigr)$.  The morphism $f$ represents a point 
of $\overline \cM _{bir} \bigl( X, -K_{X_7} \bigr)$ above which the two 
components of $S_{bir}'$ must meet.  Such a point is smooth thanks to Proposition 
\ref{grafico}.  This concludes the proof in this case.

\noindent
{\bf \boldmath Case 4: $k_*[\bar K_2] = -K_{X_7}$ and $k_*[\bar K_2] = -2K_{X_8}$.}  
The dual graphs of $k$ and $k_1$ are 
$$ \xygraph {[] !~:{@{.}} 
!{<0pt,0pt>;<20pt,0pt>:} 
{\bullet} [rr] {\bullet} 
*\cir<2pt>{}
!{\save +<0pt,8pt>*\txt{$\scriptstyle \bar K_2$}  \restore}
- [ll]
*\cir<2pt>{}
!{\save +<0pt,8pt>*\txt{$\scriptstyle \bar K_1$}  \restore} } 
\hspace{30pt}
\xygraph {[] !~:{@{.}} 
!{<0pt,0pt>;<20pt,0pt>:} 
{\bullet} [rr] {\bullet} 
*\cir<2pt>{}
!{\save +<0pt,8pt>*\txt{$\scriptstyle \bar K_2 ^1$}  \restore}
- [ll]
*\cir<2pt>{}
!{\save +<0pt,8pt>*\txt{$\scriptstyle \bar K_1$}  \restore} } $$
\vglue0pt {\centerline {Dual graphs of $k$ and $k_1$} \vspace{5pt}}

We reduce this case to the previous on with the following construction.  
We deform $k$ and $k_1$ inside $K_\beta $ to morphisms $k'$ and $k_1'$ 
respectively with dual graphs 
$$ \xygraph {[] !~:{@{.}} 
!{<0pt,0pt>;<20pt,0pt>:} 
{\bullet} [rr] {\bullet} [rr] {\bullet}
*\cir<2pt>{}
!{\save +<1pt,8pt>*\txt{$\scriptstyle \bar K_2'$}  \restore}
- [ll]
*\cir<2pt>{}
!{\save +<0pt,8pt>*\txt{$\scriptstyle \bar C_2$}  \restore}
- [ll]
*\cir<2pt>{}
!{\save +<0pt,8pt>*\txt{$\scriptstyle \bar C_1$}  \restore} } 
\hspace{30pt}
\xygraph {[] !~:{@{.}} 
!{<0pt,0pt>;<20pt,0pt>:} 
{\bullet} [rr] {\bullet} [rr] {\bullet}
*\cir<2pt>{}
!{\save +<3pt,8pt>*\txt{$\scriptstyle (\bar K_2 ^1)'$}  \restore}
- [ll]
*\cir<2pt>{}
!{\save +<0pt,8pt>*\txt{$\scriptstyle \bar C_2$}  \restore}
- [ll]
*\cir<2pt>{}
!{\save +<0pt,8pt>*\txt{$\scriptstyle \bar C_1$}  \restore} } $$
\vglue0pt {\centerline {Dual graphs of $k'$ and $k_1'$} \vspace{5pt}}
\noindent
where $\bar C_1$ and $\bar C_2$ are mapped to two given distinct rational 
divisors $M_1$ and $M_2$ in $|-K_{X_8}|$.

The strategy is the same for $k$ and for $k_1$, therefore we will 
only describe the deformation for $k$.  We may deform $k|_{\bar K_1}$ 
to the morphism $k'|_{\bar C_1 \cup \bar C_2}$, thanks to the irreducibility 
of $\overline \cM _{bir} \bigl( X, -2K_{X_8} \bigr)$.  This means we 
may deform $k$ to a morphism $\bar k$ which is either $k'$ or it has dual 
graph 
$$ \xygraph {[] !~:{@{.}} 
!{<0pt,0pt>;<20pt,0pt>:} 
{\bullet} [rr] {\bullet} [rr] {\bullet}
*\cir<2pt>{}
!{\save +<0pt,8pt>*\txt{$\scriptstyle \bar C_2$}  \restore}
- [ll]
*\cir<2pt>{}
!{\save +<0pt,8pt>*\txt{$\scriptstyle \bar C_1$}  \restore} 
- [ll]
*\cir<2pt>{}
!{\save +<0pt,8pt>*\txt{$\scriptstyle \bar K_2$}  \restore} } $$
\vglue0pt {\centerline {Dual graph of $\bar k$} \vspace{5pt}}

Since $\bigl( -K_{X_8} \bigr) \cdot \bigl( -K_{X_7} \bigr) = 2$, there 
are at most two irreducible components of the space of morphisms with dual 
graph as above.  Thanks to the previous case, we know that this space is 
connected.

Let $\tilde k' : \bar L \cup \bar C_3 \rightarrow X$ be a stable map 
birational to its image, where $\bar L$ is mapped to the $(-1)-$curve with 
divisor class $e_8$, $\bar C_3$ is mapped to a rational divisor in 
$|-K_{X_8}|$ different from the images of both $\bar C_1$ and $\bar C_2$.  
By the connectedness established above, the (closure of the) same connected 
component of $K_\beta$ containing $\bar k$ contains a morphism 
$\tilde k : \tilde C \rightarrow X$ 
with dual graph 
$$ \xygraph {[] !~:{@{.}} 
!{<0pt,0pt>;<20pt,0pt>:} 
{\bullet} [rr] {\bullet} [rr] {\bullet} [ur] {\bullet} [dd] {\bullet}
*\cir<2pt>{}
!{\save +<8pt,0pt>*\txt{$\scriptstyle \bar C_2$}  \restore}
- [ul] - [ur]
*\cir<2pt>{}
!{\save +<8pt,0pt>*\txt{$\scriptstyle \bar C_1$}  \restore}
 [dl] 
*\cir<2pt>{}
!{\save +<-2pt,8pt>*\txt{$\scriptstyle \bar E$}  \restore}
- [ll]
*\cir<2pt>{}
!{\save +<0pt,8pt>*\txt{$\scriptstyle \bar C_3$}  \restore} 
- [ll] 
*\cir<2pt>{}
!{\save +<0pt,8pt>*\txt{$\scriptstyle \bar L$}  \restore} } $$
\vglue0pt {\centerline {Dual graph of $\tilde k$} \vspace{5pt}}
\noindent
where $\bar E$ is contracted to the base-point of $|-K_{X_8}|$.  To check 
that $\tilde k$ is in $K_\beta$ it is enough to check that $\tilde k$ 
represents a smooth point of 
$\overline \cM _{0,0} \bigl( X, k_*[\bar K] \bigr)$.

The point represented by $\tilde k : \tilde C \rightarrow X$ in 
$\overline \cM _{0,0} \bigl( X, k_*[\bar K] \bigr)$ is smooth if 
${\rm H}^0 \bigl( \tilde C , \cC _{\tilde k} \bigr) = 0$ (we are 
using the notation of (\ref{conor})).  We have a natural inclusion 
$$ {\rm H}^0 \bigl( \tilde C , \cC _{\tilde k} \bigr) \subset 
{\rm H}^0 \bigl( \tilde C , \tilde k ^* \Omega ^1 _X \otimes \omega _{\tilde C} \bigr) $$

We prove first that any global section of $\cC _{\tilde k}$ is zero on 
$\bar L \cup \bar C_3$ and then that any global section of 
$\tilde k ^* \Omega ^1 _X \otimes \omega _{\tilde C}$ vanishing on 
$\bar L \cup \bar C_3$ is the zero section.

The first assertion is clear from Proposition \ref{graficone}: there are 
no non-zero global sections of $\cC _{\tilde k}$ on $\bar L$, and since the 
sheaf $\cC _{\tilde k}$ is locally free near the node between $\bar L$ and 
$\bar C_3$, it follows that a global section of $\cC _{\tilde k}$ must vanish 
at the node.  Since the degree of the sheaf $\cC _{\tilde k}$ on $\bar C_3$ is 
zero, it follows that a global section of $\cC _{\tilde k}$ must vanish on 
$\bar L \cup \bar C_3$.

The second assertion is a consequence of the fact that $h^0 \bigl( \tilde C , 
\tilde k ^* \Omega ^1 _X \otimes \omega _{\tilde C} \bigr) = 1$, and that a 
non-zero section of the sheaf $\tilde k ^* \Omega ^1 _X \otimes \omega _{\tilde C}$ 
is not identically zero on $\bar C_3$.

To compute $h^0 \bigl( \tilde C , 
\tilde k ^* \Omega ^1 _X \otimes \omega _{\tilde C} \bigr)$, we use Serre 
duality to deduce that 
$$ h^0 \bigl( \tilde C , \tilde k ^* \Omega ^1 _X \otimes \omega _{\tilde C} \bigr) 
= h^1 \bigl( \tilde C , \tilde k ^* \cT _X \bigr) $$
There is a short exact sequence of sheaves 

\vspace{-10pt} \centerline{
\( \xymatrix { 0 \ar[r] & \tilde k ^* \cT _X \ar[r] & k_{\bar L} ^* \cT_X \oplus 
k_{\bar E} ^* \cT_X \oplus \bigoplus _i k_i ^* \cT_X \ar[r] & \!\!\!\!
\mathop {\bigoplus } \limits _{
\raisebox{5pt} {\begin{tabular} {c} 
\scriptsize $\nu $ a node \\[-5pt]
\scriptsize of $\tilde C$
\end{tabular} } }
^{\vphantom{{
\raisebox{-5pt} {\begin{tabular} {c} 
\scriptsize $\nu $ a node \\[-5pt]
\scriptsize of $\tilde C$
\end{tabular} }}}} \!\!\!\!
\cT _{X, \tilde k (\nu)} \ar[r] & 0 } \) }
\noindent
Note that the restriction of $\tilde k ^* \cT_X$ to every non contracted 
component is isomorphic to $\cO_{\bP^1} (2) \oplus \cO_{\bP^1} (-1)$, where 
the subsheaf of degree two is canonically the tangent sheaf of the component.  
The associated long exact sequence to the sequence above is 
$$ 0 \longrightarrow k^{h^0} \longrightarrow k^{14} 
\longrightarrow k^8 \longrightarrow k^{h^1} \longrightarrow 0 $$

Since the images of all the non-contracted components are pairwise 
transverse (all the intersection numbers are one), and since the only global 
sections come from the tangent vector fields, it follows that any global 
section must vanish at all nodes.  Thus, there are two global sections 
coming from the each of the curves $\bar L$, $\bar C_1$ and $\bar C_2$ and 
only one coming from $\bar C_3$.  We deduce that $h^0 = 7$ and finally 
$h^1 = 1$, as asserted above.

Let us go back to the sheaf 
$\tilde k ^* \Omega ^1_X \otimes \omega _{\tilde C}$.  We just computed that 
this sheaf has exactly one global section.  We have the following 
decomposition for the degrees of the restrictions of the sheaf 
$\tilde k ^* \Omega ^1_X \otimes \omega _{\tilde C}$ to each component: 
$$ \xygraph {[] !~:{@{.}} 
!{<0pt,0pt>;<20pt,0pt>:} 
{\bullet} [rr] {\bullet} [rr] {\bullet} [ur] {\bullet} [dd] {\bullet}
*\cir<2pt>{}
!{\save +<0pt,-8pt>*\txt{$\scriptstyle (-3,0)$}  \restore}
!{\save +<8pt,2pt>*\txt{$\scriptstyle \bar C_2$}  \restore}
- [ul] - [ur]
*\cir<2pt>{}
!{\save +<0pt,8pt>*\txt{$\scriptstyle (-3,0)$}  \restore}
!{\save +<8pt,-2pt>*\txt{$\scriptstyle \bar C_1$}  \restore}
 [dl] 
*\cir<2pt>{}
!{\save +<-4pt,8pt>*\txt{$\scriptstyle (1,1)$}  \restore}
!{\save +<-2pt,-8pt>*\txt{$\scriptstyle \bar E$}  \restore}
- [ll]
*\cir<2pt>{}
!{\save +<0pt,8pt>*\txt{$\scriptstyle (-2,1)$}  \restore} 
!{\save +<0pt,-8pt>*\txt{$\scriptstyle \bar C_3$}  \restore} 
- [ll] 
*\cir<2pt>{}
!{\save +<0pt,8pt>*\txt{$\scriptstyle (-3,0)$}  \restore} 
!{\save +<0pt,-8pt>*\txt{$\scriptstyle \bar L$}  \restore} } $$
\vglue0pt {\centerline {Degrees of the sheaf 
$\tilde k ^* \Omega ^1_X \otimes \omega _{\tilde C}$} \vspace{5pt}}
\noindent
where the pair of numbers next to a vertex represent the degrees of 
$\tilde k ^* \Omega ^1_X \otimes \omega _{\tilde C}$ restricted to the 
component represented by the corresponding vertex.  We examine the 
vertex of valence three in the dual graph.  Necessary conditions for a 
section of $\tilde k ^* \Omega ^1_X \otimes \omega _{\tilde C}$ on 
$\bar E$ to extend to a global section are that the section ``points 
in the right direction'' at the nodes.  These are clearly linear 
conditions and there are three such conditions.  Moreover, every 
section satisfying the stated conditions extends uniquely to a global 
section: this is obvious on the components $\bar C_1$ and $\bar C_2$.  
For the remaining components, note that every global section must 
vanish at the node between $\bar L$ and $\bar C_3$, since the 
intersection number $\tilde k (\bar L) \cdot \tilde k (\bar C_3)$ 
equals one, and therefore the intersection is transverse.  Thus 
every global section of 
$\tilde k ^* \Omega ^1_X \otimes \omega _{\tilde C}$ is uniquely 
determined by its restriction to $\bar E$.  Thus the only way a 
section can be identically zero on $\bar C_3$ is if the sections 
on $\bar E$ all vanish at the node $\bar C_3 \cap \bar E$.  Note that 
the three tangent directions of the images of $\bar C_1$, $\bar C_2$ 
and $\bar C_3$ at their common point $p$ are pairwise independent.  
Choose homogeneous coordinates $E_1, E_2$ on $\bar E$ such that 
$[0,1] = \bar E \cap C_1$, $[1,0] = \bar E \cap C_2$.  Choose local 
coordinates $u,v$ on $X$ near $p$ such that the zero set of $u$ is 
tangent to the image of $\bar C_1$ and the zero set of $v$ is tangent 
to the image of $\bar C_2$.  Rescaling by a non-zero constant $u$ and 
$v$ we may also assume that the zero set of $u+v$ is tangent to the 
image of $\bar C_3$.  The restrictions of the global sections of 
$\tilde k ^* \Omega ^1_X \otimes \omega _{\tilde C}$ to $\bar E$ are 
multiples of the section 
$$ \sigma := E_0 du + E_1 dv $$

In particular, if a section vanishes at one of the nodes between 
$\bar E$ and $\bar C_i$, then it vanishes identically.  This concludes 
the proof that the sheaf $\cC _{\tilde k}$ has no global sections and 
thus we conclude that $\tilde k$ is a smooth point of 
$\overline \cM _{0,0} \bigl( X, k_*[\bar K] \bigr)$.

We now resume our argument.  It is clear that $\tilde k$ is also a limit 
of morphisms $\tilde k'$ with dual graph 
$$ \xygraph {[] !~:{@{.}} 
!{<0pt,0pt>;<20pt,0pt>:} 
{\bullet} [rr] {\bullet} [rr] {\bullet}
*\cir<2pt>{}
!{\save +<1pt,8pt>*\txt{$\scriptstyle \tilde K_2'$}  \restore}
- [ll]
*\cir<2pt>{}
!{\save +<0pt,8pt>*\txt{$\scriptstyle \bar C_2$}  \restore}
- [ll]
*\cir<2pt>{}
!{\save +<0pt,8pt>*\txt{$\scriptstyle \bar C_1$}  \restore} } $$
\vglue0pt {\centerline {Dual graph of the morphisms $\tilde k'$} \vspace{5pt}}
\noindent
which is precisely what we wanted to prove.  This completes the reduction 
of this case to {\bf Case 3}, and thus this case is proved.

\noindent
{\bf \boldmath Case 5: $k_*[\bar K_2] = -K_{X_7}$ and $k_*[\bar K_1] = -3K_{X_8}$.}  
We also reduce this case to {\bf Case 3}.  As before, thanks to the 
irreducibility of $\overline \cM _{bir} \bigl( X, -3K_{X_8} \bigr)$ we may deform 
the morphism $k$ so that $k|_{\bar K_1}$ is replaced by the birational morphism 
$k' : \bar C_1 \cup \bar C_2 \rightarrow X$, where $\bar C_1$ is immersed and 
represents $|-2K_{X_8}|$, and $\bar C_2$ is mapped to a given rational divisor in 
$|-K_{X_8}|$.  After possibly sliding the component $\bar K_2$ along $\bar C_1$, 
we may suppose that the dual graph of $k'$ is the following:
$$ \xygraph {[] !~:{@{.}} 
!{<0pt,0pt>;<20pt,0pt>:} 
{\bullet} [rr] {\bullet} [rr] {\bullet}
*\cir<2pt>{}
!{\save +<1pt,8pt>*\txt{$\scriptstyle \bar K_2'$}  \restore}
- [ll]
*\cir<2pt>{}
!{\save +<0pt,8pt>*\txt{$\scriptstyle \bar C_2$}  \restore}
- [ll]
*\cir<2pt>{}
!{\save +<0pt,8pt>*\txt{$\scriptstyle \bar C_1$}  \restore} } $$
\vglue0pt {\centerline {Dual graph of $k'$} \vspace{5pt}}
\noindent
Similar remarks apply to $k_2$.  This completes the reduction to {\bf Case 3} 
and the proof in this case.

\noindent
{\bf \boldmath Case 6: $k_*[\bar K_2] = -2K_{X_8}$ and $k_*[\bar K_1] = -2K_{X_8}$.}  
Since the intersection product $\bigl( -2K_{X_8} \bigr) ^2$ equals four, and the 
space $\overline \cM _{bir} \bigl( X, -2K_{X_8} \bigr)$ is irreducible, it follows 
that there are at most four irreducible components of morphisms in standard form 
representing the divisor class $-4K_{X_8}$.  Let 
$c : \bar C_1 \cup \bar C_2 \longrightarrow X$ be a stable map birational to its 
image such that $\bar C_i$ is mapped to a $(-1)-$curve $C_i$ and 
$C_1 + C_2 = -2K_{X_8}$.   Consider the morphism 
$$ \xymatrix @C=40pt { {\rm Sl} _k (\bar K_1) \ar[r]^{\pi \hspace{20pt}} & 
\overline \cM _{bir} \bigl( X, -2K_{X_8} \bigr) } $$

The morphism $\pi $ is dominant.  Thus we may find a morphism 
$k' : \bar C_1 \cup \bar C_2 \cup \bar K_2 \longrightarrow X$ such that 
$k' |_{\bar C_1 \cup \bar C_2} \simeq c$, lying in the same irreducible component 
of $K_\beta $ as $k$.  We have two possibilities for the dual graph of $k'$: 
$$ \xygraph {[] !~:{@{.}} 
!{<0pt,0pt>;<20pt,0pt>:} 
{\bullet} [rr] {\bullet} [rr] {\bullet}
*\cir<2pt>{}
!{\save +<1pt,8pt>*\txt{$\scriptstyle \bar K_2$}  \restore}
- [ll]
*\cir<2pt>{}
!{\save +<0pt,8pt>*\txt{$\scriptstyle \bar C_2$}  \restore}
- [ll]
*\cir<2pt>{}
!{\save +<0pt,8pt>*\txt{$\scriptstyle \bar C_1$}  \restore} } 
\hspace{30pt}
\xygraph {[] !~:{@{.}} 
!{<0pt,0pt>;<20pt,0pt>:} 
{\bullet} [rr] {\bullet} [rr] {\bullet}
*\cir<2pt>{}
!{\save +<3pt,8pt>*\txt{$\scriptstyle \bar C_2$}  \restore}
- [ll]
*\cir<2pt>{}
!{\save +<0pt,8pt>*\txt{$\scriptstyle \bar C_1$}  \restore}
- [ll]
*\cir<2pt>{}
!{\save +<0pt,8pt>*\txt{$\scriptstyle \bar K_2$}  \restore} } $$
\vglue0pt {\centerline {Possible dual graphs of $k'$} \vspace{5pt}}

We want to reduce to the case in which $\bar K_2$ is adjacent to 
$\bar C_2$.  Consider the morphism 
$$ \xymatrix @C=40pt { 
{\rm Sl} _{k'} (\bar K_2) \ar[r]^{\hspace{10pt}a} & \bar C_1 } $$
and as usual this morphism is dominant.  This means that we may slide 
the node between $\bar K_2$ and $\bar C_1$ until it reaches the node 
between $\bar C_1$ and $\bar C_2$.  The resulting morphism $\tilde k$ 
has dual graph 
$$ \xygraph {[] !~:{@{.}} 
!{<0pt,0pt>;<20pt,0pt>:} 
{\bullet} [rr] {\bullet} [ur] {\bullet} [dd] {\bullet} 
*\cir<2pt>{}
!{\save +<8pt,0pt>*\txt{$\scriptstyle \bar C_2$}  \restore}
- [ul]
*\cir<2pt>{}
!{\save +<-3pt,8pt>*\txt{$\scriptstyle \bar E$}  \restore}
- [ur]
*\cir<2pt>{}
!{\save +<8pt,0pt>*\txt{$\scriptstyle \bar C_1$}  \restore}
 [dl] - [ll]
*\cir<2pt>{}
!{\save +<0pt,8pt>*\txt{$\scriptstyle \tilde K_2$}  \restore} } $$
\vglue0pt {\centerline {Dual graph of $\tilde k$} \vspace{5pt}}
\noindent
where $\bar E$ is contracted by $\tilde k$.  It is easy to check that 
this morphism represents a smooth point of 
$\overline \cM _{0,0} \bigl( X, -4K_{X_8} \bigr)$ and that it is also 
a limit of morphisms $\tilde k'$ with dual graph 
$$ \xygraph {[] !~:{@{.}} 
!{<0pt,0pt>;<20pt,0pt>:} 
{\bullet} [rr] {\bullet} [rr] {\bullet}
*\cir<2pt>{}
!{\save +<1pt,8pt>*\txt{$\scriptstyle \bar K_2$}  \restore}
- [ll]
*\cir<2pt>{}
!{\save +<0pt,8pt>*\txt{$\scriptstyle \bar C_2$}  \restore}
- [ll]
*\cir<2pt>{}
!{\save +<0pt,8pt>*\txt{$\scriptstyle \bar C_1$}  \restore} } $$
\vglue0pt {\centerline {Dual graph of $\tilde k'$} \vspace{5pt}}

Thus we may indeed assume that $\bar K_2$ is adjacent to $\bar C_2$.  
Note that since $-2K_{X_8} \cdot C_1 = 2$, it follows that there are at 
most two connected components in the space of morphisms in standard form 
representing the divisor class $-4K_{X_8}$.  To conclude, it is enough to 
show that we may ``exchange'' the two intersection points 
$C_2 \cap k'(\bar K_2)$ by a connected path contained in $K_\beta $.

Consider the morphism $\varphi : X \rightarrow \bP^3$ induced by the 
linear system $|-2K_{X_8}|$.  We have already seen that the image is a 
quadric cone $Q$ and that the morphism is ramified along a smooth curve 
$R$ which is the complete intersection of $Q$ with a cubic surface.  The 
$(-1)-$curves $C_1$ and $C_2$ have as image the intersection of $Q$ with 
a plane which is everywhere tangent to the curve $R$ (and does not contain 
the vertex of the cone).  Let $p$ be one of the intersection points of 
$\varphi (C_2)$ with $R$.  Projection away from the tangent line $L$ to $R$ 
at $p$ determines a morphism $\pi _L : R \rightarrow \bP^1$ of degree four.  
Since the genus of $R$ is four, it follows from the Hurwitz formula that 
the degree of the ramification divisor of $\pi _L$ is 14.  It is immediate 
to check that $\pi _L$ ramifies above the tangent plane to $Q$ at $p$, and 
that the ramification index is two.  It is also immediate that above the 
plane containing $\varphi (C_2)$ the ramification index is two.  We deduce 
that there are planes in the pencil containing $L$ which are tangent to $R$ 
and are not the tangent plane to $Q$ at $p$ nor the plane containing 
$\varphi (C_2)$.  Such planes correspond to rational divisors $H$ in 
$|-2K_{X_8}|$ with the property that $H \cap C_2$ consists of the unique 
point $\varphi ^{-1} (p)$.  Let $\nu : \bP^1 \rightarrow X$ be a birational 
morphism whose image is one of the divisors $H$ constructed above.  The 
morphism $\nu$ represents a morphism in 
$\overline \cM _{bir} \bigl( X, -2K_{X_8} \bigr)$, and since this space is 
irreducible, we may deform $\tilde k'$ to a morphism 
$\bar k : \bar C_1 \cup \bar C_2 \cup \bar H \longrightarrow X$ with dual 
graph 
$$ \xygraph {[] !~:{@{.}} 
!{<0pt,0pt>;<20pt,0pt>:} 
{\bullet} [rr] {\bullet} [rr] {\bullet}
*\cir<2pt>{}
!{\save +<1pt,8pt>*\txt{$\scriptstyle \bar H$}  \restore}
- [ll]
*\cir<2pt>{}
!{\save +<0pt,8pt>*\txt{$\scriptstyle \bar C_2$}  \restore}
- [ll]
*\cir<2pt>{}
!{\save +<0pt,8pt>*\txt{$\scriptstyle \bar C_1$}  \restore} } $$
\vglue0pt {\centerline {Dual graph of $\bar k$} \vspace{5pt}}
\noindent
and such that $\bar k|_{\bar H} \simeq \nu$.  The morphism $\bar k$ 
represents a smooth point of the space 
$\overline \cM _{0,0} \bigl( X, -4K_{X_8} \bigr)$, thanks to 
Proposition \ref{grafico}.  Thus $\bar k \in K_\beta $ and it lies 
in the same connected component of $K_\beta $ as $k$.

Applying the same construction to the morphism $k_2$, we obtain that 
also $k_2$ lies in the same connected component of $K_\beta $ as 
$\bar k$.  This concludes the proof of this case.

\noindent
{\bf \boldmath Case 7: $k_*[\bar K_2] = -2K_{X_8}$ and $k_*[\bar K_1] = -3K_{X_8}$.}  
Let $c : \bar C_1 \cup \bar C_2 \longrightarrow X$ be a morphism 
birational to its image, such that $c (\bar C_2)$ is a rational 
divisor in $|-K_{X_8}|$ and $c (\bar C_1)$ is a general rational 
divisor in $|-2K_{X_8}|$.  Consider the morphism 
$$ \xymatrix @C=40pt { {\rm Sl} _k (\bar K_1) \ar[r]^{\pi \hspace{20pt}} & 
\overline \cM _{bir} \bigl( X, -3K_{X_8} \bigr) } $$
and note that as usual it is dominant.  Therefore we may deform $k$ to 
a morphism $k' : \bar C_1 \cup \bar C_2 \cup \bar K_2 \longrightarrow X$ 
such that $k' |_{\bar C_1 \cup \bar C_2} \simeq c$.

As before, we may slide the component $\bar K_2$ along $\bar C_1$ until it 
reaches $\bar C_2$, and reduce to the case in which $\bar K_2$ is adjacent to 
$\bar C_2$.  The same considerations of the final step of the previous case 
allow us to conclude.

This concludes the proof of the connectedness of $K_\beta $.

{\bf Step 7.}  We now simply collect all the information we obtained, to 
conclude the proof of the theorem.  {\bf Step 4} and {\bf Step 5} imply 
(under the hypotheses of the theorem) that every irreducible component of 
$\overline \cM _{bir} \bigl( X, \beta \bigr)$ is either empty or it 
contains a point lying in $K _\beta$.  {\bf Step 6} then implies that 
there is at most one component 
of $\overline \cM _{bir} \bigl( X, \beta \bigr)$ containing $K_\beta $.  
Thus if $\overline \cM _{bir} \bigl( X, \beta \bigr)$ is not empty, then 
it consists of exactly one irreducible component.  This concludes the 
proof of the theorem.  \bo

\begin{prop}
Let $X$ be a del Pezzo surface.  If $\beta $ is a nef divisor which is not 
a multiple of a conic, then the space 
$\overline \cM _{bir} \bigl( X, \beta \bigr)$ is not empty.
\end{prop}
{\it Proof.}  We may write 
$$ \beta = n_8 (-K_{X_8}) + \ldots + n_2 (-K_{X_2}) + \beta ' $$
where $n_8, \ldots , n_2 \geq 0$ and $\beta '$ is a nef divisor on a 
del Pezzo surface $X_1$ of degree eight dominated by $X$ (if 
$X \simeq \bP^2$ the assertion is obvious).

If $n_8 \geq 2$, then we define $n$ and $r$ by the conditions 
$n_8 = 2(n-1) + r$, where $r = 2,3$; if $n_8 = 1$, then we define 
$n = r = 1$; if $n_8 = 0$, then we define $n = r = 0$.  Let 
$$ \begin{array} {rcl@{~~{\text{ where }}}l} 
{\rm Pic} \bigl( \bP^1 \times \bP^1 \bigr) & \simeq & \bZ \ell _1 \oplus \bZ \ell _2 
& \ell _1 = \{ p \} \times \bP^1 ~,~ \ell _2 = \bP^1 \times \{ p \} \\[10pt]
{\rm Pic} \bigl( Bl_p (\bP^2) \bigr) & \simeq & \bZ \ell \oplus \bZ e & 
\ell ^2 = 1 ~,~ \ell \cdot e = 0 ~,~ e^2 = -1
\end{array} $$
and write 
$$ \begin{array} {rcl@{~~{\text{ if }}}l}
\beta ' & = & n_1 \bigl( \ell _1 + \ell _2 \bigr) + n_0 \ell _2 
& X \simeq \bP^1 \times \bP^1 \\[7pt]
\beta ' & = & n_1 \ell + n_0 \bigl( \ell - e \bigr) & X \simeq Bl_p (\bP^2)
\end{array} $$
where $n_1 \geq 0$, $n_0 \geq 0$ (we may need to exchange $\ell _1$, $\ell _2$).  
Note that with this notation the divisor $\beta $ is multiple of a conic if 
and only if $n_8 = n_7 = \ldots = n_1 = 0$.

Choose 
\begin{itemize}
\item $n-1$ distinct rational integral nodal divisors $C_2^8, \ldots , C_n^8$ 
in $|-2K_{X_8}|$;
\item a rational integral nodal divisor $C_1^8$ (different from the previous ones 
if $r=2$) in $|-rK_{X_8}|$;
\item $n_i$ distinct rational integral nodal divisors $C_1^i, \ldots , C_{n_i}^i$ 
in $|-K_{X_i}|$;
\item $n_1$ distinct integral divisors $C_1 , \ldots , C_{n_1}$ lying 
in $|\ell _1 + \ell _2|$, if $X_1 \simeq \bP^1 \times \bP^1$ and lying in 
$|\ell |$, if $X_1 \simeq Bl_p (\bP^2)$;
\item an integral divisor $C'$ in $|\ell _2|$ or $|\ell - e|$.
\end{itemize}

Having made these choices, we may now consider the stable map of genus 
zero $f : \bar C \rightarrow X$, with dual graph 
$$ \xygraph {[] !~:{@{.}} 
!{<0pt,0pt>;<20pt,0pt>:} 
{\bullet} [rr] {\bullet} [rr] 
{\bullet} [rr] {\bullet} [rr] 
{\bullet} [rr] 
{\bullet} [rr] {\bullet} 
*\cir<2pt>{}
!{\save +<0pt,8pt>*\txt{$\scriptstyle \bar C'$}  \restore}
- [ll]
*\cir<2pt>{}
!{\save +<0pt,8pt>*\txt{$\scriptstyle \bar C_{n_1}^1$}  \restore}
- [l] {\cdots} - [l]
*\cir<2pt>{}
!{\save +<0pt,8pt>*\txt{$\scriptstyle \bar C_{n_7}^7$}  \restore}
- [l] {\cdots} - [l]
*\cir<2pt>{}
!{\save +<0pt,8pt>*\txt{$\scriptstyle \bar C_1^7$}  \restore}
- [ll]
*\cir<2pt>{}
!{\save +<0pt,8pt>*\txt{$\scriptstyle \bar C_n^8$}  \restore}
- [l] {\cdots} - [l]
*\cir<2pt>{}
!{\save +<0pt,8pt>*\txt{$\scriptstyle \bar C_2^8$}  \restore}
- [ll]
*\cir<2pt>{}
!{\save +<0pt,8pt>*\txt{$\scriptstyle \bar C_1^8$}  \restore} } $$
\vglue0pt {\centerline {Dual graph of $\bar f$} \vspace{5pt}}
\noindent
where of course we ignore a component if the corresponding curve 
without a bar has not been defined.  The morphism $f$ on a component 
$\bar D$ is the normalization of the curve $D$ followed by inclusion 
in $X$, if $D \neq C'$, and it is a multiple cover of degree $n_0$, 
if $D = C'$.

All the restrictions of $f$ to the irreducible components of $\bar C$ 
different from $\bar C_1^8$ are free morphisms; the cohomology group 
${\rm H}^1 \bigl( \bar C_1^8 , f^* \cT_X \bigr)$ is immediately seen 
to be zero.  Thus we may deform $f$ to a morphism lying in 
$\cM _{0,0} \bigl( X, \beta \bigr)$.  If the general deformation of 
$f$ were a morphism not birational to its image, then $f_*[\bar C]$ 
would not be reduced.  Since this is not the case, it follows that we 
may deform $f$ to a morphism with irreducible domain, which is 
birational to its image.  This proves that 
$\cM _{bir} \bigl( X, \beta \bigr) \neq \emptyset $, if $\beta $ is not a 
multiple of a conic.  This concludes the proof of the proposition.  \bo

\noindent
{\it Remark 1}.  The spaces $\cM _{0,0} \bigl( X, m C \bigr)$, where $C$ 
is the class of a conic, are easily seen to be irreducible, for $m \geq 1$.  
If $m = 1$, we have $\overline \cM _{bir} \bigl( X, C \bigr) \simeq \bP^1$.  
If $m \geq 1$, then there is a morphism 
$\cM _{0,0} \bigl( X, m C \bigr) \longrightarrow \cM _{0,0} \bigl( X, C \bigr)$, 
obtained by ``forgetting the multiple cover.''  The fibers of this morphism 
are birational to Hurwitz schemes, which are irreducible (\cite{FuHu}).  
The irreducibility of $\cM _{0,0} \bigl( X, m C \bigr)$ follows.

\noindent
{\it Remark 2}.  If $L$ is an integral divisor of anticanonical degree one, 
then either $L$ is a $(-1)-$curve, or it is the anticanonical divisor on a 
del Pezzo surface of degree one.

If $L$ is a $(-1)-$curve, the space $\overline \cM _{0,0} \bigl( X, L \bigr)$ 
has dimension zero and length one; it therefore consists of a single reduced 
point and is irreducible.

If $L = -K_X$, the three spaces $\overline \cM _{0,0} \bigl( X, -K_X \bigr)$, 
$\cM _{0,0} \bigl( X, -K_X \bigr)$ and $\overline \cM _{bir} \bigl( X, -K_X \bigr)$ 
are all equal and have dimension zero and length twelve.  They are 
not irreducible.  For a general del Pezzo surface of degree 
one, the space $\overline \cM _{bir} \bigl( X, -K_X \bigr)$ is reduced 
and consists of exactly twelve points.  This happens precisely when the rational 
divisors in $|-K_X|$ are all nodal.

\pagestyle{headings}
\section{Divisors of Small Degree on $X_8$}

\subsection{The Divisor $-K_{X_8} - K_{X_7}$}

Here we prove the irreducibility of the spaces 
$\overline { \cM } _{bir} \bigl( X , \beta ) \bigr)$, where the degree 
of $X$ is one, $\beta $ is ample and the anticanonical degree of 
$\beta $ is three.  We already saw (Theorem \ref{maschera}) that the 
space $\overline { \cM } _{bir} \bigl( X , -3K_X \bigr)$ is 
irreducible.  The following proofs are similar to the proof of 
Theorem \ref{maschera}.


\begin{lem} \label{coccode}
Let $X$ be a del Pezzo surface of degree one.  Suppose that all the rational 
divisors in $|-K_X|$ are nodal and that 
$\overline { \cM } _{bir} \bigl( X,-2K_X \bigr)$ is irreducible.  Let 
$L \subset X$ be a $(-1)-$curve and let $b : X \rightarrow X'$ be the 
contraction of $L$.  Then the space 
$\overline { \cM } _{bir} \bigl( X,-K_X -K_{X'} \bigr)$ is irreducible.
\end{lem}
{\it Proof.}  Let $f : \bP^1 \rightarrow X$ be a morphism in 
$\overline { \cM } _{bir} \bigl( X,-K_X -K_{X'} \bigr)$ and suppose 
that the image of $f$ contains the independent point $p$.  Consider the space 
of morphisms of $\overline { \cM } _{bir} \bigl( X,-K_X -K_{X'} \bigr)$ 
in the same irreducible component as $[f]$ which contain the point $p$ in 
their image, denote this space by $\overline \cM _{bir} (p)$.  It follows 
immediately from the dimension estimates (\ref{dimdibarbi}) that 
$\dim _{[f]} \overline \cM _{bir} (p) = 1$ and that $[f]$ is a smooth 
point of $\overline \cM _{bir} (p)$.  We may therefore find a smooth 
irreducible projective curve $B$, a normal surface 
$\pi : S \rightarrow B$ and a morphism $F : S \rightarrow X$ such that 
the induced morphism $B \rightarrow \overline \cM _{bir} (p)$ is 
surjective onto the component containing $[f]$.  From \cite{Ko} 
Corollary II.3.5.4, it follows immediately that the morphism $F$ is 
dominant.  We want to show that there are reducible fibers of $\pi$.  
The argument is the same that appears at the beginning of the proof of 
Theorem \ref{maschera}.

Thus there must be a morphism $f_0 : \bar C \rightarrow X$ with reducible 
domain in the family of stable maps parametrized by 
$B$, and since all such morphisms contain the general point $p$ in their 
image, the same is true of the morphism $f_0$.  In particular, since the 
point $p$ does not lie on any rational curve of anticanonical degree one, 
it follows that $\bar C$ consists of exactly two components $\bar C_1$ 
and $\bar C_2$, where each $\bar C_i$ is irreducible and we may assume 
that $f_0 (\bar C_1)$ has anticanonical degree one and $f_0 (\bar C_2)$ 
has anticanonical degree two.  Denote by $C_i$ the image of $\bar C_i$.  
It also follows from the definition of an independent point and 
Proposition \ref{grafico} that $f_0$ represents a smooth point of 
$\overline { \cM } _{bir} \bigl( X , -K_X -K_{X'} \bigr)$.

There are two possibilities for $C_1$: either it is a $(-1)-$curve or it 
is a rational divisor in the anticanonical system.  We want to prove that 
we may assume that $C_1$ is not a $(-1)-$curve.

Suppose $C_1$ is a $(-1)-$curve.  The morphism $f_0 |_{\bar C_2}$ is 
a free morphism, because the image contains a general point and has 
anticanonical degree two.  Moreover the image $C_2$, being a curve of 
anticanonical degree two, is one of the following: a conic, the 
anticanonical divisor on a del Pezzo surface of degree two dominated by 
$X$ or a divisor in $|-2K_X|$.  In all these cases we know that the 
space $\overline { \cM } _{bir} \bigl( X , (f_0)_* [\bar C_2] \bigr)$ 
is irreducible.  Thus we may deform $f_0|_{\bar C_2}$ to a curve with two 
irreducible components, both mapped to $(-1)-$curves.  Considering the 
space ${\rm Sl}_{f_0} (\bar C_2)$ we conclude that we may deform $f_0$ to 
a morphism $f_1 : \bar L_1 \cup \bar L_2 \cup \bar L_3 \rightarrow X$ 
where each component $\bar L_i$ is mapped to a different $(-1)-$curve 
$L_i$ on $X$.

We deduce that we have $L_1 + L_2 + L_3 = -2K_X + L$ and $L_1$, $L_2$ 
and $L_3$ are distinct $(-1)-$curves.  Thanks to Lemma \ref{tremendi} we 
conclude that there is a standard basis $\{ \ell , e_1 , \ldots , e_8 \}$ 
of ${\rm Pic} (X)$ such that 
\begin{eqnarray} \label{noduno}
\left\{ \begin{array} {rcl}
L_1 & = & -2K_X - e_1 \\[5pt]
L_2 & = & e _8 \\[5pt]
L_3 & = & e_1 
\end{array} 
\right. &
{\text{ ~,~ }} &
\left\{ \begin{array} {rcl}
L_1 & = & -2K_X - (\ell - e_7 - e_8 ) \\[5pt]
L_2 & = & \ell - e_7 - e _8 \\[5pt]
L_3 & = & e_8 
\end{array} 
\right. 
\end{eqnarray}
or
\begin{eqnarray} \label{nodini}
&\left\{ \begin{array} {rcl}
L_1 & = & -2K_X - (\ell - e_7 - e_8) \\[5pt]
L_2 & = & e_1 \\[5pt]
L_3 & = & \ell - e_1 - e_7 
\end{array} \right. 
\end{eqnarray}
after possibly permuting the indices 1, 2 and 3.

The next step in the deformation is to produce a component mapped to 
the divisor class $-2K_X$.

In the first case of (\ref{noduno}), the component $\bar L_1$ is 
adjacent to both $\bar L_2$ and $\bar L_3$, since $L_2 \cdot L_3 = 0$.  
We may therefore consider ${\rm Sl}_{f_1} (\bar L_1 \cup \bar L_3)$ to 
smooth $\bar L_1 \cup \bar L_3$ to a single component $\bar K$ mapped 
to $-2K_X$.

In the second case of (\ref{noduno}), either $\bar L_1$ and $\bar L_2$ 
are adjacent and it is enough to smooth their union to conclude, or 
$\bar L_2$ is adjacent to $\bar L_3$ and not to $\bar L_1$.  If this 
happens, then we may smooth the union $\bar L_2 \cup \bar L_3$ to a 
single irreducible component $\bar Q$, mapped to the conic $\ell - e_7$.  
Denote the resulting morphism by $f_1'$.  We may consider the dominant 
morphism 
$$ a: {\rm Sl} _{f_1'} (\bar Q) \longrightarrow \bar L_1 $$
and let $\bar e \in \bar L_1$ be a point mapped to a point lying on 
the $(-1)-$curve with divisor class $L_2 = \ell - e_7 - e_8$.  Since 
$a$ is dominant, we may find a morphism $f_1''$ such that 
$a(f_1'') = \bar e$.  By construction, the dual graph of the morphism 
$f_1''$ is 
$$ \xygraph {[] !~:{@{.}} 
!{<0pt,0pt>;<20pt,0pt>:} 
{\bullet} [rr] {\bullet} [rr] {\bullet} 
*\cir<2pt>{}
!{\save +<0pt,8pt>*\txt{$\scriptstyle \bar L_3$}  \restore}
- [ll]
*\cir<2pt>{}
!{\save +<0pt,8pt>*\txt{$\scriptstyle \bar L_2$}  \restore}
- [ll]
*\cir<2pt>{}
!{\save +<0pt,8pt>*\txt{$\scriptstyle \bar L_1$}  \restore} } $$
\vglue0pt {\centerline {Dual graph of $f_1''$} \vspace{5pt}}
\noindent
and we may now smooth $\bar L_1 \cup \bar L_2$ to conclude.

In the case of (\ref{nodini}), we first prove that we may assume 
that $\bar L_2$ and $\bar L_3$ are adjacent.  If $\bar L_2$ and $\bar L_3$ 
are not adjacent, then $\bar L_1$ is adjacent to both $\bar L_2$ and 
$\bar L_3$ and we may consider ${\rm Sl}_{f_1} (\bar L_1 \cup \bar L_2)$ 
to smooth $\bar L_1 \cup \bar L_2$ to a single irreducible component 
$\bar K$ mapped to a curve with divisor class 
$K:= \bigl( 5 \,;\, 1 , 2 , 2 , 2 , 2 , 2 , 1 , 1 \bigr)$.  Note that 
the divisor class of $K$ is the divisor class of the anticanonical 
divisor on a del Pezzo surface of degree two dominated by $X$.  Thus 
we know that the space $\overline \cM _{bir} \bigl( X,K \bigr)$ is 
irreducible and it contains a point whose image consists of the union 
of the two $(-1)-$curves with divisor classes 
$L_1' := \bigl( 5 \,;\, 1 , 2 , 2 , 2 , 2 , 2 , 1 , 2 \bigr)$ and 
$L_2' := \bigl( 0 \,;\, 0 , 0 , 0 , 0 , 0 , 0 , 0 ,-1 \bigr)$.  Considering 
${\rm Sl}_{f_1} (\bar L_1 \cup \bar L_2)$ we may therefore deform $f_1$ to 
a morphism $f_2 : \bar L_1' \cup \bar L_2' \cup \bar L_3 \rightarrow X$ such 
that the image of $\bar L_i'$ is the $(-1)-$curve $L_i'$.  Thus we have 

$$ \left\{ \begin{array} {r@{\,=\,}l}
L_1' & \bigl( 5 \,;\, 1 , 2 , 2 , 2 , 2 , 2 , 1 , 2 \bigr) \\[5pt]
L_2' & \bigl( 0 \,;\, 0 , 0 , 0 , 0 , 0 , 0 , 0 ,-1 \bigr) \\[5pt]
L_3 & \bigl( 1 \,;\, 1 , 0 , 0 , 0 , 0 , 0 , 1 , 0 \bigr) 
\end{array} \right. 
\hspace{3.5pt} \stackrel {T_{127}} {-\!\!\!-\!\!\!-\!\!\!\longrightarrow } 
\hspace{4.3pt} 
\left\{ \begin{array} {r@{\,=\,}l}
L_1' & \bigl( 6 \,;\, 2 , 3 , 2 , 2 , 2 , 2 , 2 , 2 \bigr) \\[5pt]
L_2' & \bigl( 0 \,;\, 0 , 0 , 0 , 0 , 0 , 0 , 0 ,-1 \bigr) \\[5pt]
L_3 &  \bigl( 0 \,;\, 0 ,-1 , 0 , 0 , 0 , 0 , 0 , 0 \bigr) 
\end{array} \right. $$
which is (up to a permutation) the first case of (\ref{nodini}).

We still need to examine the case in which $\bar L_2$ and $\bar L_3$ are 
adjacent and are given by the second set of equalities in (\ref{nodini}).  
Smoothing the union $\bar L_2 \cup \bar L_3$ to a single irreducible 
component $\bar Q$ we obtain a morphism $f_2$ with dual graph 
$$ \xygraph {[] !~:{@{.}} 
!{<0pt,0pt>;<20pt,0pt>:} 
{\bullet} [rr] {\bullet} 
*\cir<2pt>{}
!{\save +<0pt,8pt>*\txt{$\scriptstyle \bar Q$}  \restore}
- [ll]
*\cir<2pt>{}
!{\save +<0pt,8pt>*\txt{$\scriptstyle \bar L_1$}  \restore} } $$
\vglue0pt {\centerline {Dual graph of $f_2$} \vspace{5pt}}
\noindent
and 
$$ \left\{ \begin{array} {rcl}
(f_2)_* \bar L_1 & = & \bigl( 5 \,;\, 2 , 2 , 2 , 2 , 2 , 2 , 1 , 1 \bigr) \\[5pt]
(f_2)_* \bar Q   & = & \bigl( 1 \,;\, 0 , 0 , 0 , 0 , 0 , 0 , 1 , 0 \bigr)
\end{array} \right. $$

Let $\bar p \in \bar L_1$ be a point such that $f_2(\bar p) \in M$, 
where $M \subset X$ is the $(-1)-$curve with divisor class 
$\ell - e_7 - e_8$.  Considering the morphism 
$$ \xymatrix @C=40pt { {\rm Sl} _{f_2} (\bar Q) 
\ar[r]^{\hspace{10pt} a} & \bar L_1 } $$
we deduce that we may slide $\bar Q$ along $\bar L_1$ until the node between 
these two components reaches the point $\bar p$.  When this happens, the image 
of the limiting position of the image of $\bar Q$ contains a point of $M$.  
Since the intersection product $(f_2)_* [\bar Q] \cdot M$ equals zero, it follows 
that the image of the limiting position of $\bar Q$ must contain $M$.  Thus the 
limit of the morphism $f_2$ under this deformation is a morphism $f_3$ whose 
dual graph is one of the graphs 
$$ \xygraph {[] !~:{@{.}} 
!{<0pt,0pt>;<20pt,0pt>:} 
{\bullet} [rr] {\bullet} [rr] {\bullet} 
*\cir<2pt>{}
!{\save +<0pt,8pt>*\txt{$\scriptstyle \bar E_8$}  \restore}
- [ll]
*\cir<2pt>{}
!{\save +<0pt,8pt>*\txt{$\scriptstyle \bar M$}  \restore}
- [ll]
*\cir<2pt>{}
!{\save +<0pt,8pt>*\txt{$\scriptstyle \bar L_1$}  \restore} } 
\hspace{30pt}
\xygraph {[] !~:{@{.}} 
!{<0pt,0pt>;<20pt,0pt>:} 
{\bullet} [rr] {\bullet} [ur] {\bullet} [dd] {\bullet} 
*\cir<2pt>{}
!{\save +<8pt,0pt>*\txt{$\scriptstyle \bar E_8$}  \restore}
- [ul]
*\cir<2pt>{}
!{\save +<0pt,8pt>*\txt{$\scriptstyle \bar E$}  \restore}
- [ur]
*\cir<2pt>{}
!{\save +<0pt,8pt>*\txt{$\scriptstyle \bar M$}  \restore}
 [dl] - [ll]
*\cir<2pt>{}
!{\save +<0pt,8pt>*\txt{$\scriptstyle \bar L_1$}  \restore} } $$
\vglue0pt {\centerline {Possible dual graphs of $f_3$} \vspace{5pt}}
\noindent
where $\bar E_8$ is mapped to the divisor class $e_8$, and the 
component $\bar E$ is contracted by $f_3$.  The second case happens 
if the $(-1)-$curve with divisor class $e_8$ contains the point 
$f_2(\bar p)$.  In both cases the point represented by $f_3$ lies 
in a unique irreducible component of 
$\overline \cM _{0,0} \bigl( X, -K_X -K_{X'} \bigr)$: in the first 
case thanks to Proposition \ref{grafico}; in the second case thanks 
to Lemma \ref{pizzica} and the fact that the intersection number 
$(f_2)_* [\bar L_1] \cdot E_8$ is one.

We may therefore deform $\bar L_1 \cup \bar E \cup \bar M$ to a 
unique irreducible component $\bar K$ mapped to the divisor $K$ with 
class $\bigl( 6 \,;\, 2 , 2 , 2 , 2 , 2 , 2 , 2 , 2 \bigr) = -2K_X$.

Thus in all cases we found a morphism in the same irreducible 
component of $\overline \cM _{bir} \bigl( X, -K_X -K_{X'} \bigr)$ as 
$f$ whose image contains a nodal integral divisor in $|-2K_X|$.  
Let $E \subset \overline \cM _{0,0} \bigl( X, -K_X -K_{X'} \bigr)$ be 
the subspace consisting of those morphisms containing a component mapped 
birationally to an irreducible divisor in $-2K_X$.  We are going to prove 
that the space $E$ is connected and contained in the smooth locus of 
$\overline \cM _{0,0} \bigl( X, -K_X -K_{X'} \bigr)$.  This concludes the 
proof of the irreducibility of 
$\overline \cM _{bir} \bigl( X, -K_X -K_{X'} \bigr)$.

Any morphism $[f: \bar K \cup \bar E_8 \rightarrow X] \in E$ is determined 
by its image together with one of the points 
$f(\bar K \cap \bar E_8) \in K \cap E_8$.  Since $-2K_X \cdot E_8 = 2$, it 
follows that $E$ has at most two irreducible components.

Suppose $E$ has two irreducible components.  
Consider the morphism $\varphi : X \rightarrow \bP^3$ induced by the 
linear system $|-2K_X|$.  We have already seen that the image is a 
quadric cone $Q$ and that the morphism is ramified along a smooth curve 
$R$ which is the complete intersection of $Q$ with a cubic surface.  The 
$(-1)-$curve $E_8$ has as image the intersection of $Q$ with a plane which 
is everywhere tangent to the curve $R$ (and does not contain the vertex of 
the cone).  Let $e$ be one of the intersection points of $\varphi (E_8)$ 
with $R$.  Projection away from the tangent line $L$ to $R$ at $e$ 
determines a morphism $\pi _L : R \rightarrow \bP^1$ of degree four.  
Since the genus of $R$ is four, it follows from the Hurwitz formula that 
the degree of the ramification divisor of $\pi _L$ is 14.  It is immediate 
to check that $\pi _L$ ramifies above the tangent plane to $Q$ at $e$, and 
that the ramification index is two.  It is also immediate that above the 
plane containing $\varphi (E_8)$ the ramification index is two.  We deduce 
that there are planes in the pencil containing $L$ which are tangent to $R$ 
and are not the tangent plane to $Q$ at $e$ nor the plane containing 
$\varphi (E_8)$.  Such planes correspond to rational divisors $H$ in 
$|-2K_X|$ with the property that $H \cap E_8$ consists of the unique 
point $\varphi ^{-1} (e)$.  Let $\nu : \bar H \rightarrow X$ be a 
birational morphism whose image is one of the divisors $H$ constructed 
above.  The morphism $\nu$ represents a morphism in 
$\overline \cM _{bir} \bigl( X, -2K_X \bigr)$, and since this space is 
irreducible by assumption, we may deform $f$ to a morphism 
$\bar f : \bar H \cup \bar E_8 \longrightarrow X$, such that 
$\bar f|_{\bar H} \simeq \nu$.  Thus $\bar f \in E$ and it clearly lies 
in the intersection of the two irreducible components of $E$.  The space 
$E$ is therefore connected.

Applying Proposition \ref{grafico} we immediately see that all the points 
of $E$ are smooth in $\overline \cM _{0,0} \bigl( X, -K_X -K_{X'} \bigr)$, 
and thus we conclude that the space 
$\overline \cM _{bir} \bigl( X, -K_X -K_{X'} \bigr)$ is irreducible.  \bo

\subsection{The Divisor $-K_{X_8} + Q$}

We prove now a similar result for the Divisor $-K_{X_8} + Q$.

\begin{lem} \label{chicchirichi}
Let $X$ be a del Pezzo surface of degree one and suppose that all the 
rational divisors in $|-K_X|$ are nodal.  Let $Q$ be the divisor class 
of a conic, then the space 
$\overline { \cM } _{bir} \bigl( X, -K_X + Q \bigr)$ is irreducible.
\end{lem}
{\it Proof.}  Let $f: \bP^1 \rightarrow X$ be a free morphism birational to 
its image, such that $f_*[\bP^1] = -K_X + Q$.  As before, we may assume 
that the image of $f$ contains a general point $p$.  Since there is a 
one parameter family of deformations of $f$ whose image contains the 
general point $p$, we may deform $f$ to a morphism 
$f' : \bar C_1 \cup \bar C_2 \longrightarrow X$ such that 
$- K_X \cdot f'_*[\bar C_i] = i$.  Since $p$ is general and contained 
in the image of $f'$ and there are no rational cuspidal divisors in 
$|-K_X|$, it follows that the point represented by $f'$ in 
$\overline { \cM } _{bir} \bigl( X, -K_X + Q \bigr)$ is smooth.

Our next step is to show that we may assume that $f'_*[\bar C_1] = -K_X$.  
Suppose that $f'_*[\bar C_1]$ is a $(-1)-$curve $L \subset X$.  We may 
choose a standard basis $\{ \ell , e_1 , \ldots , e_8 \}$ such that 
$Q = \ell - e_1$ and thus $-K_X + Q = \bigl( 4 \,;\, 2,1,1,1,1,1,1,1 \bigr)$.  
By examining \cite{Ma} Table IV.8, we see that the only possible ways of 
writing $-K_X + Q$ as a sum of a $(-1)-$curve $C_1$ and a nef divisor class 
$C_2$ are (up to permutation of the coordinates 2, \ldots, 8): 
\begin{equation} \label{ancoratre}
 -K_X + Q = 
\left\{ \begin{array} {rcl}
\bigl( 3 \,;\, 2,1,1,1,1,1,1,0 \bigr) &+& \bigl( 1 \,;\, 0,0,0,0,0,0,0,1 \bigr) \\[5pt]
\bigl( 2 \,;\, 1,1,1,1,1,0,0,0 \bigr) &+& \bigl( 2 \,;\, 1,0,0,0,0,1,1,1 \bigr) \\[10pt]
\bigl( 1 \,;\, 0,0,0,0,0,0,1,1 \bigr) &+& \bigl( 3 \,;\, 2,1,1,1,1,1,0,0 \bigr) \\[5pt]
\bigl( 0 \,;\,-1,0,0,0,0,0,0,0 \bigr) &+& \bigl( 4 \,;\, 3,1,1,1,1,1,1,1 \bigr) \\[10pt]
\bigl( 1 \,;\, 1,1,0,0,0,0,0,0 \bigr) &+& \bigl( 3 \,;\, 1,0,1,1,1,1,1,1 \bigr) \\[5pt]
\bigl( 0 \,;\,0,0,-1,0,0,0,0,0 \bigr) &+& \bigl( 4 \,;\, 2,1,2,1,1,1,1,1 \bigr) \\[10pt]
\end{array} \right.
\end{equation}

The automorphisms of ${\rm Pic} (X)$ of the form $T_{1jk}$ preserve the 
conic $Q$, for all $1 < j < k \leq 8$.  We use these automorphisms to 
reduce the number of cases.  We have 
\begin{eqnarray*}
T_{167} \Bigl( 
\bigl( 3 \,;\, 2,1,1,1,1,1,1,0 \bigr) + 
\bigl( 1 \,;\, 0,0,0,0,0,0,0,1 \bigr) \Bigr) & = \\
= \bigl( 2 \,;\, 1,1,1,1,1,0,0,0 \bigr) + 
\bigl( 2 \,;\, 1,0,0,0,0,1,1,1 \bigr) \\[10pt]
T_{178} \Bigl( 
\bigl( 1 \,;\, 0,0,0,0,0,0,1,1 \bigr) + 
\bigl( 3 \,;\, 2,1,1,1,1,1,0,0 \bigr) \Bigr) & = \\
= \bigl( 0 \,;\,-1,0,0,0,0,0,0,0 \bigr) + 
\bigl( 4 \,;\, 3,1,1,1,1,1,1,1 \bigr) \\[10pt]
T_{123} \Bigl( 
\bigl( 1 \,;\, 1,1,0,0,0,0,0,0 \bigr) + 
\bigl( 3 \,;\, 1,0,1,1,1,1,1,1 \bigr) \Bigr) & = \\
= \bigl( 0 \,;\,0,0,-1,0,0,0,0,0 \bigr) + 
\bigl( 4 \,;\, 2,1,2,1,1,1,1,1 \bigr) 
\end{eqnarray*}

We therefore only need to consider the first, third and fifth case in 
list (\ref{ancoratre}).  We reduce the first and third case to the 
fifth one.

If $\bar C_1$ is mapped to the divisor class 
$\bigl( 3 \,;\, 2,1,1,1,1,1,1,0 \bigr)$, then we consider the morphism 
$$ \xymatrix @C=40pt { {\rm Sl} _{f'} (\bar C_2) 
\ar[r]^{\hspace{10pt}a} & \bar C_1 } $$
Since $f'|_{\bar C_2}$ is a free morphism, $a$ is dominant.  Let 
$\bar p \in \bar C_1$ be a point such that $p := f'(\bar p)$ lies on 
the $(-1)-$curve with divisor class $e_1$.  Let 
$g : \bar C_1 \cup \bar E_1 \cup \bar C_2' \longrightarrow X$ be a 
morphism such that $a(g) = \bar p$.  By construction, the dual graph 
of $g$ is 
$$ \xygraph {[] !~:{@{.}} 
!{<0pt,0pt>;<20pt,0pt>:} 
{\bullet} [rr] {\bullet} [rr] {\bullet} 
*\cir<2pt>{}
!{\save +<0pt,8pt>*\txt{$\scriptstyle \bar C_2'$}  \restore}
- [ll]
*\cir<2pt>{}
!{\save +<0pt,8pt>*\txt{$\scriptstyle \bar E_1$}  \restore}
- [ll]
*\cir<2pt>{}
!{\save +<0pt,8pt>*\txt{$\scriptstyle \bar C_1$}  \restore} } $$
\vglue0pt {\centerline {Dual graph of $g$} \vspace{5pt}}
\noindent
where $\bar E_1$ is mapped to the divisor class $e_1$ and 
$\bar C_2'$ to the divisor class $\ell - e_1 - e_8$.  We now smooth 
the union $\bar C_1 \cup \bar E_1$ to a single irreducible component.  
Thus, after a permutation of the indices, we reduced to the fifth case 
in (\ref{ancoratre}).

If $\bar C_1$ is mapped to the divisor class $e_1$, then we proceed 
similarly: break $\bar C_2$ into a component mapping to the divisor 
class $\bigl( 3 \,;\, 2,1,1,1,1,1,1,0 \bigr)$ adjacent to $\bar C_1$, 
and a component mapped to the divisor class $\ell - e_1 - e_8$.  
Smoothing the union of the component $\bar C_1$ with the component 
mapped to $\bigl( 3 \,;\, 2,1,1,1,1,1,1,0 \bigr)$ reduces us to the 
fifth case in (\ref{ancoratre}).

Suppose therefore that the component $\bar C_1$ is mapped to 
$\ell - e_1 - e_2$ and the component $\bar C_2$ is mapped to 
$3 \ell - e_1 - e_3 - \ldots - e_8$.  As above, we may deform the 
morphism $f'$ to a morphism $g$ so that the component $\bar C_2$ 
breaks into a component $\bar E_2$ adjacent to $\bar C_1$ and mapped 
to $e_2$, and into a component $\bar C_2'$ mapped to the divisor 
class $-K_X$.  Smoothing the union $\bar E_2 \cup \bar C_2'$ to a 
single irreducible component, we obtain a morphism 
$g' : \bar C_1 \cup \bar Q \longrightarrow X$, where $\bar C_1$ is 
mapped to $-K_X$ and $\bar Q$ is mapped to $Q$.  Note that we have 
at the moment no control over which rational divisor in the linear 
system $|-K_X|$ the component $\bar C_1$ is mapped to.  Remember that 
with our choice of standard basis we have $Q = \ell - e_1$.  We may 
write $Q = \bigl( \ell - e_1 - e_8 \bigr) + e_8$, and since 
$-K_X \cdot e_8 = 1$, there is a unique point $\bar c$ of $\bar C_1$ 
whose image $c \in X$ lies in $E_8$, the $(-1)-$curve on $X$ with 
divisor class $e_8$.  Considering the dominant morphism 
$$ \xymatrix @C=40pt { {\rm Sl} _{g'} (\bar Q) 
\ar[r]^{\hspace{10pt}a} & \bar C_1 } $$
we may find a morphism $h : \bar C_1 \cup \bar E_8 \cup \bar Q' 
\longrightarrow X$ such that $a(h) = \bar c$.  The dual graph of 
$h$ is 
$$ \xygraph {[] !~:{@{.}} 
!{<0pt,0pt>;<20pt,0pt>:} 
{\bullet} [rr] {\bullet} [rr] {\bullet} 
*\cir<2pt>{}
!{\save +<0pt,8pt>*\txt{$\scriptstyle \bar Q'$}  \restore}
- [ll]
*\cir<2pt>{}
!{\save +<0pt,8pt>*\txt{$\scriptstyle \bar E_8$}  \restore}
- [ll]
*\cir<2pt>{}
!{\save +<0pt,8pt>*\txt{$\scriptstyle \bar C_1$}  \restore} } $$
\vglue0pt {\centerline {Dual graph of $h$} \vspace{5pt}}

Smoothing the components $\bar C_1 \cup \bar E_8$ to a single 
irreducible component $\bar K'$ we obtain a morphism 
$h' : \bar K' \cup \bar Q' \longrightarrow X$ such that 
$Q' := h'_*[\bar Q'] = \ell - e_1 - e_8$ and 
$h'_*[\bar K'] = -K_X + e_8 = -K_{X'}$, where $X'$ is the del Pezzo 
surface obtained from $X$ by contracting the $(-1)-$curve $E_8$.

Let $H \subset \overline { \cM } _{bir} \bigl( X, -K_X + Q \bigr)$ 
be the space of morphisms whose image contains $E_8$ and an integral 
rational divisor in $|-K_{X'}|$.  We have a dominant morphism 
$$ \pi : H \longrightarrow \overline { \cM } _{bir} \bigl( X, -K_{X'} \bigr) $$
whose fibers have length two, since $-K_{X'} \cdot Q' = 2$.  It 
follows that $H$ has at most two irreducible components.  Note that 
the fibers of $\pi $ over the general point of 
$\overline { \cM } _{bir} \bigl( X, -K_{X'} \bigr)$ are smooth points 
of $\overline { \cM } _{bir} \bigl( X, -K_X + Q \bigr)$.  It follows 
that the space $\overline { \cM } _{bir} \bigl( X, -K_X + Q \bigr)$ 
itself has at most two irreducible components, and is irreducible if 
$H$ is.  We prove that if $H$ is reducible, then we can find a smooth 
point of $\overline { \cM } _{bir} \bigl( X, -K_X + Q \bigr)$ lying 
in the intersection of the two components of $H$.  This is enough to 
imply that $\overline { \cM } _{bir} \bigl( X, -K_X + Q \bigr)$ is 
irreducible.

Suppose thus that $H$ as two irreducible components.  Let 
$\varphi : X \rightarrow \bP^2$ be the morphism induced by the linear 
system $-K_{X'}$.  The morphism $\varphi $ is the contraction of $E_8$ 
to $X'$ followed by the anticanonical double cover of $\bP^2$ ramified 
above a smooth plane quartic $R$.  The image of $Q'$ in $\bP^2$ is a 
conic $\tilde Q$ containing the image of $E_8$ and everywhere tangent 
to $R$.  The image of $K'$ is a tangent line to $R$.  To conclude it is 
enough to find a line in $\bP^2$ which is tangent to both $R$ and 
$\tilde Q$ at a point not on $R$.  The dual curve of $R$ is a plane 
curve of degree twelve and the dual curve of $\tilde Q$ is a plane conic.  
Thus they meet along a scheme of length 24 and they are tangent at the 
points corresponding to the points where $R$ and $\tilde Q$ are tangent.  
Since there are four such points, it follows that we may find a line 
which is tangent to $R$ and $\tilde Q$ at distinct points.  Such a line 
corresponds to a point in the intersection of the two components of $H$.  
Using Proposition \ref{grafico} it is easy to check that this point is 
smooth in $\overline { \cM } _{bir} \bigl( X, -K_X + Q \bigr)$.  This 
completes the proof of the lemma.  \bo

\section{Conclusion}

\subsection{The Irreducibility of 
$\overline { \cM } _{bir} \bigl( X_\delta , \beta \bigr)$}

We are now ready to prove the main theorems of this article.

\begin{thm} \label{nonuno}
Let $X_\delta $ be a del Pezzo surface of degree $9-\delta \geq 2$.  
The spaces $\overline { \cM } _{bir} \bigl( X_\delta , \beta \bigr)$ 
are irreducible or empty for every divisor 
$\beta \in {\rm Pic} (X_\delta )$.
\end{thm}
{\it Proof.}  Suppose $\overline { \cM } _{bir} \bigl( X_\delta , \beta \bigr)$ 
is not empty.  Then $\beta $ is represented by an effective integral 
curve on $X_\delta $.

If $\beta $ is not nef, then it follows that 
$\beta ^2 < 0$.  We deduce that $\beta $ is a positive multiple $d$ of a 
$(-1)-$curve.  If $d = 1$, then 
$\overline { \cM } _{bir} \bigl( X_\delta , \beta \bigr)$ consists of a 
single point.  If $d > 1$, then the space 
$\overline { \cM } _{bir} \bigl( X_\delta , \beta \bigr)$ is empty.  In 
this case, the space 
$\overline { \cM } _{0,0} \bigl( X_\delta , \beta \bigr)$ is irreducible, 
since it is dominated by the space of triples of homogeneous polynomials 
of degree $d$ in two variables.

Suppose now that $\beta $ is a nef divisor.  Thanks to Theorem 
\ref{passo}, we simply need to check that on a del Pezzo surface of 
degree at least two, the spaces 
$\overline { \cM } _{bir} \bigl( X_\delta , \beta \bigr)$ are 
irreducible for all effective integral divisor classes $\beta $ such 
that $-K_{X_\delta } \cdot \beta $ equals two or three.  The divisors 
of degree two on $X_\delta $ are the conics and, if $\delta = 7$, the 
divisor $-K_{X_7}$.  If $\beta $ is a conic, then 
$\overline { \cM } _{bir} \bigl( X_\delta , \beta \bigr)$ is isomorphic 
to $\bP^1$.  If $\beta = -K_{X_7}$, then 
$\overline { \cM } _{bir} \bigl( X_7 , -K_{X_7} \bigr)$ is 
isomorphic to a smooth plane quartic curve, Proposition \ref{cadute}.

The nef divisors of degree three on $X_\delta $ are $-K_{X_6}$ and $\ell $, 
where $X_6$ is a del Pezzo surface of degree three dominated by $X_\delta $ 
and $\ell $ is part of a standard basis $\{ \ell , e_1 , \ldots , e_\delta \}$.  
The first case is treated in Proposition \ref{cane}, the second case is 
treated in Theorem \ref{barbapapa}.  This concludes the proof of the 
theorem.  \bo

\begin{thm} \label{nonplusuno}
Let $X_8$ be a general del Pezzo surface of degree one.  The spaces 
$\overline { \cM } _{bir} \bigl( X_8 , \beta \bigr)$ are irreducible 
or empty for every divisor $\beta \in {\rm Pic} (X_\delta )$, with 
the unique exception of $\beta = -K_{X_8}$.  The space 
$\overline { \cM } _{bir} \bigl( X_8 , -K_{X_8} \bigr)$ is a reduced 
scheme of length twelve.
\end{thm}
{\it Proof.}  Proceeding as before, we only need to prove the irreducibility 
of $\overline { \cM } _{bir} \bigl( X_8 , \beta \bigr)$ for the nef 
divisors of anticanonical degree two or three.  The nef divisor classes 
on $X_8$ which are not ample, are the pull-back of nef divisor classes 
from del Pezzo surfaces of larger degree.  Thus we only need to consider 
ample divisor classes of anticanonical degree two or three.

The only ample divisor of degree two is $-2K_{X_8}$ and the mapping space 
$\overline { \cM } _{bir} \bigl( X_8 , -2K_{X_8} \bigr)$ is irreducible 
thanks to Theorem \ref{cabala}.

The ample divisor classes of degree three on $X_8$ are $-3K_{X_8}$, 
$-K_{X_8} - K_{X_7}$ and $-K_{X_8} + Q$, where $X_7$ is a del Pezzo 
surface of degree two dominated by $X_8$ and $Q$ is the divisor class 
of a conic.  
The space $\overline { \cM } _{bir} \bigl( X_8 , -3K_{X_8} \bigr)$ is 
irreducible thanks to Theorem \ref{maschera}.  
The space $\overline { \cM } _{bir} \bigl( X_8 , -K_{X_8} - K_{X_7} \bigr)$ 
is irreducible thanks to Lemma \ref{coccode}.  
The space $\overline { \cM } _{bir} \bigl( X_8 , -K_{X_8} + Q \bigr)$ is 
irreducible thanks to Lemma \ref{chicchirichi}.  This concludes the 
proof of the theorem.  \bo

\noindent
{\it Remark}.  The genericity assumption on $X_8$ in the statement 
of the Theorem \ref{nonuno} can be made more explicit.  Our argument 
requires the surface $X_8$ to have only nodal rational divisors in 
$|-K_{X_8}|$ and the space 
$\overline { \cM } _{bir} \bigl( X_8 , -2K_{X_8} \bigr)$ to be 
irreducible.  This last condition in turn is certainly satisfied (cf. 
Theorem \ref{cabala} and its proof) if the ramification curve 
$R \subset \bP^3$ of the morphism $\varphi : X_8 \rightarrow \bP^3$ 
induced by $-2K_{X_8}$ does not admit planes $P \subset \bP^3$ 
transverse to the image of $\varphi$ and intersecting $R$ along a 
divisor of the form $3 \bigl( (p) + (q) \bigr)$.

As a corollary of the above Theorems, we deduce the irreducibility of 
the Severi varieties of rational curves on the del Pezzo surfaces.  Let 
$\beta $ be a divisor class in ${\rm Pic} (X_\delta )$ and let 
$V_{0,\beta } \subset |\beta |$ be the closure of the set of points 
corresponding to integral rational divisors.  We call $V_{0,\beta }$ 
the Severi variety of rational curves on $X$ with divisor class $\beta $.

\begin{cor}
Let $X_\delta $ be a del Pezzo surface of degree $9-\delta \geq 2$.  
The Severi varieties $V_{0,\beta }$ of rational curves on $X_\delta $ 
are either empty of irreducible for every divisor 
$\beta \in {\rm Pic} (X_\delta )$.  \bo
\end{cor}

\begin{cor}
Let $X_8$ be a general del Pezzo surface of degree one.  The Severi 
varieties $V_{0,\beta }$ of rational curves on $X_8$ are either 
empty or irreducible for every divisor $\beta \in {\rm Pic} (X_8)$, 
with the unique exception of $\beta = -K_{X_8}$.  \bo
\end{cor}

\end{document}